\documentclass{article}

\usepackage{arxiv}

\usepackage[utf8]{inputenc} 
\usepackage[T1]{fontenc}    
\usepackage{hyperref}       
\usepackage{url}            
\usepackage{booktabs}       
\usepackage{amsfonts}       
\usepackage{nicefrac}       
\usepackage{microtype}      
\usepackage{lscape}
\usepackage[authoryear, sort]{natbib}

\usepackage{morefloats}
\usepackage{color}
\usepackage{multirow}
\usepackage{latexsym}
\usepackage{amsmath,amssymb,amsthm,graphicx}
\usepackage{dsfont}
\usepackage{enumitem}
\usepackage{mathtools}

\allowdisplaybreaks
\newtheoremstyle{thm}
{9pt}
{9pt}
{\itshape}
{}
{\bfseries}
{.}
{ }
{}
\theoremstyle{thm}

\newtheorem{theorem}{Theorem}[section]
\newtheorem{lemma}[theorem]{Lemma}
\newtheorem{corollary}[theorem]{Corollary}
\newtheorem{prop}[theorem]{Proposition}

\newtheoremstyle{def}
{9pt}
{9pt}
{}
{}
{\bfseries}
{.}
{ }
{}
\theoremstyle{def}

\newtheorem{definition}[theorem]{Definition}
\newtheorem{remark}[theorem]{Remark}
\newtheorem{example}[theorem]{Example}

\newenvironment{prf}{\noindent\textbf{\emph{Proof.}}}{\qed}

\newcommand{\R}{\mathbb{R}} 
\newcommand{\N}{\mathbb{N}} 
\newcommand{\X}{\mathbb{X}} 
\newcommand{\NN}{\mathbf{N}} 
\newcommand{\E}{\mathbb{E}} 
\newcommand{\PP}{\mathbb{P}} 
\newcommand{\dd}{\mathrm{d}}

\renewcommand{\footnoterule}{%
	\kern -3.5pt
	\hrule width \textwidth height 1pt
	\kern 3.5pt
}

\makeatletter
\def\blfootnote{\xdef\@thefnmark{}\@footnotetext}
\makeatother

\title{Structural properties of Gibbsian point processes in abstract spaces}

\author{Steffen Betsch\\
Institute of Stochastics,\\ Karlsruhe Institute of Technology (KIT),\\ Germany.\\
\href{mailto:steffen.betsch@kit.edu}{steffen.betsch@kit.edu}\\
}

\begin{document}

\date{\today}
\maketitle


\begin{abstract}
	In the language of random counting measures many structural properties of the Poisson process can be studied in arbitrary measurable spaces. We provide a similarly general treatise of Gibbs processes. With the GNZ equations as a definition of these objects, Gibbs processes can be introduced in abstract spaces without any topological structure. In this general setting, partition functions, Janossy densities, and correlation functions are studied. While the definition covers finite and infinite Gibbs processes alike, the finite case allows, even in abstract spaces, for an equivalent and more explicit characterization via a familiar series expansion. Recent generalizations of factorial measures to arbitrary measurable spaces, where counting measures cannot be written as sums of Dirac measures, likewise allow to generalize the concept of Hamiltonians. The DLR equations, which completely characterize a Gibbs process, as well as basic results for the local convergence topology are also formulated in full generality. We prove a new theorem on the extraction of locally convergent subsequences from a sequence of point processes and use this statement to provide existence results for Gibbs processes in general spaces with potentially infinite range of interaction. These results are used to guarantee the existence of Gibbs processes with cluster-dependent interactions and to prove a recent conjecture concerning the existence of Gibbsian particle processes. 
\end{abstract}

\section{Introduction}
\label{SEC Introduction}

Let $(\X, \mathcal{X})$ be a measurable space which is \textit{localized}, meaning there exist sets $B_1 \subset B_2 \subset \dotso$ from $\mathcal{X}$ such that $\bigcup_{j = 1}^\infty B_j = \X$. Let $\mathcal{X}_b$ be the collection of all sets in $\mathcal{X}$ which are contained in one of the $B_j$ and call those sets \textit{bounded}. A measure $\lambda$ on $\X$ is \textit{locally finite} if $\lambda(B) < \infty$ for all $B \in \mathcal{X}_b$. Denote by $\NN = \NN(\X)$ the set of all measures $\mu$ on $\X$ for which $\mu(B) \in \N_0$ for each $B \in \mathcal{X}_b$, and endow $\NN$ with the $\sigma$-field $\mathcal{N}$ generated by the maps $\pi_B : \mu \mapsto \mu(B)$, $B \in \mathcal{X}$. A \textit{point process} in $\X$ is a random element in $\NN$ defined on some underlying probability space $(\Omega, \mathcal{F}, \PP)$.

Let $\kappa : \X \times \NN \to [0, \infty)$ be measurable and fix a locally finite measure $\lambda$ on $\X$. A point process $\eta$ in $\X$ is called \textit{Gibbs process} with \textit{Papangelou (conditional) intensity} (PI) $\kappa$ and reference measure $\lambda$ if
\begin{equation} \label{GNZ equation}
\E \bigg[ \int_\X f(x, \eta) \, \dd \eta(x) \bigg]
= \E \bigg[ \int_\X f(x, \eta + \delta_x) \, \kappa(x, \eta) \, \dd \lambda(x) \bigg]
\end{equation}
for all measurable maps $f : \X \times \NN \to [0, \infty]$. Here $\delta_x$ denotes the Dirac measure in $x \in \X$. These defining equations are the \textit{GNZ equations} named after \cite{G:1976} and  \cite{NZ:1979}. The definition includes both finite and infinite processes and it does not require any underlying topological structure on $\X$.

We prove in Section \ref{SEC Definition} that, in just about any reasonable state space, a Gibbs process can only exist if $\kappa$ satisfies
\begin{equation} \label{cocycle assumption}
	\kappa(x, \mu) \cdot \kappa(y, \mu + \delta_x) = \kappa(y, \mu) \cdot \kappa(x, \mu + \delta_y)
\end{equation}
for all $\mu \in \NN$ and $x, y \in \X$, a property which is called the \textit{cocycle relation}. For $m \in \N$ we define $\kappa_m : \X^m \times \NN \to [0, \infty)$ by
\begin{equation*}
	\kappa_m(x_1, \dots, x_m, \mu)
	= \kappa(x_m, \mu) \cdot \kappa(x_{m - 1}, \mu + \delta_{x_m}) \cdot \dotso \cdot \kappa(x_1, \mu + \delta_{x_2} + \dotso + \delta_{x_{m}}).
\end{equation*}
By the cocycle property, these functions are symmetric in their first $m$ components.

One of the main objects of interest in the context of Gibbs processes is the \textit{partition function}. It is immediately clear that its definition can be given on arbitrary state spaces. Indeed, for $B \in \mathcal{X}$ we define $Z_B : \NN \to [0, \infty]$ as
\begin{equation} \label{definition of partition function}
	Z_B(\psi)
	= 1 + \sum_{m = 1}^\infty \frac{1}{m!} \int_{B^m} \kappa_m(x_1, \dots, x_m, \psi) \, \dd \lambda^m(x_1, \dots, x_m) .
\end{equation}
The function $Z_B$ is measurable and satisfies $Z_B(\psi) \geq 1$ for every $\psi \in \NN$. The partition function can be used to give a more explicit definition of finite Gibbs processes which is equivalent to \eqref{GNZ equation}. Usually these finite Gibbs processes are considered on bounded sets, but we extend the conventional knowledge ever so slightly by showing that boundedness (or finite $\lambda$-measure) of the domain of the process is inessential, what matters is that the partition function is finite. More precisely, if $C \in \mathcal{X}$ and $\psi \in \NN$, then a finite Gibbs process in $C$ with \textit{boundary condition} $\psi$ exists if, and only if, $Z_C(\psi) < \infty$. If $\xi$ is such a process, then the distribution $\PP^\xi$ of $\xi$ is given through
\begin{equation*}
	\PP^\xi(\cdot)
	= \frac{1}{Z_C(\psi)} \bigg( \mathds{1}\{ \mathbf{0} \in \cdot \, \} + \sum_{m = 1}^\infty \frac{1}{m!} \int_{C^m} \mathds{1}\Big\{ \sum_{j = 1}^m \delta_{x_j} \in \cdot \, \Big\} \, \kappa_m(x_1, \dots, x_m, \psi) \, \dd\lambda^m(x_1, \dots, x_m) \bigg) ,
\end{equation*}
where $\mathbf{0} \in \NN$ denotes the zero measure.

A major difficulty in dealing with point processes, and Gibbs processes in particular, in arbitrary measurable spaces is that a counting measure $\mu \in \NN$ cannot (in general) be written as a sum of Dirac measures, that is, more intuitively, the atoms of $\mu$ cannot be identified. It is thus not at all clear how energy functions can be defined, one of the most basic notions that underlie the theory of Gibbs processes. Recent constructions of factorial measures in general spaces by \cite{LP:2017} allow for such a definition. In particular, they show how to construct, for each $\mu \in \NN$ and $m \in \N$, a measure $\mu^{(m)}$ in $\NN(\X^m)$ that behaves essentially like the factorial measure $\sum_{j_1, \dots, j_m \leq k}^{\neq} \delta_{(x_{j_1}, \dots, x_{j_m})}$ which is defined for a sum of Dirac measures $\sum_{j = 1}^k \delta_{x_j} \in \NN$. Technicalities aside, these concepts allow us to define the \textit{Hamiltonian} $H : \NN \times \NN \to (- \infty, \infty]$ as
\begin{align*}
	H(\mu, \psi)
	= \infty \cdot \mathds{1}\big\{ \mu(\X) = \infty \big\} - \sum_{m = 1}^\infty &\mathds{1}\big\{ \mu(\X) = m \big\} \cdot \log\bigg( \frac{1}{m!} \int_{\X^m} \kappa_m(x_1, \dots, x_m, \psi) \, \dd \mu^{(m)}(x_1, \dots, x_m) \bigg)  .
\end{align*}
Full detail is given in Section \ref{SEC Hamiltonians}, but let us observe that for $\mu = \sum_{j = 1}^m \delta_{x_j}$ we have
\begin{equation} \label{hamiltonian for Dirac sum}
	\mathrm{e}^{- H(\mu, \psi)} = \kappa_m(x_1, \dots, x_m, \psi) .
\end{equation}
With this fundamental object available on arbitrary measurable spaces, the well-known DLR equations, and the characterization of the Gibbs processes they provide, generalize to this abstraction. For the readers convenience and to provide a reference for others, details on the DLR equations are added in Section \ref{SEC DLR equations}.

In Section \ref{SEC prelinimaries} we collect properties of the local convergence topology introduced by \cite{GZ:1993}, a mode of convergence that is, in alignment with our abstract setting, not bound to any topological structures on the space $\X$ itself. With Theorem \ref{THM existence of locally convergent subsequences} we provide a new construction to extract from a sequence of point processes a locally convergent subsequence as well as a corresponding limit process. For the theorem to hold, we essentially need the uniform integrability of the Janossy densities associated to the point processes. 
In Section \ref{SEC existence result Gibbs processes} this result is used to derive a general existence result for Gibbs processes which is then applied to more explicit models. For pair potentials we obtain very minor improvements over the existing literature, but the abstraction also paves the way for two new classes of results. For one thing, we provide existence results for cluster-dependent interactions as they are considered by \cite{LO:2021}. In order to rigorously construct the disagreement coupling, \cite{LO:2021} consider on a Borel space $(\X, \mathcal{X})$ a symmetric and binary relation and assume that the PI $\kappa(x, \mu)$ depends only on the cluster that $x$ forms with the points in $\mu$ via the relation. In this context, a usual prerequisite is a suitable kind of subcriticality, meaning that the clusters of Poisson process points via the given relation are finite almost surely. Under this assumption, we prove in Section \ref{SEC Gibbs point processes with cluster-dependent interaction} that a unique Gibbs process does exists. Note that the interaction between the points of such a process depends purely on the corresponding clusters and that this interaction need not have a finite range. The uniqueness result we provide in this setting covers the results by \cite{HtH:2019} and \cite{BHLV:2020}, who consider particle processes with the binary relation on the state space given through the intersection of particles, but is much more general. Also we ensure existence of the processes from \cite{HtH:2019} within the region of uniqueness. Moreover, our approach leads to manageable conditions for the existence (and uniqueness) of Gibbsian particle processes, detailed in Section \ref{SEC existence and uniqueness of Gibbs particle processes}. Though the existence result is substantially more general, it particularly covers the conjecture by \cite{BHLV:2020} who study a special class of such particle processes and emphasize that their existence is not guaranteed by the available literature. 

To prove existence of Gibbs processes in $\R^d$, \cite{D:1969} uses a different compactness criterion. Moreover, \cite{D:2009}, \cite{DDG:2012}, \cite{DV:2020}, and \cite{RZ:2020} use that \cite{GZ:1993} have established the compactness of level sets of the specific entropy in the local convergence topology and provide existence results for stationary Gibbs processes in $\R^d$. As we elaborate in Section \ref{SEC existence result Gibbs processes}, the latter approach leads to better existence results for Gibbs processes, but only in $\R^d$ and not in the abstract setting handled here. There also exists an analytical approach to the theory of Gibbs processes as laid out by \cite{KPR:2012} and \cite{CDKP:2018}, but the focus is on pair potentials in $\R^d$ and the theory is quite technical. Another class of existence results comes from the description of Gibbs measures through specifications. The corresponding results, like ours, are available in very abstract spaces, cf. \cite{P:1976}. Also they use some variant of the Dunford-Pettis property \citep[see Chapter 5 of][]{P:2005:specif}, just like we do. A main difference between this last class of existence results and ours lies in the general approach toward Gibbs process, namely between introducing the Gibbs process via the DLR equations in terms of specifications or considering Gibbs processes as solutions to the GNZ equations in the modern point process theoretic framework that we introduce. The latter approach leads to a transparent theory that is better compatible with the point process theoretical foundations of stochastic geometry and which fits neatly into the interpretation of Gibbs processes in spatial statistics.

For the proof of uniqueness in distribution of Gibbs processes there exists an even more diverse pool of methods. One classical method is due to \cite{D:1968:3} and used, for instance, by \cite{GH:1996}, \cite{PZ:1999}, \cite{BP:2002}, and \cite{CDKP:2018}, and very recently also by \cite{HZ:2021}. Another classical approach based on the Kirkwood-Salsburg equations, which appears in the context of cluster expansions, goes back to \cite{R:1969}. Apart from its novel contributions, the paper by \cite{J:2019} provides an overview of the latter method and the corresponding literature. The recent papers by \cite{MP:2020} and \cite{Z:2021} also use this approach. A third approach is based on suitable Markovian dynamics, see \cite{SY:2013}, and a forth method, which also appears in the paper at hand, is the so-called disagreement percolation which goes back to \cite{vdBM:1994} and was further developed by \cite{Ht:2019}, \cite{HtH:2019}, and \cite{LO:2021}. An adaptation of this method in combination with the random connection model was recently used by \cite{BL:2021}. To classify the uniqueness result which is provided in the article at hand, let us first mention that most of the above papers \citep[and also the preprint by][]{MP:2021} focus on pair potentials. An exception are the results obtained via disagreement percolation, which is used here as well. Following the general frame of the present paper, we provide a rather abstract result and also deliver the corresponding existence results which were previously missing in the context of disagreement percolation.

As regards the uniqueness question, it has to be mentioned that one of the major goals in the theory of (continuous) Gibbs measures and statistical mechanics is the proof and understanding of a phase transition from uniqueness to non-uniqueness. Even for rather restricted classes of interactions, like repulsive pair interactions with finite range, no general result is proven, though such a result is widely believed to hold. Rigorous results for phase transition focus on explicit special cases like the Potts model, cf. \cite{GH:1996} and \cite{GLL:2005}, a Kac type potential, cf. \cite{LMP:1999}, or the Widom-Rowlinson model, cf. \cite{WR:1970}, \cite{R:1971}, and \cite{DH:2019, DH:2021}.

\section{Basic properties}
\label{SEC Definition}

In this section, as well as in Sections \ref{SEC Jan., Cor. and Palm meas. of Gibbs processes} -- \ref{SEC DLR equations}, we assume the most general setting from the beginning of Section \ref{SEC Introduction}. Note that the terms in equation \eqref{GNZ equation} and in the following lemma are well defined by Lemma \ref{LEMMA measurability of factorial-measure-integrals}.
\begin{lemma}[\bf Multivariate GNZ equations] \label{LEMMA mult. GNZ}
	Let $\eta$ be a Gibbs process with PI $\kappa$ and $m \in \N$. Then, for any measurable function $f : \X^m \times \NN \to [0, \infty]$,
	\begin{align*}
		&\E \bigg[ \int_{\X^m} f(x_1, \dots, x_m, \eta) \, \dd \eta^{(m)}(x_1, \dots, x_m) \bigg] \\
		&\quad = \E \bigg[ \int_{\X^m} f(x_1, \dots, x_m, \eta + \delta_{x_1} + \dotso + \delta_{x_m}) \, \kappa_m(x_1, \dots, x_m, \eta) \, \dd \lambda^m(x_1, \dots, x_m) \bigg] .
	\end{align*}
\end{lemma}
As for the notation in Lemma \ref{LEMMA mult. GNZ}, we write $\nu^{(m)}$ for the $m$-th factorial measure of $\nu \in \NN$, the existence of which is guaranteed in this full generality by Proposition 4.3 of \cite{LP:2017}, see also Appendix \ref{Appendix factorial measures}. Throughout this work, we adopt the usual $\infty$-conventions, that is, $\infty + c = c + \infty = \infty$ for all $c \in (- \infty, \infty]$, and $\infty \cdot c = c \cdot \infty = \infty$ for all $c \in (0, \infty]$, as well as $\infty \cdot c = c \cdot \infty = - \infty$ for all $c \in [-\infty, 0)$. We also use the standard measure theory convention $\infty \cdot 0 = 0 \cdot \infty = 0$, and put $e^{- \infty} = 0$ as well as $\log(0) = -\infty$. Other than these, we do not define any terms involving $+\infty$ or $-\infty$.

\begin{prf}
	We prove the assertion by induction on $m$. For $m = 1$, the claim corresponds to the GNZ equations. We proceed to prove that if the claim is true for some fixed $m$, it also holds for $m + 1$. Let $D \in \mathcal{X}^{\otimes (m + 1)}$ and $A \in \mathcal{N}$. We define the measurable map $F : \X^m \times \NN \to (- \infty, \infty]$,
	\begin{equation*}
		F(x_1, \dots, x_m, \mu)
		= \mathds{1}_A(\mu) \cdot \bigg( \int_\X \mathds{1}_D(x_1, \dots, x_m, y) \, \dd\mu(y) - \sum_{j = 1}^m \mathds{1}_D(x_1, \dots, x_m, x_j) \bigg) .
	\end{equation*}
	It follows from the proof of Proposition A.18 of \cite{LP:2017} that $F^-(x_1, \dots, x_m, \mu) = 0$ for $\mu^{(m)}$-almost every (a.e.) $(x_1, \dots, x_m) \in \X^m$ and all $\mu \in \NN$, where $F^-$ denotes the negative part of $F$ (and similarly $F^+$ is the positive part of $F$). By the characterizing equation \eqref{def. equation factorial measures} of the factorial measure, we obtain
	\begin{equation*}
		\E\bigg[ \int_{\X^{m + 1}} \mathds{1}_D(x_1, \dots, x_{m + 1}) \, \mathds{1}_A(\eta) \, \dd\eta^{(m + 1)}(x_1, \dots, x_{m + 1}) \bigg] 
		= \E\bigg[ \int_{\X^m} F^+(x_1, \dots, x_m, \eta) \, \dd\eta^{(m)}(x_1, \dots, x_m) \bigg] .
	\end{equation*}
	By the induction hypothesis, the multivariate GNZ equations hold for every measurable map $\X^m \times \NN \to [0, \infty]$, so we can apply it to $F^+$. Thus, the previous term equals
	\begin{align*}
		&\E\bigg[ \int_{\X^m} F^+\Big( x_1, \dots, x_m, \eta + \sum_{i = 1}^m \delta_{x_i} \Big) \, \kappa_m(x_1, \dots, x_m, \eta) \, \dd\lambda^m(x_1, \dots, x_m) \bigg] \\
		&\quad = \E\bigg[ \int_{\X^m} \mathds{1}_A\Big( \eta + \sum_{i = 1}^m \delta_{x_i} \Big) \bigg( \int_\X \mathds{1}_D(x_1, \dots, x_{m + 1}) \, \dd\eta(x_{m + 1}) \bigg) \, \kappa_m(x_1, \dots, x_m, \eta) \, \dd\lambda^m(x_1, \dots, x_m) \bigg] .
	\end{align*}
	With the GNZ equations and the definition of $\kappa_{m + 1}$ we further calculate this term as
	\begin{align*}	
		&\E\bigg[ \int_{\X^m} \int_\X \mathds{1}_D(x_1, \dots, x_{m + 1}) \, \mathds{1}_A\Big( \eta + \sum_{i = 1}^{m + 1} \delta_{x_i} \Big) \, \kappa_m(x_1, \dots, x_m, \eta + \delta_{x_{m + 1}}) \, \kappa(x_{m + 1}, \eta) \, \dd\lambda(x_{m + 1}) \, \dd\lambda^m(x_1, \dots, x_m) \bigg] \\
		&\quad = \E\bigg[ \int_{\X^{m + 1}} \mathds{1}_D(x_1, \dots, x_{m + 1}) \, \mathds{1}_A\Big( \eta + \sum_{i = 1}^{m + 1} \delta_{x_i} \Big) \, \kappa_{m + 1}(x_1, \dots, x_{m + 1}, \eta) \, \dd\lambda^{m+1}(x_1, \dots, x_{m + 1}) \bigg] .
	\end{align*}
	Define on $\X^{m + 1} \times \NN$ the measures
	\begin{equation*}
		\mathrm{C}(F)
		= \E\bigg[ \int_{\X^{m + 1}} \mathds{1}_F(x_1, \dots, x_{m + 1}, \eta) \, \dd\eta^{(m + 1)}(x_1, \dots, x_{m + 1}) \bigg]
	\end{equation*}
	and
	\begin{equation*}
		\tilde{\mathrm{C}}(F)
		= \E\bigg[ \int_{\X^{m + 1}} \mathds{1}_F\Big( x_1, \dots, x_{m + 1}, \eta + \sum_{i = 1}^{m + 1} \delta_{x_i} \Big) \, \kappa_{m + 1}(x_1, \dots, x_{m + 1}, \eta) \, \dd\lambda^{m + 1}(x_1, \dots, x_{m + 1}) \bigg] ,
	\end{equation*}
	for $F \in \mathcal{X}^{\otimes (m + 1)} \otimes \mathcal{N}$. By the calculation above, these measures are equal on the $\pi$-system $\{ D \times A : D \in \mathcal{X}^{\otimes (m + 1)}, \, A \in \mathcal{N} \}$ that generates $\mathcal{X}^{\otimes (m + 1)} \otimes \mathcal{N}$. Moreover, the measures are $\sigma$-finite as
	\begin{equation*}
		\tilde{\mathrm{C}}\big( B_\ell^{m + 1} \times \big\{ \mu \in \NN : \mu(B_\ell) \leq n \big\} \big)
		= \mathrm{C}\big( B_\ell^{m + 1} \times \big\{ \mu \in \NN : \mu(B_\ell) \leq n \big\} \big)
		= \E\Big[ \mathds{1}\big\{ \eta(B_\ell) \leq n \big\} \cdot \eta^{(m + 1)}(B_\ell^{m + 1}) \Big]
		\leq n!
	\end{equation*}
	for all $\ell, n \in \N$. Thus, the uniqueness theorem for measures yields $\mathrm{C} = \tilde{\mathrm{C}}$. Monotone approximation with simple functions implies the claim for general measurable functions $f : \X^{m + 1} \times \NN \to [0, \infty]$, and the induction is complete.
\end{prf}

\vspace{3mm}

Notice that if the underlying measurable space has a measurable diagonal, that is, $\{ (x, x) : x \in \X \} \in \mathcal{X}^{\otimes 2}$, then the map $\X \times \NN \ni (x, \mu) \mapsto \mu \setminus \delta_x \in \NN$ is measurable, where
\begin{equation*}
	\mu \setminus \delta_x
	= \mu - \delta_x \, \mathds{1}\big\{ \mu(\{ x \}) > 0 \big\} .
\end{equation*}
This fact is discussed in Appendix \ref{Appendix measurable diagonals}. In particular, if $(\X, \mathcal{X})$ is separable (i.e. $\{ x \} \in \mathcal{X}$ for each $x \in \X$) and (the $\sigma$-field $\mathcal{X}$ is) countably generated then it has a measurable diagonal. Any Borel space, in the sense that is used by \cite{LP:2017} or \cite{K:2017}, has a measurable diagonal.

\begin{corollary} \label{COR ad mult. GNZ}
	Let $(\X, \mathcal{X})$ have a measurable diagonal. Further, let $\eta$ be a Gibbs process with PI $\kappa$, and fix $m \in \N$. Then, for any measurable function $f : \X^m \times \NN \to [0, \infty]$,
	\begin{align*}
	&\E \bigg[ \int_{\X^m} f(x_1, \ldots, x_m, \eta \setminus \delta_{x_1} \setminus \ldots \setminus \delta_{x_m}) \, \dd \eta^{(m)}(x_1, \ldots, x_m) \bigg] \\
	&\quad = \E \bigg[ \int_{\X^m} f(x_1, \ldots, x_m, \eta) \, \kappa_m(x_1, \ldots, x_m, \eta) \, \dd \lambda^m(x_1, \ldots, x_m) \bigg] .
	\end{align*}
\end{corollary}
The proof is immediate upon applying Lemma \ref{LEMMA mult. GNZ} to the map
\begin{equation*}
	(x_1, \ldots, x_m, \mu) 
	\mapsto f(x_1, \ldots, x_m, \mu \setminus \delta_{x_1} \setminus \ldots \setminus \delta_{x_m}) ,
\end{equation*}
which is measurable by Lemma \ref{LEMMA Appen. meas. of reducing a counting measure}, and using that $(\mu + \delta_{x_1} + \ldots + \delta_{x_m}) \setminus \delta_{x_1} \setminus \ldots \setminus \delta_{x_m} = \mu$ for $\mu \in \NN$ and $x_1, \ldots, x_m \in \X$.

We apply Corollary \ref{COR ad mult. GNZ} to provide a first observation about $\kappa_m$.
\begin{lemma} \label{LEMMA symmetry of kappa_m}
	Assume that $(\X, \mathcal{X})$ has a measurable diagonal. Let $\eta$ be a Gibbs process with PI $\kappa$, and fix $m \in \N$. Then, for any permutation $\tau$ of $\{1, \dots, m\}$,
	\begin{equation*}
		\kappa_m(x_1, \dots, x_m, \mu)
		= \kappa_m\big( x_{\tau(1)}, \dots, x_{\tau(m)}, \mu \big)
	\end{equation*}
	for $(\lambda^m \otimes \PP^\eta)$-a.e.\ $(x_1, \dots, x_m, \mu) \in \X^m \times \NN$.
\end{lemma}
\begin{prf}
	Let $f : \X^m \times \NN \to [0, \infty]$ be measurable. Then, by Corollary \ref{COR ad mult. GNZ}, we have
	\begin{align*}
		&\int_\NN \int_{\X^m} f(x_1, \ldots, x_m, \mu) \, \kappa_m(x_{\tau(1)}, \ldots, x_{\tau(m)}, \mu) \, \dd \lambda^m(x_1, \ldots, x_m) \, \dd \PP^\eta(\mu) \\
		&\quad = \E \bigg[ \int_{\X^m} f\big( x_{\tau^{-1}(1)}, \ldots, x_{\tau^{-1}(m)}, \eta \big) \, \kappa_m(x_1, \ldots, x_m, \eta) \, \dd \lambda^m(x_1, \ldots, x_m) \bigg] \\
		&\quad = \E \bigg[ \int_{\X^m} f\big(x_{\tau^{-1}(1)}, \ldots, x_{\tau^{-1}(m)}, \eta \setminus \delta_{x_1} \setminus \ldots \setminus \delta_{x_m}\big) \, \dd \eta^{(m)}(x_1, \ldots, x_m) \bigg] \\
		&\quad = \E \bigg[ \int_{\X^m} f\big(x_1, \ldots, x_m, \eta \setminus \delta_{x_{\tau(1)}} \setminus \ldots \setminus \delta_{x_{\tau(m)}}\big) \, \dd \eta^{(m)}(x_1, \ldots, x_m) \bigg] \\
		&\quad = \E \bigg[ \int_{\X^m} f\big(x_1, \ldots, x_m, \eta \setminus \delta_{x_{1}} \setminus \ldots \setminus \delta_{x_{m}}\big) \, \dd \eta^{(m)}(x_1, \ldots, x_m) \bigg] \\
		&\quad = \int_\NN \int_{\X^m} f(x_1, \ldots, x_m, \mu) \, \kappa_m(x_1, \ldots, x_m, \mu) \, \dd \lambda^m(x_1, \ldots, x_m) \, \dd \PP^\eta(\mu) .
	\end{align*}
	\citep[Note that the symmetry of factorial measures is shown at the end of Appendix A.1 of][]{LP:2017}.
\end{prf}

\vspace{3mm}

Lemma \ref{LEMMA symmetry of kappa_m} essentially states that the existence of a Gibbs process with PI $\kappa$ implies that the $\kappa_m$ corresponding to $\kappa$ are necessarily symmetric in the first $m$ components. This is equivalent to the cocycle relation \eqref{cocycle assumption}. Hence the justification to always assume that $\kappa$ obeys the cocycle condition. 

In statistical physics, the finiteness of the partition function \eqref{definition of partition function} usually has to be stated as a separate assumption. The point process theoretic definition of the Gibbs process via the GNZ equations already contains essential information implicitly, as the following lemma shows. For less general state spaces this is known from \cite{NZ:1979} and \cite{MWM:1979}, but their techniques are generic and transfer to the abstract setting virtually unchanged. As the proof is short, we provide details. We write $B^\mathsf{c} = \X \setminus B$ for the complement of a set and $\nu_B(\cdot) = \nu( \cdot \cap B)$ for the restriction of a measure $\nu$ on $\X$ to $B$.
\begin{lemma} \label{LEMMA partition function finiteness property}
	If $\eta$ is a Gibbs process with PI $\kappa$, and $B \in \mathcal{X}$ with $\PP\big(
	\eta(B) < \infty \big) = 1$, then $\PP \big( Z_B(\eta_{B^\mathsf{c}}) < \infty \big) = 1$.
\end{lemma}
\begin{prf}
	For any measurable function $g : \NN \to [0, \infty]$ and any $D \in \mathcal{X}$, Proposition \ref{PROP properties of factorial measures due to Last/Penrose} (ii) gives
	\begin{equation*}
		\E \Big[ g(\eta) \, \mathds{1}\{\eta(D) < \infty\} \Big]
		= \E \Big[ g(\eta) \, \mathds{1}\{\eta(D) = 0\} \Big] + \sum_{m = 1}^\infty \frac{1}{m!} \, \E \bigg[ \int_{D^m} g(\eta) \, \mathds{1}\{\eta(D) = m\} \, \dd \eta^{(m)}(x_1, \dots, x_m) \bigg] .
	\end{equation*}
	By Lemma \ref{LEMMA mult. GNZ}, this equals
	\begin{equation*}	
		\E \bigg[ \mathds{1}\{\eta(D) = 0\} \bigg( g(\eta) + \sum_{m = 1}^\infty \frac{1}{m!} \int_{D^m} g(\eta + \delta_{x_1} + \dotso + \delta_{x_m}) \, \kappa_m(x_1, \dots, x_m, \eta) \, \dd \lambda^m(x_1, \dots, x_m) \bigg) \bigg] .
	\end{equation*}
	Applied to $D = B$ and $g(\mu) = \mathds{1}\{ \mu_{B^\mathsf{c}} \in A \}$, for $A \in \mathcal{N}$, the previous calculation yields
	\begin{equation*}
		\PP(\eta_{B^\mathsf{c}} \in A)
		= \E \Big[ \mathds{1}\{ \eta(B) = 0 \} \cdot \mathds{1}_A(\eta_{B^\mathsf{c}}) \cdot Z_B(\eta_{B^\mathsf{c}}) \Big]
		= \E \Big[ \mathds{1}_A(\eta_{B^\mathsf{c}}) \cdot Z_B(\eta_{B^\mathsf{c}}) \cdot \PP\big( \eta(B) = 0 \mid \eta_{B^\mathsf{c}} \big) \Big].
	\end{equation*}
	As $A \in \mathcal{N}$ was arbitrary, we find that
	\begin{equation*}
		Z_B(\eta_{B^\mathsf{c}}) \cdot \PP\big( \eta(B) = 0 \mid \eta_{B^\mathsf{c}} \big)
		= 1 \quad \PP\text{-almost surely (a.s.)}
	\end{equation*}
	Minding the convention $\infty \cdot 0 = 0$, we conclude that $\PP\big( \PP(\eta(B) = 0 \mid \eta_{B^\mathsf{c}}) > 0 \big) = 1$ as well as $\PP\big( Z_B(\eta_{B^\mathsf{c}}) < \infty \big) = 1$.
\end{prf}

\vspace{3mm}

By our construction of $\NN$, every set $B \in \mathcal{X}_b$ qualifies for Lemma \ref{LEMMA partition function finiteness property}, as $\PP\big( \eta(B) < \infty \big) = 1$ is satisfied for any point process $\eta$. An immediate corollary from the proof reads as follows.
\begin{corollary} \label{COR properties of the partition function}
	Let $\eta$ be a Gibbs process with PI $\kappa$, and $B \in \mathcal{X}$ with $\PP\big(
	\eta(B) < \infty \big) = 1$. Then
	\begin{equation*}
		\PP\big( \eta(B) = 0 \mid \eta_{B^\mathsf{c}} \big)
		= \frac{1}{Z_B(\eta_{B^\mathsf{c}})} \quad \PP\text{-a.s.}
	\end{equation*}
\end{corollary}

\section{Janossy densities and correlation functions of Gibbs processes}
\label{SEC Jan., Cor. and Palm meas. of Gibbs processes}

In this section, we calculate the Janossy densities and correlation functions of a general Gibbs process. A primer and some basic facts about these quantities are given in Appendix \ref{Appendix facts about Jan. and Cor. meas.}. Consider the setting from Section \ref{SEC Definition} and suppose that $\kappa$ satisfies the cocycle property \eqref{cocycle assumption}.

\begin{lemma} \label{LEMMA Janossy densities of Gibbs processes}
	Let $\eta$ be a Gibbs process in $\X$ with PI $\kappa$. Then $\eta$ has the Janossy densities
	\begin{equation*}
		j_{\eta, B, m}(x_1, \dots, x_m)
		= \frac{1}{m!} \, \E \Big[ \mathds{1}\big\{ \eta(B) = 0 \big\} \, \kappa_m(x_1, \dots, x_m, \eta) \Big] \, \mathds{1}_{B^m}(x_1, \dots, x_m) ,
	\end{equation*}
	where $(x_1, \dots, x_m) \in \X^m$, $m \in \N$, and $B \in \mathcal{X}$.
\end{lemma}
\begin{prf}
	Fix $B \in \mathcal{X}$ and $m \in \N$. For any $D \in \mathcal{X}^{\otimes m}$, Lemma \ref{LEMMA mult. GNZ} yields
	\begin{align*}
		J_{\eta, B, m}(D)
		&= \frac{1}{m!} \, \E\bigg[ \int_{\X^m} \mathds{1}\big\{ \eta(B) = m \big\} \, \mathds{1}_{D}(x_1, \dots, x_m) \, \dd \eta_B^{(m)}(x_1, \dots, x_m) \bigg] \\
		&= \int_{D \cap B^m} \frac{1}{m!} \, \E\Big[ \mathds{1}\big\{ \eta(B) = 0 \big\} \, \kappa_m(x_1, \dots, x_m, \eta) \Big] \, \dd \lambda^m(x_1, \dots, x_m) .
	\end{align*}
\end{prf} 

We now turn to the correlation functions of a Gibbs process.
\begin{lemma} \label{LEMMA correlation functions of Gibbs processes}
	Let $\eta$ be a Gibbs process in $\X$ with PI $\kappa$. Then $\eta$ has the correlation functions
	\begin{equation*}
		\rho_{\eta, m}(x_1, \dots, x_m)
		= \E \big[ \kappa_m(x_1, \dots, x_m, \eta) \big] ,
	\end{equation*}
	where $(x_1, \dots, x_m) \in \X^m$ and $m \in \N$.
\end{lemma}
\begin{prf}
	Fix $m \in \N$ and $D \in \mathcal{X}^{\otimes m}$. Then Lemma \ref{LEMMA mult. GNZ} implies that the factorial moment measure can be written as
	\begin{equation*}
		\alpha_{\eta, m}(D)
		= \E \bigg[ \int_{\X^m} \mathds{1}_D(x_1, \dots, x_m) \, \dd \eta^{(m)}(x_1, \dots, x_m) \bigg]
		= \int_{D} \E \big[ \kappa_m(x_1, \dots, x_m, \eta) \big] \, \dd \lambda^m(x_1, \dots, x_m) .
	\end{equation*}
\end{prf}
\begin{remark} \label{RMK on the correlation functions of the Gibbs process}
	With the knowledge about the correlation functions we can add to the result on the finiteness of the partition function in Lemma \ref{LEMMA partition function finiteness property}. Let $\eta$ be a Gibbs process in $\X$ with PI $\kappa$, and fix $B \in \mathcal{X}$. By definition of the partition function, we have
	\begin{equation*}
		\E \big[ Z_B(\eta) \big]
		= 1 + \sum_{m = 1}^\infty \frac{1}{m!} \int_{B^m} \E\big[ \kappa_m(x_1, \dots, x_m, \eta) \big] \, \dd \lambda^m(x_1, \dots, x_m) ,
	\end{equation*}
	and, using Lemma \ref{LEMMA correlation functions of Gibbs processes}, we obtain (in case $\eta(B)$ is finite almost surely)
	\begin{equation*}
		\E \big[ Z_B(\eta) \big]
		= 1 + \sum_{m = 1}^\infty \frac{\alpha_{\eta, m}(B^m)}{m!}
		= \E\big[ 2^{\eta(B)} \big]
	\end{equation*}
	as in the proof of Theorem \ref{THM repres. of Janossy meas. via factorial mom. meas.}. Hence, $\PP\big( Z_B(\eta) < \infty \big) = 1$ whenever $\E\big[ 2^{\eta(B)} \big] < \infty$. Remark \ref{RMK convergence of fact. moment meas. series} gives an overview as to when this property is satisfied.
\end{remark}

\section{Finite Gibbs processes}
\label{SEC finite Gibbs processes}

Under the assumption that the partition function is finite, we can explicitly specify a probability distribution which qualifies as that of a finite Gibbs process, even in an abstract measurable space $(\X, \mathcal{X})$. In fact, we show in Lemma \ref{LEMMA characterization of finite Gibbs process via kappa_m} that all finite Gibbs processes are necessarily of this form. For $C \in \mathcal{X}$ we denote by $\NN_C$ the set of all measures $\mu \in \NN$ with $\mu(C^\mathsf{c}) = 0$, and we denote by $\NN_f$ the set of all finite measures from $\NN$.

For $C \in \mathcal{X}$ and $\psi \in \NN$ define
\begin{equation*}
	\kappa^{(C, \psi)}(x, \mu)
	= \kappa(x, \psi + \mu) \cdot \mathds{1}_C(x), \quad x \in \X, ~ \mu \in \NN .
\end{equation*}
The map $\kappa^{(C, \psi)} : \X \times \NN \to [0, \infty)$ is measurable and inherits the cocycle property from $\kappa$. As in Section \ref{SEC Definition} we define the symmetric functions $\kappa_m^{(C, \psi)}$ ($m \in \N$) and the partition functions $Z_B^{(C, \psi)}$ ($B \in \mathcal{X}$) corresponding to $\kappa^{(C, \psi)}$. By definition, the partition functions relate in the following way.
\begin{lemma} \label{LEMMA relation of the partition functions}
	For $B, C \in \mathcal{X}$ and $\psi, \nu \in \NN$, it holds that $Z_B^{(C, \psi)}(\nu) = Z_{B \cap C}(\psi + \nu)$.
\end{lemma}

In the case of finite partition functions, we can specify the distribution of Gibbs processes on sets of finite $\lambda$-measure.
\begin{lemma} \label{LEMMA finite Gibbs process via kappa_m}
	Let $C \in \mathcal{X}$ with $\lambda(C) < \infty$ and $\psi \in \NN$ be such that $Z_C(\psi) < \infty$. Consider a point process $\xi$ in $\X$ with distribution
	\begin{equation} \label{distribution of finite Gibbs process}
		\PP^\xi(\cdot)
		= \frac{1}{Z_C(\psi)} \bigg( \mathds{1}\{ \mathbf{0} \in \cdot \, \} + \sum_{m = 1}^\infty \frac{1}{m!} \int_{C^m} \mathds{1}\Big\{ \sum_{j = 1}^m \delta_{x_j} \in \cdot \, \Big\} \, \kappa_m(x_1, \dots, x_m, \psi) \, \dd\lambda^m(x_1, \dots, x_m) \bigg) .
	\end{equation}
	Then $\xi$ is a Gibbs process with PI $\kappa^{(C, \psi)}$ and reference measure $\lambda$. Moreover, $\PP(\xi \in \NN_C \cap \NN_f) = 1$.
\end{lemma}
\begin{prf}
	By definition of $Z_C(\psi)$, $\PP^\xi$ is indeed a probability measure on $\NN$, and we clearly have $\PP^\xi(\NN_C \cap \NN_f) = 1$. Let $\mathrm{D}_{C, \psi} : \NN \to [0, \infty]$ be defined by
	\begin{equation*}
		\mathrm{D}_{C, \psi}(\mu)
		= \frac{\mathrm{e}^{\lambda(C)}}{Z_C(\psi)} \bigg( \mathds{1}\big\{ \mu(C) = 0 \big\} + \sum_{m = 1}^\infty \frac{1}{m!} \, \mathds{1}\big\{ \mu(C) = m \big\} \int_{C^m} \kappa_m(x_1, \dots, x_m, \psi) \, \dd\mu^{(m)}(x_1, \dots, x_m) \bigg) .
	\end{equation*}
	The map $\mathrm{D}_{C, \psi}$ is measurable by Lemma \ref{LEMMA measurability of factorial-measure-integrals} and we have
 	\begin{equation*}
	 	\mathrm{D}_{C, \psi}\Big( \sum_{j = 1}^m \delta_{x_j} \Big) 
	 	= \frac{\mathrm{e}^{\lambda(C)}}{Z_C(\psi)} \cdot \kappa_m(x_1, \dots, x_m, \psi)
 	\end{equation*}
	for all $x_1, \dots, x_m \in C$ and each $m \in \N$. For $\mu \in \NN$ and $x \in C$, we have $\mathrm{D}_{C, \psi}(\mu + \delta_x) = \mathrm{D}_{C, \psi}(\mu) \cdot \kappa(x, \psi + \mu_C)$. Indeed, this follows immediately if $\mu(C) = 0$ or if $\mu(C) = \infty$. For $\mu \in \NN$ with $\mu(C) = m \in \N$, observe that, by Lemma \ref{LEMMA general lemma on factorial measures} and Proposition \ref{PROP representation of NN to R functions via factorial measures}, 
	\begin{align*}
		\mathrm{D}_{C, \psi}(\mu + \delta_x)
		&= \frac{\mathrm{e}^{\lambda(C)}}{Z_C(\psi)} \cdot \frac{1}{(m + 1)!} \int_{C^{m + 1}} \kappa_{m + 1}(x_1, \dots, x_{m + 1}, \psi) \, \dd(\mu + \delta_x)^{(m + 1)}(x_1, \dots, x_{m + 1}) \\
		&= \frac{\mathrm{e}^{\lambda(C)}}{Z_C(\psi)} \cdot \frac{1}{m!} \int_{C^m} \kappa_{m + 1}(x, x_1, \dots, x_{m}, \psi) \, \dd\mu^{(m)}(x_1, \dots, x_m) \\
		&= \frac{1}{m!} \int_{C^m} \mathrm{D}_{C, \psi}\Big( \sum_{j = 1}^m \delta_{x_j} \Big) \cdot \kappa\Big( x, \psi + \sum_{j = 1}^m \delta_{x_j} \Big) \, \dd\mu^{(m)}(x_1, \dots, x_m) \\
		&= \mathrm{D}_{C, \psi}(\mu) \cdot \kappa(x, \psi + \mu_C) .
	\end{align*}
	Let $\Phi$ be a Poisson process in $\X$ with intensity measure $\lambda$. The definition of $\PP^\xi$ in \eqref{distribution of finite Gibbs process} and Corollary \ref{COR Janossy representation of Poisson processes} imply that, for any measurable map $F : \NN \to [0, \infty]$,
	\begin{equation} \label{finite Gibbs processes construction, map D}
		\E\big[ F(\xi) \big]
		= \E\big[ F(\Phi_C) \cdot \mathrm{D}_{C, \psi}(\Phi) \big] .
	\end{equation}
	Applied to $F(\mu) = \int_\X f(x, \mu) \, \dd\mu(x)$, where $f : \X \times \NN \to [0, \infty]$ is a measurable function, and combined with the Mecke equation \citep[Theorem 4.1 of][]{LP:2017}, this yields
	\begin{equation*}
		\E\bigg[ \int_\X f(x, \xi) \, \dd\xi(x) \bigg]
		= \E\bigg[ \int_C f(x, \Phi_C + \delta_x) \, \mathrm{D}_{C, \psi}(\Phi + \delta_x) \, \dd\lambda(x) \bigg].
	\end{equation*}
	With the properties of $\mathrm{D}_{C, \psi}$, the term on the right hand side equals
	\begin{equation*}	
		\E\bigg[ \int_C f(x, \Phi_C + \delta_x) \, \kappa(x, \psi + \Phi_C) \, \dd\lambda(x) \cdot \mathrm{D}_{C, \psi}(\Phi) \bigg]
	\end{equation*}
	and another application of \eqref{finite Gibbs processes construction, map D} to $F(\mu) = \int_C f(x, \mu + \delta_x) \, \kappa(x, \psi + \mu) \, \dd\lambda(x)$ gives
	\begin{equation*}
		\E\bigg[ \int_\X f(x, \xi) \, \dd\xi(x) \bigg]
		= \E\bigg[ \int_\X f(x, \xi + \delta_x) \, \kappa^{(C, \psi)}(x, \xi) \, \dd\lambda(x) \bigg] ,
	\end{equation*}
	which concludes the proof.
\end{prf}

\vspace{3mm}

It is obvious that the distribution defined in Lemma \ref{LEMMA finite Gibbs process via kappa_m} makes sense as soon as $Z_C(\psi) < \infty$ even without the assumption that $\lambda(C) < \infty$. Indeed, the previous result easily generalizes.
\begin{corollary} \label{COR general finite Gibbs process via kappa_m}
	Let $C \in \mathcal{X}$ and $\psi \in \NN$ be such that $Z_C(\psi) < \infty$. A point process $\xi$ in $\X$ with distribution given through \eqref{distribution of finite Gibbs process} is a finite Gibbs process with PI $\kappa^{(C, \psi)}$ and reference measure $\lambda$.
\end{corollary}
\begin{prf}
	By the definition of partition functions, and monotone convergence, we have
	\begin{equation*}
		\lim_{\ell \to \infty} Z_{C \cap B_\ell}(\psi)
		= Z_C(\psi) ,
	\end{equation*}
	where $1 \leq Z_{C \cap B_\ell}(\psi) \leq Z_C(\psi) < \infty$ for all $\ell \in \N$. For each $\ell \in \N$, let $\xi_\ell$ be a point process with distribution as in Lemma \ref{LEMMA finite Gibbs process via kappa_m}, but with $C$ replaced by $C \cap B_\ell$. As such, each $\xi_\ell$ is a Gibbs process with PI $\kappa^{(C \cap B_\ell, \psi)}$. Let $F : \NN \to [0, \infty]$ be an arbitrary measurable function, and $F_1, F_2, \ldots : \NN \to [0, \infty]$ any sequence of measurable functions with $F_\ell(\mu) \nearrow F(\mu)$ for each $\mu \in \NN$ as $\ell \to \infty$. By construction of the processes, and monotone convergence, we have
	\begin{align*}
		\lim_{\ell \to \infty} \E\big[ F_\ell(\xi_\ell) \big]
		&= \lim_{\ell \to \infty} \frac{1}{Z_{C \cap B_\ell}(\psi)} \bigg( F_\ell(\mathbf{0}) + \sum_{m = 1}^\infty \frac{1}{m!} \int_{(C \cap B_\ell)^m} F_\ell\Big( \sum_{j = 1}^m \delta_{x_j} \Big) \, \kappa_m(x_1, \dots, x_m, \psi) \, \dd\lambda^m(x_1, \dots, x_m) \bigg) \\
		&= \frac{1}{Z_{C}(\psi)} \bigg( F(\mathbf{0}) + \sum_{m = 1}^\infty \frac{1}{m!} \int_{C^m} F\Big( \sum_{j = 1}^m \delta_{x_j} \Big) \, \kappa_m(x_1, \dots, x_m, \psi) \, \dd\lambda^m(x_1, \dots, x_m) \bigg) \\
		&= \E\big[ F(\xi) \big] .
	\end{align*}
	Applied twice, this observation, together with the GNZ equations, yields
	\begin{align*}
		\E\bigg[ \int_\X f(x, \xi) \, \dd\xi(x) \bigg]
		= \lim_{\ell \to \infty} \E\bigg[ \int_\X f(x, \xi_\ell) \, \dd\xi_\ell(x) \bigg]
		&= \lim_{\ell \to \infty} \E\bigg[ \int_\X f(x, \xi_\ell + \delta_x) \, \kappa^{(C, \psi)}(x, \xi_\ell) \, \mathds{1}_{B_\ell}(x) \, \dd\lambda(x) \bigg] \\
		&= \E\bigg[ \int_\X f(x, \xi + \delta_x) \, \kappa^{(C, \psi)}(x, \xi) \, \dd\lambda(x) \bigg]
	\end{align*}
	for every measurable map $f : \X \times \NN \to [0, \infty]$, which concludes the proof.
\end{prf}

\vspace{3mm}

The measure $\psi$ in Lemma \ref{LEMMA finite Gibbs process via kappa_m} or Corollary \ref{COR general finite Gibbs process via kappa_m} is often called boundary condition. In particular, if $\psi = \mathbf{0}$, the construction gives a Gibbs process whose PI is $\kappa$ restricted to the domain $C$. With the upcoming result we establish that any finite Gibbs process in a window $C \in \mathcal{X}$ has to have the distribution from Lemma \ref{LEMMA finite Gibbs process via kappa_m}. In particular, finite Gibbs processes are unique in distribution. While it was known that finite partition functions are essential to produce finite Gibbs processes, the equivalence seems not to be stated in the literature and finite Gibbs processes were always constructed on bounded sets. The given result shows that boundedness of the domain is inessential, finite partition functions are what matters.
\begin{lemma} \label{LEMMA characterization of finite Gibbs process via kappa_m}
	Let $C \in \mathcal{X}$ and $\psi \in \NN$. A finite Gibbs process with PI $\kappa^{(C, \psi)}$ exists if, and only if, $Z_C(\psi) < \infty$. If $\xi$ is such a finite Gibbs process in $\X$ with PI $\kappa^{(C, \psi)}$ and reference measure $\lambda$, then the distribution of $\xi$ is given by \eqref{distribution of finite Gibbs process}.
\end{lemma}
\begin{prf}
	If $Z_C(\psi) < \infty$, a finite Gibbs process with PI $\kappa^{(C, \psi)}$ exists by Corollary \ref{COR general finite Gibbs process via kappa_m}. Now, let $\xi$ be a finite Gibbs process with PI $\kappa^{(C, \psi)}$. The GNZ equations yield
	\begin{equation*}
		\E\big[ \xi(C^\mathsf{c}) \big]
		= \E\bigg[ \int_\X \mathds{1}_{C^\mathsf{c}}(x) \, \dd\xi(x) \bigg]
		= \E\bigg[ \int_\X \mathds{1}_{C^\mathsf{c}}(x) \, \kappa(x, \psi + \xi) \, \mathds{1}_C(x) \, \dd\lambda(x) \bigg]
		= 0 ,
	\end{equation*}
	so $\PP^\xi(\NN_C) = 1$. By Lemmata \ref{LEMMA relation of the partition functions} and \ref{LEMMA partition function finiteness property}, we get
	\begin{equation*}
		Z_C(\psi)
		= Z^{(C, \psi)}_C(\mathbf{0})
		= Z^{(C, \psi)}_C(\xi_{C^\mathsf{c}})
		< \infty ,
	\end{equation*}
	with the statements involving $\xi_{C^\mathsf{c}}$ holding almost surely. By Corollary \ref{COR properties of the partition function} and Lemma \ref{LEMMA Janossy densities of Gibbs processes}, the Janossy densities of $\xi$ on the full domain $C$ are $J_{\xi, C, 0} = Z_C(\psi)^{-1}$ and
	\begin{equation*}
		j_{\xi, C, m}(x_1, \dots, x_m)
		= \frac{1}{m!} \cdot \frac{\kappa_m(x_1, \dots, x_m, \psi)}{Z_C(\psi)} \cdot \mathds{1}_{C^m}(x_1, \dots, x_m) ,
	\end{equation*}
	and Lemma \ref{LEMMA expec. repres. via Janossy measures} concludes the proof.
\end{prf}

\vspace{3mm}

In the proof of the previous lemma we have used a conclusion from Corollary \ref{COR properties of the partition function} and Lemma \ref{LEMMA Janossy densities of Gibbs processes} concerning the Janossy densities of a finite Gibbs process on its full domain. This we state explicitly at this point, and we add a further conclusion from the proof of Lemma \ref{LEMMA finite Gibbs process via kappa_m}. We denote by $\Pi_\nu$ the distribution of a Poisson process in $\X$ with locally finite intensity measure $\nu$ on $\X$.
\begin{corollary} \label{COR full domain Janossy densities of Gibbs processes}
	Let $C \in \mathcal{X}$ and $\psi \in \NN$. If $\xi$ is a finite Gibbs process with PI $\kappa^{(C, \psi)}$, then 
	\begin{equation*}
		J_{\xi, C, 0}
		= \PP\big( \xi(C) = 0 \big)
		= \frac{1}{Z_C(\psi)} .
	\end{equation*}
	Moreover, the Janossy densities of $\xi$ on the full domain $C$ are
	\begin{equation*}
		j_{\xi, C, m}(x_1, \dots, x_m)
		= \frac{1}{m!} \cdot \frac{\kappa_m(x_1, \dots, x_m, \psi)}{Z_C(\psi)} \cdot \mathds{1}_{C^m}(x_1, \dots, x_m) ,
	\end{equation*}
	where $x_1, \dots, x_m \in \X$ and $m \in \N$. If, additionally, $\lambda(C) < \infty$, then $\PP^\xi$ is absolutely continuous with respect to $\Pi_{\lambda_{C}}$ with density function $\mathrm{D}_{C, \psi}$ from the proof of Lemma \ref{LEMMA finite Gibbs process via kappa_m}.
\end{corollary}

\section{Hamiltonians in abstract measurable spaces}
\label{SEC Hamiltonians}

Let $\kappa : \X \times \NN \to [0, \infty)$ be measurable and satisfy the cocycle assumption \eqref{cocycle assumption}. Based on $\kappa$ we define the \textit{Hamiltonian} $H : \NN \times \NN \to (- \infty, \infty]$ as
\begin{equation*}
	H(\mu , \psi)
	= \infty \cdot \mathds{1}\big\{ \mu(\X) = \infty \big\} - \sum_{m = 1}^\infty \mathds{1}\big\{ \mu(\X) = m \big\} \cdot \log\bigg( \frac{1}{m!} \int_{\X^m} \kappa_m(x_1, \ldots, x_m, \psi) \, \dd \mu^{(m)}(x_1, \ldots, x_m) \bigg) .
\end{equation*}
\begin{lemma} \label{LEMMA Hamiltonian is well-defined}
	The Hamiltonian is well-defined and measurable.
\end{lemma}
\begin{prf}
	For any $m \in \N$ and $\psi \in \NN$, the map
	\begin{equation*}
		G(\mu) 
		= \int_{\X^m} \kappa_m(x_1, \ldots, x_m, \psi) \, \dd \mu^{(m)}(x_1, \ldots, x_m)
	\end{equation*}
	is measurable by Lemma \ref{LEMMA measurability of factorial-measure-integrals}, and Lemma \ref{LEMMA finite functions on Dirac measures are finite on NN_f} implies that $G(\mu) < \infty$ for all $\mu \in \NN_f$. Hence, taking indicator functions into account, the Hamiltonian cannot attain the value $- \infty$ and is well-defined. Using that $[0, \infty) \ni s \mapsto -\log(s) \in (- \infty, \infty]$ is measurable, Lemma \ref{LEMMA measurability of factorial-measure-integrals} implies that $H$ is measurable.
\end{prf}

\vspace{3mm}

Notice that $H(\mathbf{0}, \psi) = 0$ for each $\psi \in \NN$, and $H(\mu, \psi) = \infty$ if $\mu(\X) = \infty$. Moreover, if $\mu = \sum_{j = 1}^m \delta_{x_j}$, for some $x_1, \dots, x_m \in \X$ and $m \in \N$, then, by the symmetry of $\kappa_m$,
\begin{equation*}
	H(\mu, \psi)
	= - \log\big( \kappa_m(x_1, \dots, x_m, \psi) \big) .
\end{equation*}
Solving for $\kappa_m$, the last line is \eqref{hamiltonian for Dirac sum}. For $\mu = \sum_{i = 1}^{k} \delta_{x_i}$, $\nu = \sum_{j = 1}^{m} \delta_{y_j}$ (with $x_1, \dots, x_k, y_1, \dots, y_m \in \X$, $k, m \in \N$) and $\psi \in \NN$, the Hamiltonian satisfies
\begin{align*}
	H(\mu, \psi) + H(\nu, \psi + \mu)
	&= - \log\big( \kappa_k(x_1, \dots, x_k, \psi) \cdot \kappa_m(y_1, \dots, y_m, \psi + \delta_{x_1} + \dotso + \delta_{x_k}) \big) \\
	&= - \log\big( \kappa_{k + m}(x_1, \dots, x_k, y_1, \dots, y_m, \psi) \big) \\
	&= H(\mu + \nu, \psi) ,
\end{align*}
which is the well known \textit{energy function property}. The property gives rise to the interpretation that $H(\mu, \psi)$ captures the interactions between points in the configuration $\mu$ as well as between points in $\mu$ and points in the boundary configuration $\psi$. The question arises whether the general definition of $H$ in measurable spaces admits the energy function property more generally.
\begin{lemma}[\bf Energy function property] \label{LEMMA energy function property}
	Let $\mu, \nu, \psi \in \NN$. The Hamiltonian satisfies $H(\mu + \nu, \psi) = H(\mu, \psi) + H(\nu, \psi + \mu)$.
\end{lemma}
\begin{prf}
	If $\mu(\X) = 0$ or $\nu(\X) = 0$, the claim is trivially true. The same can be said if $\mu(\X) = \infty$ or $\nu(\X) = \infty$. Thus, we assume that $\mu(\X) = k \in \N$ and $\nu(\X) = m \in \N$. Applying Proposition \ref{PROP representation of NN to R functions via factorial measures} to the map $\mathrm{e}^{- H(\,\cdot\, , \psi)} : \NN \to [0, \infty)$ gives
	\begin{equation*}
		\exp\big( - H(\mu + \nu, \psi) \big)
		= \frac{1}{(k + m)!} \int_{\X^{k + m}} \exp\bigg( - H\Big( \sum_{j = 1}^{k + m} \delta_{z_j}, \, \psi \Big) \bigg) \, \dd(\mu + \nu)^{(k + m)}(z_1, \dots, z_{k + m}) .
	\end{equation*}
	By Lemma \ref{LEMMA general lemma on factorial measures} and the observation ahead of the lemma at hand, the right hand side is equal to
	\begin{align*}
		&\frac{1}{k! \cdot m!} \int_{\X^{k}} \int_{\X^{m}} \exp\bigg( - H\Big( \sum_{i = 1}^{k} \delta_{x_i} + \sum_{j = 1}^{m} \delta_{y_j} , \, \psi \Big) \bigg) \, \dd\nu^{(m)}(y_1, \dots, y_m) \, \dd\mu^{(k)}(x_1, \dots, x_k) \\
		&~~ = \frac{1}{k!} \int_{\X^{k}} \exp\bigg( - H\Big( \sum_{i = 1}^{k} \delta_{x_i} , \, \psi \Big) \bigg) \bigg[ \frac{1}{m!} \int_{\X^{m}} \exp\bigg( - H\Big( \sum_{j = 1}^{m} \delta_{y_j} , \, \psi + \sum_{i = 1}^{k} \delta_{x_i} \Big) \bigg) \, \dd\nu^{(m)}(y_1, \dots, y_m) \bigg] \, \dd\mu^{(k)}(x_1, \dots, x_k) .
	\end{align*}
	Proposition \ref{PROP representation of NN to R functions via factorial measures} applied to the map $\mathrm{e}^{- H(\,\cdot\, , \psi \, + \, \delta_{x_1} \, + \, \ldots \, + \, \delta_{x_k})} : \NN \to [0, \infty)$ yields
	\begin{equation*}
		\exp\big( - H(\mu + \nu, \psi) \big)
		= \frac{1}{k!} \int_{\X^{k}} \exp\bigg( - H\Big( \sum_{i = 1}^{k} \delta_{x_i} , \, \psi \Big) \bigg) \, \exp\bigg( - H\Big( \nu , \psi + \sum_{i = 1}^{k} \delta_{x_i} \Big) \bigg) \, \dd\mu^{(k)}(x_1, \dots, x_k) ,
	\end{equation*}
	and a final application of Proposition \ref{PROP representation of NN to R functions via factorial measures} to the map $\mathrm{e}^{- H(\,\cdot\, , \psi)} \, \mathrm{e}^{- H(\nu, \psi \, + \, \cdot\,)} : \NN \to [0, \infty)$ gives
	\begin{equation*}
		\exp\big( - H(\mu + \nu , \psi) \big)
		= \exp\big( - H(\mu , \psi) \big) \, \exp\big( - H(\nu , \psi + \mu) \big)
		= \exp \big( - H(\mu , \psi) - H(\nu , \psi + \mu) \big) .
	\end{equation*}
	Taking logarithms concludes the proof.
\end{prf}

\vspace{3mm}

Note that the \textit{energy function}, which is very prominent in the context of Gibbs measures, is simply given by $\mu \mapsto H(\mu , \mathbf{0})$.
\begin{remark}[\bf Hereditarity] \label{RMK hereditarity property of H}
	Let $\psi \in \NN$. The Hamiltonian $H$ is hereditary, that is, if $\mu \in \NN$ is such that $H(\mu, \psi) = \infty$, then
	\begin{equation*}
		H(\mu + \nu , \psi)
		= H(\mu , \psi) + H(\nu , \mu + \psi)
		= \infty
	\end{equation*}
	for every $\nu \in \NN$, by Lemma \ref{LEMMA energy function property}.
\end{remark}

We now use the Hamiltonian to rewrite several of our previous observations about Gibbs processes. The formulae we obtain are perhaps more familiar to some readers. First of all, we write the partition function as an integral with respect to the distribution of a Poisson process. Note that while the partition function is defined on general measurable sets, the following representation only works for sets of finite $\lambda$-measure. The result is immediate from \eqref{hamiltonian for Dirac sum} and Corollary \ref{COR Janossy representation of Poisson processes}.
\begin{corollary} \label{LEMMA partition function representation via Poisson process}
	Fix a set $B \in \mathcal{X}$ with $\lambda(B) < \infty$. For each $\psi \in \NN$, the partition function $Z_B(\psi)$ can be written as
	\begin{equation*}
		Z_B(\psi)
		= \mathrm{e}^{\lambda(B)} \int_\NN \mathrm{e}^{- H(\mu , \psi)} \, \dd \Pi_{\lambda_B}(\mu) .
	\end{equation*}
\end{corollary}

In Lemma \ref{LEMMA characterization of finite Gibbs process via kappa_m} we learned that the existence of finite Gibbs processes is invariably linked to the finiteness of partition functions and that the distribution of such a finite Gibbs process can be stated explicitly. On sets of finite $\lambda$-measure this distribution can be written in terms of the distribution of a Poisson process, similar to the observation in the previous lemma. The following result is an immediate consequence of Lemmata \ref{LEMMA finite Gibbs process via kappa_m} and \ref{LEMMA characterization of finite Gibbs process via kappa_m} and Corollary \ref{COR Janossy representation of Poisson processes}.
\begin{corollary} \label{COR finite Gibbs process via H}
	Let $C \in \mathcal{X}$ with $\lambda(C) < \infty$ and $\psi \in \NN$ be such that $Z_C(\psi) < \infty$. A point process $\xi$ in $\X$ is a finite Gibbs process with PI $\kappa^{(C, \psi)}$ if, and only if, the distribution of $\xi$ is given by
	\begin{equation*}
		\PP^\xi(\cdot) 
		= \frac{\mathrm{e}^{\lambda(C)}}{Z_C(\psi)} \int_\NN \mathds{1}\{ \mu \in \cdot \} \, \mathrm{e}^{- H(\mu , \psi)} \dd \Pi_{\lambda_C}(\mu).
	\end{equation*}
\end{corollary}

Both Corollaries \ref{COR full domain Janossy densities of Gibbs processes} and \ref{COR finite Gibbs process via H} express the distribution of a finite Gibbs process in terms of a density function with respect to the Poisson process distribution. Of course these density functions have to agree almost everywhere, but from their definitions and Proposition \ref{PROP representation of NN to R functions via factorial measures} it even follows that they are identical.

\section{The DLR equations}
\label{SEC DLR equations}

In this section we state, for the sake of completeness, the so-called DLR equations, named after Dobrushin, Lanford, and Ruelle, see \cite{D:1968:1} \citep[also][]{D:1968:2, D:1969}, \cite{LR:1969} and \cite{R:1970}. Even though we state the result in hitherto unprecedented generality, namely for an arbitrary measurable state space, we do not give a proof. This is due to the fact that the techniques from \cite{NZ:1979} and \cite{MWM:1979}, who established the equivalence of the DLR and GNZ equations, also work in this greater generality with only minor adaptations. In fact, the calculation which shows that a Gibbs process satisfies the DLR equation is similar to the proof of Lemma \ref{LEMMA partition function finiteness property} using Corollary \ref{COR Janossy representation of Poisson processes}. 

We want to emphasize that the DLR equations are of extraordinary importance in statistical physics and are often used as a definition of the Gibbs process (instead of the GNZ equations). The fact that no topological structure is required on the state space shows how the characterization captures essential structural features of Gibbs processes. A reader who is entirely unfamiliar with the DLR equations might find it helpful to consult an introductory course on Gibbs processes, like \cite{J:2017:script} or \cite{D:2019}.

\begin{theorem}[\bf The DLR equations] \label{THM DLR equations via H}
	A point process $\eta$ in $\X$ is a Gibbs process with PI $\kappa$ if, and only if, for all $B \in \mathcal{X}_b$, the process satisfies $\PP\big( Z_B(\eta_{B^\mathsf{c}}) < \infty \big) = 1$ as well as
	\begin{equation*}
		\E \big[ F(\eta_B) \mid \eta_{B^\mathsf{c}} \big] = \frac{\mathrm{e}^{\lambda(B)}}{Z_B(\eta_{B^\mathsf{c}})} \int_\NN F(\mu) \, \mathrm{e}^{- H(\mu , \eta_{B^\mathsf{c}})} \, \dd \Pi_{\lambda_B}(\mu) \quad \PP\text{-a.s.}
	\end{equation*}
	for each measurable function $F : \NN \to [0, \infty]$.
\end{theorem}
Recalling that we can express $\kappa_m$ via the Hamiltonian $H$ as in Section \ref{SEC Hamiltonians}, a straightforward application of Corollary \ref{COR Janossy representation of Poisson processes} allows one to rewrite the DLR equations in terms of $\kappa_m$, but we do without an explicit statement of this reformulation. However, we use the formalism of finite Gibbs processes to rewrite the equations which brings us back to the equilibrium equations of \cite{R:1970}. Let $B \in \mathcal{X}$ with $\lambda(B) < \infty$ and $\psi \in \NN$ such that $Z_B(\psi) < \infty$. Denote by $\mathrm{P}_{B, \psi}$ the distribution from Corollary \ref{COR finite Gibbs process via H} of a finite Gibbs process with PI $\kappa^{(B, \psi)}$ as a measure on $(\NN, \mathcal{N})$. In particular, $\mathrm{P}_{B, \psi}(\NN_B) = 1$ and $\mathrm{P}_{B, \psi}$ is absolutely continuous with respect to $\Pi_{\lambda_B}$ with density function
\begin{equation*}
	\frac{\dd \mathrm{P}_{B, \psi}}{\dd \Pi_{\lambda_B}}(\mu)
	= \frac{e^{\lambda(B)}}{Z_B(\psi)} \cdot \mathrm{e}^{- H(\mu , \psi)} .
\end{equation*}
If $Z_B(\psi) = \infty$, Lemma \ref{LEMMA characterization of finite Gibbs process via kappa_m} implies that no finite Gibbs process with PI $\kappa^{(B, \psi)}$ exists, so we put $\mathrm{P}_{B, \psi} \equiv 0$. Let $p_B(\mu) = \mu_B$ be the restriction mapping on $\NN$ and $\mathcal{N}_B = \sigma(p_B)$ the $\sigma$-field of $B$-local events as discussed in Appendix \ref{Appendix local events and functions}.

\begin{lemma} \label{LEMMA measurability of finite Gibbs measure wrt. boundary cond.}
	Fix $B \in \mathcal{X}$ with $\lambda(B) < \infty$. The map
	\begin{equation*}
		\NN \times \NN \ni (\psi, \nu) \mapsto
		\int_\NN G(\mu, \nu) \, \dd \mathrm{P}_{B, \psi_{B^\mathsf{c}}}(\mu)
	\end{equation*}
	is $\mathcal{N}_{B^\mathsf{c}} \otimes \mathcal{N}$-measurable (hence $\mathcal{N} \otimes \mathcal{N}$-measurable) for every measurable map $G : \NN \times \NN \to [0, \infty]$.
\end{lemma}
\begin{prf}
	Since the partition function $Z_B : \NN \to [1, \infty]$ is measurable, the map $\psi \mapsto Z_B(\psi_{B^\mathsf{c}}) = Z_B \circ p_{B^\mathsf{c}}(\psi)$ is $\mathcal{N}_{B^\mathsf{c}}$-measurable. Thus, Fubini's theorem implies that
	\begin{equation*}
		(\psi, \nu) \mapsto
		\int_\NN G(\mu, \nu) \, \dd \mathrm{P}_{B, \psi_{B^\mathsf{c}}}(\mu)
		= \mathds{1}\big\{ Z_B(\psi_{B^\mathsf{c}}) < \infty \big\} \cdot \frac{\mathrm{e}^{\lambda(B)}}{Z_B(\psi_{B^\mathsf{c}})} \int_\NN G(\mu, \nu) \,  \mathrm{e}^{- H(\mu , \psi_{B^\mathsf{c}})} \, \dd \Pi_{\lambda_B}(\mu)
	\end{equation*}
	is $\mathcal{N}_{B^\mathsf{c}} \otimes \mathcal{N}$-measurable.
\end{prf}

\vspace{3mm}

With this new formalism we present a reformulation of the DLR equations in terms of finite Gibbs processes depending on the boundary condition. Notice that, while considerations via conditional probabilities go way back to \cite{D:1968:1}, the following corollary comes closest to the equilibrium equations stated by \cite{LR:1969} and \cite{R:1970}. The result follows from Theorem \ref{THM DLR equations via H} with Lemma \ref{LEMMA measurability of finite Gibbs measure wrt. boundary cond.} guaranteeing that all terms are well-defined.
\begin{corollary} \label{COR DLR reformulation in terms of finite Gibbs proc. distr.}
	A point process $\eta$ in $\X$ is a Gibbs process with PI $\kappa$ if, and only if,
	\begin{equation*}
	\E \big[ F(\eta) \big] = \E \bigg[ \int_\NN F(\mu + \eta_{B^\mathsf{c}}) \, \dd \mathrm{P}_{B, \eta_{B^\mathsf{c}}}(\mu) \bigg]
	\end{equation*}
	for all measurable functions $F : \NN \to [0, \infty]$ and each $B \in \mathcal{X}_b$.
\end{corollary}

Note that in the context of the DLR equations many authors use the notion of specifications, some even use it to define Gibbs measures, see \cite{P:1976} and \cite{NZ:1979}. As this representation does not fit too well in our point process theoretic notation, and is extensively detailed in the given references \cite[see also][]{B:2022}, we do not reiterate this notion.

\section{Some facts about the local convergence topology and new convergence results}
\label{SEC prelinimaries}

In this section we discuss important fundamentals which are used throughout. We stick to the setting from the beginning of Section \ref{SEC Introduction}. We discuss the concept of local convergence introduced by \cite{GZ:1993} and provide connections to Janossy and factorial moment measures. In the definition of local convergence we follow the recent publications by \cite{J:2019} and \cite{RZ:2020}, where the concept includes the use of local and tame functions as proposed by \cite{GZ:1993}, but other authors also use the term local convergence when only including local and bounded functions, see \cite{DV:2020} for one instance.

A function $F : \NN \to [- \infty, \infty]$ is called \textit{local} if there exists a set $B \in \mathcal{X}_b$ such that $F(\mu) = F(\mu_B)$ for every $\mu \in \NN$. We call such a function \textit{tame} if there exists a constant $C \geq 0$ and a set $B \in \mathcal{X}_b$ such that $|F(\mu)| \leq C \big( 1 + \mu(B) \big)$ for all $\mu \in \NN$. In particular, every bounded function $F$ is tame. Some properties of local functions are discussed in Appendix \ref{Appendix local events and functions}. If $\eta, \eta_1, \eta_2, \ldots$ are point processes in $\X$, we say that $(\eta_n)_{n \in \N}$ converges locally to $\eta$ if $\E[F(\eta_n)] \to \E[F(\eta)]$ as $n \to \infty$ for every measurable, local, and tame function $F : \NN \to [0, \infty)$ for which the expectations are finite. For short, we write $\eta_n \stackrel{loc}{\longrightarrow} \eta$. There is no need for the functions $F$ to be continuous as in many other modes of convergence, and thus no topological structure is needed on $\NN$.

It is easy to verify that the class of local and tame functions is a measure determining class and that local limits are unique in distribution. Moreover, the following result is an immediate consequence of Theorem 11.1.VII of \cite{DV:2008}.

\begin{prop} \label{PROP basic properties of local convergence}
	Let $\eta, \eta_1, \eta_2, \ldots$ be point processes in a complete separable metric space $\X$. Then local convergence is stronger than convergence in law, that is, if $\eta_n \stackrel{loc}{\longrightarrow} \eta$, then $\eta_n$ converges to $\eta$ in law.
\end{prop}

It is well-known that local convergences is equivalent to a suitable weak* convergence (in the functional analytic sense) of correlations functions. We state this result later on in this subsection, but first provide the connection between local convergence and convergence of the Janossy measures. This connection is used (rather implicitly) in the literature, for instance in the appendix of \cite{J:2019}, but to our knowledge this is the first time the results are stated separately. We state them in full abstraction and under weaker assumptions than were used previously. Note that a definition and basic properties of Janossy measures are discussed in Appendix \ref{Appendix facts about Jan. and Cor. meas.}.
\begin{lemma} \label{LEMMA conv. of Janossy meas. implies local conv.}
	Let $\eta, \eta_1, \eta_2, \ldots$ be point processes in $\X$ and fix a bounded set $B \in \mathcal{X}_b$. Assume there exists a map $c_B : \N \to [0, \infty)$ with $\sum_{m = 1}^\infty c_B(m) < \infty$ such that the Janossy measures of $\eta_n$ restricted to $B$ satisfy
	\begin{equation*}
		\sup_{n \in \N} J_{\eta_n, B, m}(B^m)
		\leq c_B(m)
	\end{equation*}
	for each $m \in \N$. Further, suppose that, as $n \to \infty$,
	\begin{equation*}
		\int_{\X^m} f \, \dd J_{\eta_n, B, m} 
		\longrightarrow \int_{\X^m} f \, \dd J_{\eta, B, m}
	\end{equation*}
	for all measurable and bounded functions $f : \X^m \to [0, \infty)$ and all $m \in \N$. Then $\E \big[ F(\eta_n) \big] \to \E \big[ F(\eta) \big]$ as $n \to \infty$ for all measurable, $B$-local, and bounded maps $F : \NN \to [0, \infty)$.
\end{lemma}
\begin{prf}
	Let $F : \NN \to [0, \infty)$ be measurable, $B$-local, and bounded by a constant $C \geq 0$. For now, assume that $F(\mathbf{0}) = 0$. For each $m \in \N$, define the function $f_m : \X^m \to [0, \infty)$ as $f_m(x_1, \dots, x_m) = F(\delta_{x_1} + \dotso + \delta_{x_m})$. These are measurable and symmetric functions which are bounded by $C$. Therefore,
	\begin{equation*}
		\sup_{n \in \N} \int_{\X^m} f_m(x_1, \dots, x_m) \, \dd J_{\eta_n, B, m}(x_1, \dots, x_m)
		\leq C \cdot c_B(m)
	\end{equation*}
	for all $m \in \N$ and the right hand side is summable over $m$. Thus, we can use Lemma \ref{LEMMA expec. repres. via Janossy measures} (minding $F(\mathbf{0}) = 0$), dominated convergence, and the assumption, to conclude that
	\begin{align*}
		\lim_{n \to \infty} \E \big[ F(\eta_n) \big]
		&= \lim_{n \to \infty} \sum_{m = 1}^\infty \int_{\X^m} f_m(x_1, \dots, x_m) \, \dd J_{\eta_n, B, m}(x_1, \dots, x_m) \\
		&= \sum_{m = 1}^\infty \int_{\X^m} f_m(x_1, \dots, x_m) \, \dd J_{\eta, B, m}(x_1, \dots, x_m) \\
		&= \E \big[ F(\eta) \big] .
	\end{align*}
	In case $F(\mathbf{0}) > 0$, the assertion follows by applying the first part of the proof to the decomposition
	\begin{equation*}
		F(\mu)
		= \max\big\{ F(\mu) - F(\mathbf{0}), \, 0 \big\} - \max\big\{ F(\mathbf{0}) - F(\mu), \, 0 \big\} + F(\mathbf{0}) .
	\end{equation*}
\end{prf}
\begin{remark} \label{RMK on the connection between local conv. and conv. of Janossy measures}
	If the map $c_B$ in Lemma \ref{LEMMA conv. of Janossy meas. implies local conv.} satisfies
	\begin{equation*}
		\sum_{m = 1}^\infty m \cdot c_B(m) < \infty,
	\end{equation*}
	the result holds for all measurable, $B$-local, and tame maps $F$. The stronger assumption becomes necessary since the maps $f_m$ defined in the proof only satisfy $f_m(\cdot) \leq C^\prime (1 + m)$ if $F$ is tame with constant $C^\prime$. Note that the expectations are finite for such maps by the assumption on $c_B$ since Lemma \ref{LEMMA expec. repres. via Janossy measures} implies
	\begin{equation*}
		\E\big[ F(\eta) \big]
		= \E\big[ F(\eta_B) \big]
		\leq C^\prime + \sum_{m = 1}^\infty \int_{\X^m} C^\prime (1 + m) \, \dd J_{\eta, B, m}(x_1, \dots, x_m)
		\leq C^\prime + C^\prime \sum_{m = 1}^\infty (1 + m) \cdot c_B(m) ,
	\end{equation*}
	and similarly for $\eta_n$ ($n \in \N$). If the boundedness- and convergence assumption on the Janossy measures holds for all $B \in \mathcal{X}_b$ (with the slightly stronger assumption on $c_B$ indicated above), the expectations converge for all local and tame $F$, so $\eta_n$ converges locally to $\eta$.
\end{remark}

\begin{example}[\bf Local convergence of Poisson processes] \label{EXA local conv. of Poisson proc.}
	Let $D_1, D_2, \dotso \in \mathcal{X}$ with $D_1 \subset D_2 \subset \dotso$ and $\bigcup_{n = 1}^\infty D_n = \X$. For each $n \in \N$ let $\Phi_n$ be a Poisson process in $\X$ with intensity measure $\lambda_{D_n}$ and let $\Phi$ be a Poisson process in $\X$ with intensity measure $\lambda$. Then $\Phi_n \stackrel{loc}{\longrightarrow} \Phi$ as $n \to \infty$. This follows readily from Appendix \ref{Appendix facts about Jan. and Cor. meas.} and Lemma \ref{LEMMA conv. of Janossy meas. implies local conv.}.
\end{example}

The converse of Lemma \ref{LEMMA conv. of Janossy meas. implies local conv.} is also true (without any integrability assumption).
\begin{lemma} \label{LEMMA local conv. implies conv. of Janossy meas.}
	Let $\eta, \eta_1, \eta_2, \ldots$ be point processes in $\X$ and fix a set $B \in \mathcal{X}_b$. Suppose that, as $n \to \infty$, $\E \big[ F(\eta_n) \big] \to \E \big[ F(\eta) \big]$ for all measurable, $B$-local, and bounded functions $F : \NN \to [0, \infty)$. Then, as $n \to \infty$,
	\begin{equation*}
		\int_{\X^m} f \, \dd J_{\eta_n, B, m} 
		\longrightarrow \int_{\X^m} f \, \dd J_{\eta, B, m}
	\end{equation*}
	for all measurable and bounded functions $f : \X^m \to [0, \infty)$ and all $m \in \N$.
\end{lemma}
\begin{prf}
	Fix $m \in \N$ and let $f : \X^m \to [0, \infty)$ be measurable and bounded by a constant $C \geq 0$. Define
	\begin{equation*}
		F(\mu) 
		= \frac{1}{m!} \cdot \mathds{1}\big\{ \mu(B) = m \big\} \int_{B^m} f(x_1, \dots, x_m) \, \dd\mu^{(m)}(x_1, \dots, x_m) .
	\end{equation*}
	The map $F$ is measurable by Lemma \ref{LEMMA measurability of factorial-measure-integrals}, it is $B$-local since $(\mu^{(m)})_{B^m} = \mu_B^{(m)}$, and bounded as
	\begin{equation*}
		F(\mu) 
		\leq \frac{C}{m!} \cdot \mathds{1}\big\{ \mu(B) = m \big\} \cdot \mu^{(m)}(B^m)
		\leq C , \quad \mu \in \NN,
	\end{equation*}
	using that $\mu \in \NN$ with $\mu(B) = m$ satisfies $\mu^{(m)}(B^m) = m!$ (see Proposition \ref{PROP properties of factorial measures due to Last/Penrose}). Thus, by assumption and the definition of the Janossy measures, we get
	\begin{equation*}
		\int_{\X^m} f \, \dd J_{\eta_n, B, m}
		= \frac{1}{m!} \, \E\bigg[ \mathds{1}\big\{ \eta_n(B) = m \big\} \int_{\X^m} f \, \dd(\eta_n)_B^{(m)} \bigg]
		= \E \big[ F(\eta_n) \big]
		\longrightarrow \E \big[ F(\eta) \big]
		= \int_{\X^m} f \, \dd J_{\eta, B, m} ,
	\end{equation*}
	as $n \to \infty$.
\end{prf}

\vspace{3mm}

We now state a result which implies the well-known connection between local convergence and weak* convergence of correlation functions. This convergence of the correlation functions is equivalent to the convergence of the Janossy measures in the previous lemmata, but only under a stronger assumption, namely a version of Ruelle's condition.
\begin{definition}[\bf Ruelle's condition] \label{DEF Ruelle's condition}
	A point process $\eta$ in $\X$ is said to satisfy \textit{Ruelle's condition} if
	\begin{equation*}
		\alpha_{\eta, m}(\cdot)
		\leq (\vartheta \lambda)^m(\cdot)
	\end{equation*}
	for each $m \in \N$ and some measurable function $\vartheta : \X \to [0, \infty)$. Here we denote by $\vartheta \lambda$ the measure on $(\X, \mathcal{X})$ with $\lambda$-density $\vartheta$.
\end{definition}
If a sequence $(\eta_n)_{n \in \N}$ of point processes satisfies the classical version of Ruelle's condition, where $\vartheta \equiv c$ with the same constant $c \geq 0$, the corresponding correlation functions exist and give a bounded sequence in $L^\infty$, and weak* convergence in $L^\infty$ is precisely what occurs in the following lemma if one observes that the factorial moment measures are given through the corresponding $\lambda^m$-densities, that is, the correlation functions. Notice that our new existence proof for Gibbs processes works via Lemma \ref{LEMMA conv. of Janossy meas. implies local conv.} and allows for a way around the $L^1$ test functions.

We call a function $\vartheta : \X \to [- \infty, \infty]$ \textit{locally $\lambda$-integrable} if it is integrable over bounded sets.

\begin{lemma} \label{LEMMA equiv. of local conv. and conv. of correl. functions}
	Let $\eta, \eta_1, \eta_2, \ldots$ be point processes in $\X$. Assume that there exists a measurable, locally $\lambda$-integrable map $\vartheta : \X \to [0, \infty)$ such that $\alpha_{\eta_n, m} \leq (\vartheta \lambda)^m$ for all $m \in \N$ and $n \in \N$. Then
	\begin{equation*}
		\int_{\X^m} f \, \dd J_{\eta_n, B, m} 
		\longrightarrow \int_{\X^m} f \, \dd J_{\eta, B, m} , \quad \text{as } n \to \infty,
	\end{equation*}
	for all measurable and bounded functions $f : \X^m \to [0, \infty)$, all $B \in \mathcal{X}_b$, and all $m \in \N$, if, and only if,
	\begin{equation*}
		\int_{\X^m} g \, \dd \alpha_{\eta_n, m}
		\longrightarrow \int_{\X^m} g \, \dd \alpha_{\eta, m} , \quad \text{as } n \to \infty,
	\end{equation*}
	for all $g \in L^1\big( \X^m, (\vartheta \lambda)^m \big)$ and every $m \in \N$.
\end{lemma}

Up to minor technical details, the necessity part of the lemma follows from Theorem \ref{THM repres. of factorial mom. meas. via Janossy meas.} and an adaptation of the proof of Theorem 2.51 from \cite{J:2017:script}. Sufficiency follows from Theorem \ref{THM repres. of Janossy meas. via factorial mom. meas.}. Lemmata \ref{LEMMA equiv. of local conv. and conv. of correl. functions} and \ref{LEMMA conv. of Janossy meas. implies local conv.} (together with Remark \ref{RMK on the connection between local conv. and conv. of Janossy measures}) imply that if the factorial moment measures (or correlation functions) converge in the weak sense of Lemma \ref{LEMMA equiv. of local conv. and conv. of correl. functions}, then $\eta_n$ converges locally to $\eta$. The converse statement can be formulated via Lemma \ref{LEMMA local conv. implies conv. of Janossy meas.}. Of course, if $\vartheta \equiv c$ for some $c \geq 0$, the functions $g$ in the lemma are from $L^1(\X^m, \lambda^m)$.

\begin{remark}[\bf A note on previous existence proofs for Gibbs processes] \label{RMK old exist. proof for Gibbs proc.}
	A rather modern version of an existence proof, which we refine in this manuscript, uses the following observation about the local convergence topology. Assume that $\eta_1, \eta_2, \ldots$ are point processes in a complete separable metric space $\X$ which satisfy Ruelle's condition for some universal constant $c \geq 0$. Then the corresponding correlation functions exist and $\rho_{n, m} = \rho_{\eta_n, m} \in L^\infty(\X^m, \lambda^m) = L^1(\X^m, \lambda^m)^\prime$. With the Banach-Alaoglu theorem from functional analysis and a diagonal sequence construction, it is possible to extract a subsequence $\{ n_k : k \in \N \} \subset \N$ such that
	\begin{equation*}
		\lim_{k \to \infty} \int_{\X^m} g \, \rho_{n_k, m} \, \dd \lambda^m 
		= \int_{\X^m} g \, \rho_{m} \, \dd \lambda^m 
	\end{equation*}
	for all $g \in L^1(\X^m, \lambda^m)$, each $m \in \N$, and some functions $\rho_m \in L^\infty(\X^m, \lambda^m)$. If one can show that the functions $\rho_m$ are the correlation functions of some point process $\eta$ then Lemmata \ref{LEMMA equiv. of local conv. and conv. of correl. functions} and \ref{LEMMA conv. of Janossy meas. implies local conv.} imply that $\eta_{n_k} \stackrel{loc}{\longrightarrow} \eta$ as $k \to \infty$. In order to show the existence of an infinite Gibbs process, the idea is to start with a suitable sequence of finite Gibbs processes, to guarantee that they satisfy Ruelle's bound, to construct the local limit as above, and to prove that the limit is itself Gibbs. 
	
	To obtain the limit process $\eta$ in the construction above, it is used that a family $\{ \rho_m : m \in \N \}$ of symmetric functions $\rho_m :\X^m \to [0, \infty)$, which satisfy the Ruelle condition, are the correlation functions of some point process if, and only if, for all $B \in \mathcal{X}_b$, all $m \in \N_0$, and $\lambda^m$-a.e.\ $(x_1, \dots, x_m) \in B^m$,
	\begin{equation*}
		\sum_{k = 0}^\infty \frac{(-1)^k}{k!} \int_{B^k} \rho_{m + k}(x_1, \dots, x_{m + k}) \, \dd \lambda^{k}(x_{m + 1}, \dots, x_{m + k}) 
		\geq 0 .
	\end{equation*}
	A comprehensive proof for locally stable energy functions is given in the lectures notes by \cite{J:2017:script}. There, the whole construction, including the proof of a less general version of Theorem \ref{THM repres. of Janossy meas. via factorial mom. meas.}, are given in terms of the so-called $K$-transform, as (implicitly) introduced by \cite{L:1973} and (explicitly) used by \cite{KK:2002}. These proofs based on the $K$-transform are technical and the whole construction relies on Ruelle's bound. Requiring a Ruelle type condition to assure that the correlation functions are a bounded sequence in $L^\infty$ is somewhat unnatural. In aiming at using Lemma \ref{LEMMA conv. of Janossy meas. implies local conv.}, we do not have to make sure the correlation functions are bounded, but can neatly employ that Janossy densities naturally form a bounded sequence in $L^1$. More precisely, for any sequence $(\eta_n)_{n \in \N}$ of point processes in $\X$, any $m \in \N$, and all $B \in \mathcal{X}_b$, we have 
	\begin{equation*}
		\sup_{n \in \N} \lVert j_{\eta_n, B, m} \rVert_{L^1(\X^m, \lambda^m)}
		= \sup_{n \in \N} J_{\eta_n, B, m}(B^m)
		= \sup_{n \in \N} \PP\big( \eta_n(B) = m \big)
		\leq 1,
	\end{equation*}
	given that the Janossy densities exist.
	
	As we show in our proof of Theorem \ref{THM repres. of Janossy meas. via factorial mom. meas.} and our other results, we focus on the point process theoretic perspective which leads to neat proofs. We use this perspective to provide a new result which allows for the extraction of locally convergent subsequences under weaker conditions.
\end{remark}

The following result forms the foundation of our existence proof for infinite Gibbs processes in Section \ref{SEC existence result Gibbs processes}. For the proof we need the Kolmogorov extension theorem for probability measures on $\NN$ as recalled in Appendix \ref{Appendix Kolmog. extension result on N}. In particular, to apply this theorem we have to restrict our attention to substandard Borel spaces (see Appendix \ref{Appendix Kolmog. extension result on N}). Notice that no claim is made about the uniqueness of the constructed subsequence and the limit process, we merely provide an existence result.
\begin{theorem} \label{THM existence of locally convergent subsequences}
	Let $(\X, \mathcal{X})$ be a substandard Borel space. Let $\eta_1, \eta_2, \ldots$ be point processes in $\X$ such that the Janossy densities $j_{n, B, m} = j_{\eta_n, B, m}$ corresponding to these processes exist and satisfy, for each $B \in \mathcal{X}_b$ and $m \in \N$, 
	\begin{equation*}
		\lim_{c \to \infty} \sup_{n \in \N} \int_{\X^m} j_{n, B, m}(x_1, \dots, x_m) \, \mathds{1}\big\{ j_{n, B, m}(x_1, \dots, x_m) \geq c \big\} \, \dd\lambda_B^m(x_1, \dots, x_m)
		= 0 .
	\end{equation*}
	Also assume that there exist maps $c_B : \N \to [0, \infty)$ (for each $B \in \mathcal{X}_b$) with $\sum_{m = 1}^\infty c_B(m) < \infty$ such that
	\begin{equation*}
		\sup_{n \in \N} \int_{\X^m} j_{n, B, m}(x_1, \dots, x_m) \, \dd\lambda_B^m(x_1, \dots, x_m)
		\leq c_B(m)
	\end{equation*}
	for each $m \in \N$. Then there exists a point process $\eta$ in $\X$ and a subsequence $(\eta_{n_k})_{k \in \N}$ such that $\E[ F(\eta_{n_k}) ] \to \E[ F(\eta) ]$ as $k \to \infty$, for all measurable, local, and bounded functions $F : \NN \to [0, \infty)$. If the maps $c_B$ are such that $\sum_{m = 1}^\infty m \cdot c_B(m) < \infty$, then even $\eta_{n_k} \stackrel{loc}{\longrightarrow} \eta$ as $k \to \infty$.
\end{theorem}
\begin{prf}
	We use the shorthand notations $J_{n, B, m} = J_{\eta_n, B, m}$ and $j_{n, B, m} = j_{\eta_n, B, m}$ for Janossy measures and densities. For the following construction, fix $B \in \mathcal{X}_b$. Since $J_{n, B, 0} = \PP\big( \eta_n(B) = 0 \big)$ is a bounded sequence in $[0, 1]$, there exists a strictly increasing map $r_0 : \N \to \N$, corresponding to the selection of a subsequence, such that
	\begin{equation*}
		\lim_{k \to \infty} J_{r_0(k), B, 0} = J_{B, 0}
	\end{equation*}
	for some $J_{B, 0} \in [0, 1]$. Iteratively applying the Dunford-Pettis lemma \citep[Corollary 4.7.19 of][]{B:2007} in the spaces $L^1(\X^m, \lambda_B^m)$ gives, for each $m \in \N$, a strictly increasing map $r_m : \N \to \N$ and a function $j_{B, m} \in L^1(\X^m, \lambda_B^m)$ which is set to $0$ outside of $B^m$ such that
	\begin{equation*}
		\lim_{k \to \infty} \int_{\X^m} f \, j_{r_0 \circ \dotso \circ r_m(k), B, m} \, \dd \lambda_B^m
		= \int_{\X^m} f \, j_{B, m} \, \dd \lambda_B^m
	\end{equation*}
	for all $f \in L^\infty(\X^m, \lambda_B^m)$. If we put $n_k = r_0 \circ \dotso \circ r_k(k)$, $k \in \N$, which corresponds to taking the diagonal sequence, then
	\begin{equation} \label{weak* convergence of Janossy densities}
		\lim_{k \to \infty} J_{n_k, B, 0} = J_{B, 0}
		\quad \text{and} \quad
		\lim_{k \to \infty} \int_{\X^m} f \, j_{n_k, B, m} \, \dd \lambda_B^m
		= \int_{\X^m} f \, j_{B, m} \, \dd \lambda_B^m
	\end{equation} 
	for all $f \in L^\infty(\X^m, \lambda_B^m)$ and each $m \in \N$. For each $m \in \N$, the limit function $j_{B, m}$ is non-negative $\lambda^m$-a.e. Moreover, the measures $J_{B, m}$ ($m \in \N$) defined as $J_{B, m}(D) = \int_{\X^m} \mathds{1}_D \cdot j_{B, m} \, \dd \lambda^m$, $D \in \mathcal{X}^{\otimes m}$, are symmetric. Indeed, for sets $D_1, \dots, D_m \in \mathcal{X}$ and any permutation $\tau$ of $\{ 1, \dots, m \}$, equation \eqref{weak* convergence of Janossy densities} gives
	\begin{align*}
		J_{B, m}\big( D_{\tau(1)} \times \dotso \times D_{\tau(m)} \big)
		= \lim_{k \to \infty} J_{n_k, B, m}\big( D_{\tau(1)} \times \dotso \times D_{\tau(m)} \big)
		&= \lim_{k \to \infty} J_{n_k, B, m}\big( D_{1} \times \dotso \times D_{m} \big) \\
		&= J_{B, m}\big( D_{1} \times \dotso \times D_{m} \big) .
	\end{align*}
	Notice that the so constructed subsequence depends on the chosen set $B$.
	
	We now proceed to apply the above construction to the sets $B_\ell$. Applied to $B_1$, the previous arguments provide a subsequence $\{ n_k^1 : k \in \N \} \subset \N$ as well as $J_{B_1, 0} \in [0, 1]$ and ($\lambda^m$-a.e.) non-negative functions $j_{B_1, m} \in L^1(\X^m, \lambda^m)$, $m \in \N$, which vanish outside of $B_1^m$, such that
	\begin{equation*}
		\lim_{k \to \infty} J_{n_k^1, B_1, 0} = J_{B_1, 0}
		\quad \text{and} \quad
		\lim_{k \to \infty} \int_{\X^m} f \, j_{n_k^1, B_1, m} \, \dd \lambda^m
		= \int_{\X^m} f \, j_{B_1, m} \, \dd \lambda^m
	\end{equation*} 
	for all $f \in L^\infty(\X^m, \lambda^m)$ and each $m \in \N$. Iteratively applying this scheme gives (in the $\ell$-th step) a subsequence 
	\begin{equation*}
		\{ n_k^\ell : k \in \N \} \subset \{ n_k^{\ell - 1} : k \in \N \} \subset \dotso \subset \{ n_k^1 : k \in \N \}
	\end{equation*}
	as well as $J_{B_\ell, 0} \in [0, 1]$ and ($\lambda^m$-a.e.) non-negative functions $j_{B_\ell, m} \in L^1(\X^m, \lambda^m)$, $m \in \N$, which vanish outside of $B_\ell^m$, such that
	\begin{equation*}
		\lim_{k \to \infty} J_{n_k^\ell, B_\ell, 0} = J_{B_\ell, 0}
		\quad \text{and} \quad
		\lim_{k \to \infty} \int_{\X^m} f \, j_{n_k^\ell, B_\ell, m} \, \dd \lambda^m
		= \int_{\X^m} f \, j_{B_\ell, m} \, \dd \lambda^m
	\end{equation*} 
	for all $f \in L^\infty(\X^m, \lambda^m)$ and each $m \in \N$. Thus, choosing the diagonal sequence $n_k = n_k^k$, we have
	\begin{equation} \label{convergence of janossy densities in exist. result}
		\lim_{k \to \infty} J_{n_k, B_\ell, 0} = J_{B_\ell, 0}
		\quad \text{and} \quad
		\lim_{k \to \infty} \int_{\X^m} f \, j_{n_k, B_\ell, m} \, \dd \lambda^m
		= \int_{\X^m} f \, j_{B_\ell, m} \, \dd \lambda^m
	\end{equation} 
	for all $f \in L^\infty(\X^m, \lambda^m)$, each $m \in \N$, and every $\ell \in \N$. On each of the substandard Borel spaces $(\NN, \mathcal{N}_{B_\ell})$, $\ell \in \N$, we define the measures
	\begin{equation*}
		\mathrm{P}_\ell^{(n)}(A)
		= \mathds{1}_A(\mathbf{0}) \cdot J_{n, B_\ell, 0} + \sum_{m = 1}^\infty \int_{B_\ell^m} \mathds{1}_A\Big( \sum_{i = 1}^m \delta_{x_i} \Big) \, j_{n, B_\ell, m}(x_1, \dots, x_m) \, \dd \lambda^m(x_1, \dots, x_m), \quad A \in \mathcal{N}_{B_\ell} ,
	\end{equation*}
	$n \in \N$, and
	\begin{equation*}
		\mathrm{P}_\ell(A)
		= \mathds{1}_A(\mathbf{0}) \cdot J_{B_\ell, 0} + \sum_{m = 1}^\infty \int_{B_\ell^m} \mathds{1}_A\Big( \sum_{i = 1}^m \delta_{x_i} \Big) \, j_{B_\ell, m}(x_1, \dots, x_m) \, \dd \lambda^m(x_1, \dots, x_m), \quad A \in \mathcal{N}_{B_\ell} .
	\end{equation*}
	By assumption, we have
	\begin{equation*}
		\sup_{k \in \N} \int_{B_\ell^m} \mathds{1}_A\Big( \sum_{i = 1}^m \delta_{x_i} \Big) \, j_{n_k, B_\ell, m}(x_1, \dots, x_m) \, \dd \lambda^m(x_1, \dots, x_m)
		\leq c_{B_\ell}(m) ,
	\end{equation*}
	where the right hand side constitutes an integrable bound with respect to summation over $m$. Hence, we can apply dominated convergence and the limit results from equation \eqref{convergence of janossy densities in exist. result} to conclude that
	\begin{equation*}
		\lim_{k \to \infty} \mathrm{P}_\ell^{(n_k)}(A)
		= \mathrm{P}_\ell(A)
	\end{equation*}
	for all $A \in \mathcal{N}_{B_\ell}$ and each $\ell \in \N$. Notice that Lemma \ref{LEMMA expec. repres. via Janossy measures} implies that, for $A \in \mathcal{N}_{B_\ell}$, we have
	\begin{equation*}
		\mathrm{P}_\ell^{(n)}(A)
		= \PP(\eta_n \in A) .
	\end{equation*}
	It follows that $\mathrm{P}_\ell^{(n)}(\NN) = 1$ and therefore $\mathrm{P}_\ell(\NN) = \lim_{k \to \infty} \mathrm{P}_\ell^{(n_k)}(\NN) = 1$, so $\mathrm{P}_\ell, \mathrm{P}_\ell^{(1)}, \mathrm{P}_\ell^{(2)}, \ldots$ are probability measures on $(\NN, \mathcal{N}_{B_\ell})$ for each $\ell \in \N$. Moreover, if we take indices $i < \ell$ and any set $A \in \mathcal{N}_{B_i} \subset \mathcal{N}_{B_\ell}$, we get
	\begin{equation*}
		\mathrm{P}_\ell(A)
		= \lim_{k \to \infty} \mathrm{P}_\ell^{(n_k)}(A)
		= \lim_{k \to \infty} \PP(\eta_{n_k} \in A)
		= \lim_{k \to \infty} \mathrm{P}_i^{(n_k)}(A)
		= \mathrm{P}_i(A) .
	\end{equation*}
	Proposition \ref{PROP Kolmogorov extension theorem, countable version}, implies that there exists a probability measure $\mathrm{P}$ on $(\NN, \mathcal{N})$ such that $\mathrm{P}(A) = \mathrm{P}_\ell(A)$ for all $A \in \mathcal{N}_{B_\ell}$ and $\ell \in \N$. Thus, if we let $\eta$ be a point process in $\X$ with distribution $\mathrm{P}$, then by construction
	\begin{equation*}
		\PP(\eta \in A)
		= \mathrm{P}(A)
		= \mathrm{P}_\ell(A)
		= \mathds{1}_A(\mathbf{0}) \cdot J_{B_\ell, 0} + \sum_{m = 1}^\infty \int_{B_\ell^m} \mathds{1}_A\Big( \sum_{i = 1}^m \delta_{x_i} \Big) \, j_{B_\ell, m}(x_1, \dots, x_m) \, \dd \lambda^m(x_1, \dots, x_m)
	\end{equation*}
	for each $A \in \mathcal{N}_{B_\ell}$ and $\ell \in \N$, and hence Lemma \ref{LEMMA characterization of Janossy measures} implies that $(J_{B_\ell, m})_{m \in \N_0}$ are the Janossy measures of $\eta$ restricted to $B_\ell$, for each $\ell \in \N$. The limit relations in \eqref{convergence of janossy densities in exist. result} and Lemma \ref{LEMMA conv. of Janossy meas. implies local conv.} yield
	\begin{equation*}
		\lim_{k \to \infty} \E \big[ F(\eta_{n_k}) \big]
		= \E \big[ F(\eta) \big]
	\end{equation*}
	for all measurable, local, and bounded maps $F : \NN \to [0, \infty)$, where we use that any local function is $B_\ell$-local for some $\ell \in \N$. The additional claim concerning local convergence follows readily from Remark \ref{RMK on the connection between local conv. and conv. of Janossy measures}.
\end{prf}

\vspace{3mm}

It is clear from the proof that the assumptions on the Janossy densities in the theorem need only be satisfied on the sets $B_\ell$.

\begin{remark} \label{RMK conditions on the Janossy density for the conv. result}
	We now discuss a condition which is sufficient to ensure the uniform integrability and summability conditions of Theorem \ref{THM existence of locally convergent subsequences}. As in the theorem, let $\eta_1, \eta_2, \ldots$ be point processes in $\X$ with Janossy densities $j_{n, B, m}$. Assume that there exists a measurable and locally $\lambda$-integrable function $\vartheta : \X \to [0, \infty)$ such that
	\begin{equation*}
		\sup_{n \in \N} j_{n, B, m}(x_1, \dots, x_m)
		\leq \frac{\vartheta(x_1) \cdot \dotso \cdot \vartheta(x_m)}{m!}
	\end{equation*}
	for $\lambda^m$-a.e.\ $(x_1, \dots, x_m) \in \X^m$, each $m \in \N$, and all $B \in \mathcal{X}_b$. Then, for all $B \in \mathcal{X}_b$ and $m \in \N$,		
	\begin{align*}
		&\lim_{c \to \infty} \sup_{n \in \N} \int_{\X^m} j_{n, B, m}(x_1, \dots, x_m) \, \mathds{1}\big\{ j_{n, B, m}(x_1, \dots, x_m) \geq c \big\} \, \dd\lambda_B^m(x_1, \dots, x_m) \\
		&\quad \leq \frac{1}{m!} \lim_{c \to \infty} \int_{B^m} \vartheta(x_1) \cdot \dotso \cdot \vartheta(x_m) \cdot \mathds{1}\big\{ \vartheta(x_1) \cdot \dotso \cdot \vartheta(x_m) \geq c \cdot m! \big\} \, \dd\lambda^m(x_1, \dots, x_m)
		= 0
	\end{align*}
	as well as
	\begin{equation*}
		\sup_{n \in \N} \int_{\X^m} j_{n, B, m}(x_1, \dots, x_m) \, \dd\lambda_B^m(x_1, \dots, x_m)
		\leq \frac{1}{m!} \bigg( \int_B \vartheta(x) \, \dd\lambda(x) \bigg)^m ,
	\end{equation*}
	where $c_B(m) = \frac{1}{m!} \big( \int_B \vartheta(x) \, \dd\lambda(x) \big)^m$ meets
	\begin{equation*}
		\sum_{m = 1}^\infty m \cdot c_B(m)
		= \int_B \vartheta(x) \, \dd\lambda(x) \cdot \exp\bigg( \int_B \vartheta(x) \, \dd\lambda(x) \bigg)
		< \infty .
	\end{equation*}
	By Theorem \ref{THM repres. of factorial mom. meas. via Janossy meas.} and Corollary \ref{COR repres. of correlation func. via Janossy densities} this bound on the Janossy densities leads to ($B$-dependent) bounds on the factorial moment measures and correlation functions, so the previous condition can be interpreted as local versions of Ruelle's condition from Definition \ref{DEF Ruelle's condition}.
\end{remark}

\section{An existence result for Gibbs point processes in general spaces}
\label{SEC existence result Gibbs processes}

We first state and prove the abstract result and provide basic discussions on the assumptions. Afterward we turn to initial and immediate examples.

\subsection{The abstract result}
\label{SUBSEC abstract existence result for Gibbs point processe}

Consider the setting from the beginning of Section \ref{SEC Introduction}. In order to construct a (possibly infinite) Gibbs process with PI $\kappa$ and reference measure $\lambda$, we let $\xi_1, \xi_2, \ldots$ be a sequence of finite Gibbs processes with PIs  $\kappa^{(B_n, \mathbf{0})}$. These processes are given explicitly by Lemma \ref{LEMMA characterization of finite Gibbs process via kappa_m}. It is our goal to use Theorem \ref{THM existence of locally convergent subsequences} to extract from $(\xi_n)_{n \in \N}$ a locally convergent subsequence and a limit process $\eta$. We then want to prove that $\eta$ is a Gibbs process with PI $\kappa$. For the first step we need to ensure that the assumptions on the Janossy measures in Theorem \ref{THM existence of locally convergent subsequences} are satisfied. In this we focus on the assumption given in Remark \ref{RMK conditions on the Janossy density for the conv. result}. In and of itself this is an assumption on the whole construction of the finite Gibbs processes $\xi_n$ and not a mere condition for $\kappa$ and $\lambda$. Such explicit assumptions will follow later, but the generality of the following theorem will be useful.

\begin{theorem} \label{THM abstract existence result}
	Let $(\X, \mathcal{X})$ be a substandard Borel space with localizing structure $B_1 \subset B_2 \subset \dotso$ and let $\lambda$ be a locally finite measure on $\X$. Let $\kappa : \X \times \NN \to [0, \infty)$ be a measurable map which satisfies the cocycle relation \eqref{cocycle assumption} and is such that
	\begin{equation*}
		Z_{B_n}(\mathbf{0}) < \infty, \quad n \in \N .
	\end{equation*}
	Moreover, suppose that, for $\lambda$-a.e.\ $x \in \X$,
	\begin{equation*}
		\kappa(x, \mu)
		\leq \tilde{\vartheta}(x) \cdot c^{\mu(\X)}, \quad \mu \in \NN,
	\end{equation*}
	for a constant $c \geq 0$ and a measurable, locally $\lambda$-integrable map $\tilde{\vartheta} : \X \to [0, \infty)$. Let $\xi_n$ be a finite Gibbs process with PI $\kappa^{(B_n, \mathbf{0})}$, for each $n \in \N$, and assume that
	\begin{equation*}
		\sup_{n \in \N} j_{\xi_n, B, m}(x_1, \dots, x_m)
		\leq \frac{\vartheta(x_1) \cdot \dotso \cdot \vartheta(x_m)}{m!}
	\end{equation*}
	for $\lambda^m$-a.e.\ $(x_1, \dots, x_m) \in \X^m$, all $m \in \N$, every $B \in \mathcal{X}_b$, and some measurable, locally $\lambda$-integrable map $\vartheta : \X \to [0, \infty)$. Denote by $\eta$ any one of the limit processes obtainable from Theorem \ref{THM existence of locally convergent subsequences} and assume that, for all $B \in \mathcal{X}^*_b$, 
	\begin{equation*}
		\limsup_{\ell \to \infty} \int_B \E\big| \kappa(x, \eta_{B_\ell}) - \kappa(x, \eta) \big| \, \dd\lambda(x)
		= 0
	\end{equation*}
	as well as
	\begin{equation*}
		\limsup_{\ell \to \infty} \, \sup_{k \in \N} \int_B \E\big| \kappa\big( x, (\xi_{n_k})_{B_\ell} \big) - \kappa(x, \xi_{n_k}) \big| \, \dd\lambda(x)
		= 0 ,
	\end{equation*}
	where $(\xi_{n_k})_{k \in \N}$ is the subsequence of $(\xi_n)_{n \in \N}$ which converges locally to $\eta$ and where $\mathcal{X}^*_b \subset \mathcal{X}_b$ is a $\pi$-system which contains a nested sequence of sets that exhaust $\X$ and is such that
	\begin{equation*}
		\sigma\big( \{ B \times A : B \in \mathcal{X}^*_b, \, A \in \mathcal{Z} \} \big)
		= \mathcal{X} \otimes \mathcal{N} ,
	\end{equation*}
	with $\mathcal{Z}$ denoting the local events from Definition \ref{DEF local events and functions}.Then $\eta$ is a Gibbs process with PI $\kappa$.
\end{theorem}
\begin{prf}
	First of all, recall that by Lemma \ref{LEMMA characterization of finite Gibbs process via kappa_m} the Gibbs processes $\xi_n$ exist as $Z_{B_n}(\mathbf{0}) < \infty$ for $n \in \N$. The bound on the Janossy densities covers the assumptions of Theorem \ref{THM existence of locally convergent subsequences}, so we have $\xi_{n_k} \stackrel{loc}{\longrightarrow} \eta$ as $k \to \infty$, where the subsequence and the limit process $\eta$ are as in the statement of the theorem. It remains to prove that $\eta$ is a Gibbs process with PI $\kappa$.
	
	Notice that the Janossy densities of $\eta$ satisfy the same bound as those of the processes $\xi_n$. Indeed, by Lemma \ref{LEMMA local conv. implies conv. of Janossy meas.} we have, for each $m \in \N$, every $B \in \mathcal{X}_b$, and any measurable and bounded map $f : \X^m \to [0, \infty)$,
	\begin{equation*}
		\int_{\X^m} f \, j_{\eta, B, m} \, \dd\lambda^m
		= \limsup_{k \to \infty} \int_{\X^m} f \, j_{\xi_{n_k}, B, m} \, \dd\lambda^m
		\leq \frac{1}{m!} \int_{B^m} f(x_1, \dots, x_m) \cdot \vartheta(x_1) \cdot \dotso \cdot \vartheta(x_m) \, \dd\lambda^m(x_1, \dots, x_m) ,
	\end{equation*}
	so $j_{\eta, B, m}(x_1, \dots, x_m) \leq \frac{\vartheta(x_1) \cdot \dotso \cdot \vartheta(x_m)}{m!}$ for $\lambda^m$-a.e.\ $(x_1, \dots, x_m) \in \X^m$.
	
	Fix $B \in \mathcal{X}^*_b$ and $A \in \mathcal{Z}$, and let $C \in \mathcal{X}_b$ be such that $A \in \mathcal{N}_C$. Define the following measurable maps $\NN \to [0, \infty]$,
	\begin{align*}
		F(\mu) &= \int_\X \mathds{1}_B(x) \, \mathds{1}_A(\mu) \, \dd \mu(x) , \\
		\tilde{F}(\mu) &= \int_\X \mathds{1}_B(x) \, \mathds{1}_A(\mu + \delta_x) \, \kappa(x, \mu) \, \dd\lambda(x) , \\
		\tilde{F}_\ell(\mu) &= \int_\X \mathds{1}_B(x) \, \mathds{1}_A(\mu + \delta_x) \, \kappa(x, \mu_{B_\ell}) \, \dd\lambda(x), \quad \ell \in \N.
	\end{align*}
	We collect in four steps the essential properties of these maps.
	\begin{itemize}
		\item[(i)] The function $F$ is $(B \cup C)$-local and tame, as $F(\mu) \leq \mu(B)$. Thus, the local convergence applies to $F$, so
		\begin{equation*}
			\lim_{k \to \infty} \E\big[ F(\xi_{n_k}) \big]
			= \E \big[ F(\eta) \big] ,
		\end{equation*}
		where these expectations are bounded by
		\begin{align*}
			\sup_{n \in \N} \E \big[ F(\xi_{n}) \big]
			\leq \sup_{n \in \N} \E\big[ \xi_{n}(B) \big]
			&= \sup_{n \in \N} \sum_{m = 1}^\infty \int_{\X^m} \sum_{i = 1}^m \delta_{x_i}(B) \cdot j_{\xi_{n}, B, m}(x_1, \dots, x_m) \, \dd\lambda^m(x_1, \dots, x_m) \\
			&\leq \int_B \vartheta(x) \, \dd\lambda(x) \cdot \exp\bigg( \int_B \vartheta(x) \, \dd\lambda(x) \bigg) ,
		\end{align*}
		and similarly for $\E[F(\eta)]$ (see also Remark \ref{RMK on the connection between local conv. and conv. of Janossy measures}).
		
		\item[(ii)] Using that, for $\lambda$-a.e.\ $x \in \X$ and all $\ell \in \N$,
		\begin{equation} \label{bound on kappa in the existence proof}
			\E\big[ \kappa(x, \eta_{B_\ell}) \big]
			\leq \tilde{\vartheta}(x) \cdot \E\big[ c^{\eta(B_\ell)} \big]
			= \tilde{\vartheta}(x) \sum_{m = 0}^\infty c^m \cdot J_{\eta, B_\ell, m}(B_\ell^m)
			\leq \tilde{\vartheta}(x) \cdot \exp\bigg( c \int_{B_\ell} \vartheta(y) \, \dd\lambda(y) \bigg) ,
		\end{equation}
		we obtain
		\begin{equation*}
			\E\big[ \tilde{F}_\ell(\eta) \big]
			\leq \int_B \E\big[ \kappa(x, \eta_{B_\ell}) \big] \, \dd\lambda(x)
			\leq \int_B \tilde{\vartheta}(x) \, \dd\lambda(x) \cdot \exp\bigg( c \int_{B_\ell} \vartheta(x) \, \dd\lambda(x) \bigg),
		\end{equation*}
		and the very same term bounds $\sup_{n \in \N} \E\big[ \tilde{F}_\ell(\xi_n) \big]$. By assumption, we have
		\begin{equation*}
			\limsup_{\ell \to \infty} \big| \E\big[ \tilde{F}(\eta) \big] - \E\big[ \tilde{F}_\ell(\eta) \big] \big|
			\leq \limsup_{\ell \to \infty} \int_B \E\big| \kappa(x, \eta_{B_\ell}) - \kappa(x, \eta) \big| \, \dd\lambda(x)
			= 0 ,
		\end{equation*}
		and
		\begin{equation*}
			\limsup_{\ell \to \infty} \, \sup_{k \in \N} \big| \E\big[ \tilde{F}(\xi_{n_k}) \big] - \E\big[ \tilde{F}_\ell(\xi_{n_k}) \big] \big|
			\leq \limsup_{\ell \to \infty} \, \sup_{k \in \N} \int_B \E\big| \kappa\big( x, (\xi_{n_k})_{B_\ell} \big) - \kappa(x, \xi_{n_k}) \big| \, \dd\lambda(x)
			= 0 .
		\end{equation*}
		
		\item[(iii)] We now show that, despite $\tilde{F}_\ell$ not being tame, we have $\lim_{k \to \infty} \E\big[ \tilde{F}_\ell(\xi_{n_k}) \big] = \E\big[ \tilde{F}_\ell(\eta) \big]$ for each $\ell \in \N$. To this end, fix $\ell\in \N$ and define the measurable maps $\NN \to [0, \infty)$,
		\begin{equation*}
			\tilde{F}_{\ell , j}(\mu)
			= \int_\X \mathds{1}_B(x) \, \mathds{1}_A(\mu + \delta_x) \, \kappa(x, \mu_{B_\ell}) \, \mathds{1}\big\{ \kappa(x, \mu_{B_\ell}) \leq j \big\} \, \dd\lambda(x), \quad j \in \N.
		\end{equation*}
		For $j \in \N$, the map $\tilde{F}_{\ell , j}$ is $(B \cup C \cup B_\ell)$-local and bounded by $j \cdot \lambda(B)$. Therefore, the local convergence applies, so
		\begin{equation*}
			\lim_{k \to \infty} \E\big[ \tilde{F}_{\ell, j}(\xi_{n_k}) \big] 
			= \E\big[ \tilde{F}_{\ell, j}(\eta) \big]
		\end{equation*}
		for every $j \in \N$. Using \eqref{bound on kappa in the existence proof} to justify the application of dominated convergence, we have
		\begin{equation*}
			\limsup_{j \to \infty} \big| \E\big[ \tilde{F}_{\ell, j}(\eta) \big] - \E\big[ \tilde{F}_\ell(\eta) \big] \big|
			\leq \limsup_{j \to \infty} \int_B \E\Big[ \kappa(x, \eta_{B_\ell}) \, \mathds{1}\big\{ \kappa(x, \eta_{B_\ell}) > j \big\} \Big] \, \dd\lambda(x)
			= 0 .
		\end{equation*}
		Moreover, observe that, by Lemma \ref{LEMMA expec. repres. via Janossy measures} and the bounds on $\kappa$ and the Janossy densities,
		\begin{align*}
			&\sup_{k \in \N} \big| \E\big[ \tilde{F}_{\ell, j}(\xi_{n_k}) \big] - \E\big[ \tilde{F}_\ell(\xi_{n_k}) \big] \big| \\
			&\quad \leq \sup_{k \in \N} \int_B \E\Big[ \kappa\big( x, (\xi_{n_k})_{B_\ell} \big) \, \mathds{1}\big\{ \kappa\big( x, (\xi_{n_k})_{B_\ell} \big) > j \big\} \Big] \, \dd\lambda(x) \\
			&\quad = \sup_{k \in \N} \bigg( \int_B \kappa(x, \mathbf{0}) \, \mathds{1}\big\{ \kappa(x, \mathbf{0}) > j \big\} \, \dd\lambda(x) \cdot \PP\big( \xi_{n_k}(B_\ell) = 0 \big) \\ 
			&\hspace{16mm} + \sum_{m = 1}^\infty \int_{\X^m} \int_B \kappa\Big( x, \sum_{i = 1}^m (\delta_{x_i})_{B_\ell} \Big) \, \mathds{1}\bigg\{ \kappa\Big( x, \sum_{i = 1}^m (\delta_{x_i})_{B_\ell} \Big) > j \bigg\} \, \dd\lambda(x) \, \dd J_{\xi_{n_k}, B_\ell, m}(x_1, \dots, x_m) \bigg) \\
			&\quad \leq \int_B \tilde{\vartheta}(x) \, \mathds{1}\big\{ \tilde{\vartheta}(x) > j \big\} \, \dd\lambda(x) + \sum_{m = 1}^\infty \frac{c^m}{m!} \bigg( \int_{B_\ell} \vartheta(x) \, \dd\lambda(x) \bigg)^m \int_B \tilde{\vartheta}(x) \, \mathds{1}\big\{ \tilde{\vartheta}(x) \cdot c^m > j \big\} \, \dd\lambda(x) ,
		\end{align*}
		and the right hand side converges to $0$ as $j \to \infty$, by dominated convergence. Now, let $\varepsilon > 0$. Choose $j_0 \in \N$ such that
		\begin{equation*}
			\big| \E\big[ \tilde{F}_{\ell, j_0}(\eta) \big] - \E\big[ \tilde{F}_\ell(\eta) \big] \big|
			< \frac{\varepsilon}{3}
			\qquad \text{and} \qquad
			\sup_{k \in \N} \big| \E\big[ \tilde{F}_{\ell, j_0}(\xi_{n_k}) \big] - \E\big[ \tilde{F}_\ell(\xi_{n_k}) \big] \big|
			< \frac{\varepsilon}{3} .
		\end{equation*}
		Choose $k_0 \in \N$ such that, for each $k \geq k_0$,
		\begin{equation*}
			\big| \E\big[ \tilde{F}_{\ell, j_0}(\xi_{n_k}) \big] - \E\big[ \tilde{F}_{\ell, j_0}(\eta) \big] \big|
			< \frac{\varepsilon}{3} .
		\end{equation*}
		Then, for each $k \geq k_0$, the triangle inequality yields $\big| \E\big[ \tilde{F}_{\ell}(\xi_{n_k}) \big] - \E\big[ \tilde{F}_{\ell}(\eta) \big] \big| < \varepsilon$.
		
		\item[(iv)] We now use (ii) and (iii) to show that $\lim_{k \to \infty} \E\big[ \tilde{F}(\xi_{n_k}) \big] = \E\big[ \tilde{F}(\eta) \big]$. Let $\varepsilon > 0$. By (ii), we can choose $\ell_0 \in \N$ such that
		\begin{equation*}
			\big| \E\big[ \tilde{F}(\eta) \big] - \E\big[ \tilde{F}_{\ell_0}(\eta) \big] \big|
			< \frac{\varepsilon}{3}
			\qquad \text{and} \qquad
			\sup_{k \in \N} \big| \E\big[ \tilde{F}(\xi_{n_k}) \big] - \E\big[ \tilde{F}_{\ell_0}(\xi_{n_k}) \big] \big|
			< \frac{\varepsilon}{3} .
		\end{equation*}
		By (iii), we can choose $k_0 \in \N$ such that, for each $k \geq k_0$,
		\begin{equation*}
			\big| \E\big[ \tilde{F}_{\ell_0}(\xi_{n_k}) \big] - \E\big[ \tilde{F}_{\ell_0}(\eta) \big] \big|
			< \frac{\varepsilon}{3} .
		\end{equation*}
		Then, for each $k \geq k_0$, the triangle inequality gives $\big| \E\big[ \tilde{F}(\xi_{n_k}) \big] - \E\big[ \tilde{F}(\eta) \big] \big| < \varepsilon$.
	\end{itemize}
	If $k$ is large enough so that $B \subset B_{n_k}$, then the GNZ equation for $\xi_{n_k}$ reads as
	\begin{equation*}
		\E\big[ F(\xi_{n_k}) \big]
		= \E \bigg[ \int_\X \mathds{1}_B(x) \, \mathds{1}_A(\xi_{n_k}) \, \dd \xi_{n_k}(x) \bigg]
		= \E \bigg[ \int_\X \mathds{1}_B(x) \, \mathds{1}_A(\xi_{n_k} + \delta_x) \, \kappa(x, \xi_{n_k}) \, \mathds{1}_{B_{n_k}}(x) \, \dd \lambda(x) \bigg]
		= \E\big[ \tilde{F}(\xi_{n_k}) \big] .
	\end{equation*}
	Thus, the convergence results from (i) and (iv) immediately give
	\begin{equation*}
		\E \bigg[ \int_\X \mathds{1}_B(x) \, \mathds{1}_A(\eta) \, \dd \eta(x) \bigg]
		= \E \bigg[ \int_\X \mathds{1}_B(x) \, \mathds{1}_A(\eta + \delta_x) \, \kappa(x, \eta) \, \dd \lambda(x) \bigg] .
	\end{equation*}
	Consequently, the GNZ equation holds for all functions $(x, \mu) \mapsto \mathds{1}_{B \times A}(x, \mu)$ with $B \in \mathcal{X}^*_b$ and $A \in \mathcal{Z}$. A literal copy of the final step in the proof of Lemma \ref{LEMMA mult. GNZ} (for $m = 0$), using the properties of $\mathcal{X}_b^*$, extends the equality to all functions $(x, \mu) \mapsto \mathds{1}_E(x, \mu)$, $E \in \mathcal{X} \otimes \mathcal{N}$, and monotone approximation allows for any measurable function $f : \X \times \NN \to [0, \infty]$. We conclude that $\eta$ is a Gibbs process with PI $\kappa$.
\end{prf}

\begin{remark} \label{RMK general remarks to the abstract existence result}
	In order to construct Gibbs processes with a boundary condition $\psi$, one has to apply Theorem \ref{THM abstract existence result} to $\tilde{\kappa} = \kappa(\,\cdot\, , \psi + \cdot\,)$. However, apart from the cocycle assumption, it is not a given that $\tilde{\kappa}$ inherits the necessary properties from $\kappa$. In the special case where $\kappa$ is locally stable, a condition discussed below, $\tilde{\kappa}$ inherits this property and most of the assumptions in Theorem \ref{THM abstract existence result} are satisfied.
		
	Note that any Gibbs process whose Janossy densities satisfy the bound in Theorem \ref{THM abstract existence result} has a locally finite intensity measure. Indeed, this we have shown in item (i) of the proof of Theorem \ref{THM abstract existence result}.
\end{remark}

One particular assumption on $\kappa$, which covers all bounds in Theorem \ref{THM abstract existence result}, is the local stability assumption which is frequent in stochastic geometry and spatial statistics as it is an essential assumption for many simulation algorithms for Gibbsian point processes, cf. \cite{MW:2004}.
\begin{definition}[\bf Local stability] \label{DEF local stability of kappa}
	A measurable map $\kappa : \X \times \NN \to [0, \infty)$ is called \textit{($\lambda$-)locally stable} if
	\begin{equation} \label{local stability of kappa}
		\sup_{\mu \in \NN} \kappa(x, \mu)
		\leq \vartheta(x)
	\end{equation}
	for $\lambda$-a.e.\ $x \in \X$ and some measurable, locally $\lambda$-integrable map $\vartheta : \X \to [0, \infty)$.
\end{definition}

\begin{remark} \label{RMK notes on the local stability assumption}
	Despite being most handy, local stability can be a major restriction. The setting of pair potentials constitutes an example where Theorem \ref{THM abstract existence result} can be used to show existence when $\kappa$ is not locally stable. However, the local stability assumption comes as a convenient one. Just recall that, by Lemma \ref{LEMMA Janossy densities of Gibbs processes}, we have
	\begin{equation*}
		j_{\xi_n, B, m}(x_1, \dots, x_m)
		= \frac{1}{m!} \, \E\Big[ \mathds{1}\big\{ \xi_n(B) = 0 \big\} \, \kappa_m(x_1, \dots, x_m, \xi_n) \Big] \, \mathds{1}_{(B \cap B_n)^m}(x_1, \dots, x_m)
	\end{equation*}
	where $\xi_n$ is a Gibbs process with PI $\kappa^{(B_n, \mathbf{0})}$ ($n \in \N$) as in Theorem \ref{THM abstract existence result}. These densities are certainly bounded in the desired manner if we assume that
	\begin{equation*}
		\sup_{\mu \in \NN_f} \kappa_m(x_1, \dots, x_m, \mu)
		\leq \vartheta(x_1) \cdot \dotso \cdot \vartheta(x_m)
	\end{equation*}
	for $\lambda^m$-a.e.\ $(x_1, \dots, x_m) \in \X^m$, all $m \in \N$, and some locally $\lambda$-integrable map $\vartheta : \X \to [0, \infty)$. This is virtually equivalent to $\kappa$ being locally stable. Aside from the bound on the Janossy density, the local stability assumption also guarantees that
	\begin{equation*}
		Z_{B_n}(\mathbf{0})
		= 1 + \sum_{m = 1}^\infty \frac{1}{m!} \int_{B_n^m} \kappa_m(x_1, \dots, x_m, \mathbf{0}) \, \dd\lambda^m(x_1, \dots, x_m)
		\leq \exp\bigg( \int_{B_n} \vartheta(x) \, \dd\lambda(x) \bigg)
		< \infty
	\end{equation*}
	for every $n \in \N$. Moreover, local stability implies $\kappa(x, \mu) \leq \vartheta(x) \cdot c^{\mu(\X)}$ for $\lambda$-a.e.\ $x \in \X$, all $\mu \in \NN$, and any $c \geq 1$. Hence, local stability covers all assumptions from Theorem \ref{THM abstract existence result} except the two limit relations. 
	
	Another possibility to obtain the bound on the Janossy densities in terms of an explicit assumption on $\kappa$ is to study the correlation functions of the finite Gibbs processes more closely. As the Janossy densities of these processes on their full domain are given explicitly in terms of $\kappa_m$ by Corollary \ref{COR full domain Janossy densities of Gibbs processes}, Corollary \ref{COR repres. of correlation func. via Janossy densities} implies that the correlation functions are given as
	\begin{equation*}
		\rho_{\xi_n, m}(x_1, \dots, x_m)
		= \frac{1}{Z_{B_n}(\mathbf{0})} \sum_{k = m}^\infty \frac{1}{(k - m)!} \int_{B_n^{k - m}} \kappa_k(x_1, \dots, x_k, \mathbf{0}) \, \dd\lambda^{k - m}(x_{m + 1}, \dots, x_k)
	\end{equation*}
	for $\lambda^m$-a.e.\ $(x_1, \dots, x_m) \in B_n^m$ and all $m, n \in \N$. Hence, if the right hand side, which is an explicit quantity in terms of $\kappa$, is suitably bounded, we get
	\begin{equation*}
		\sup_{n \in \N} \rho_{\xi_n, m}(x_1, \dots, x_m)
		\leq \vartheta(x_1) \cdot \dotso \cdot \vartheta(x_m)
	\end{equation*}
	for $\lambda^m$-a.e.\ $(x_1, \dots, x_m) \in \X^m$ and all $m \in \N$. This is Ruelle's bound (according to Definition \ref{DEF Ruelle's condition}) and it implies the necessary bound on the Janossy densities. In the case of pair potentials, the above quantity is studied by \cite{R:1970} and the required bound is derived under suitable assumptions. We come back to this condition in Theorem \ref{THM existence result in the pair potential setting}.
\end{remark}

\begin{remark} \label{RMK local stability bound could always be 1}
	In the definition of a Gibbs process via the GNZ equations \eqref{GNZ equation} we find the term $\kappa(x, \eta) \, \dd\lambda(x)$ and this is the only place where $\kappa$ and $\lambda$ appear. Consequently, there is some freedom in the choice of these parameters. For instance, in arguments where the specific reference measure does not matter (which includes most of our arguments), it is no loss of generality to choose the local stability bound on $\kappa$ as $\vartheta \equiv 1$. A similar argument could be used to normalize bounds on the Janossy densities or correlation functions. However, we will not resort to this normalization and carry the bound for $\kappa$ with us.
	
	Note that in allowing for $\kappa$ to be bounded by general $\vartheta$ we follow \cite{LO:2021}. \cite{J:2019} also includes such an inhomogeneity but in terms of an intensity functions associated with the measure $\lambda$.
\end{remark}

Approaching the existence proof via level sets of the specific entropy requires some kind of stationary, so the resulting existence results work only in $\R^d$ (or can probably be extended to locally compact topological groups). In $\R^d$ this method seems favorable as it leads to stationary Gibbs processes and works under weak explicit assumption on the energy function. Indeed, the paper by \cite{DV:2020}, in which the literature we mentioned in the introduction culminates (in a sense), only requires stability (a minimal assumption) and intensity regularity (an assumption that is made for similar technical reasons as the limit relations in our result).

\subsection{Initial examples and some remarks}
\label{SUBSEC examples and remarks on the general existence result}

As we have seen in Remark \ref{RMK notes on the local stability assumption}, local stability is one straight forward assumption to cover all of the prerequisites of Theorem \ref{THM abstract existence result} except the two limit relations. Probably the easiest way to ensures the validity of these limit requirements is to suppose that $\kappa$ has \textit{finite range}, meaning that for each $x \in \X$ there exists a set $B_x \in \mathcal{X}_b$ with
\begin{equation*}
	\kappa(x, \mu)
	= \kappa(x, \mu_{B_x}), \quad \mu \in \NN ,
\end{equation*}
such that, for every $B \in \mathcal{X}_b$, we have $\bigcup_{x \in B} B_x \in \mathcal{X}_b$, and such that $x \in B_y$ if, and only if, $y \in B_x$, for all $x, y \in \X$.

Denoting $\tilde{B} = \bigcup_{x \in B} B_x$, for fixed $B \in \mathcal{X}_b$, we then have $\kappa(x, \mu) = \kappa(x, \mu_{\tilde{B}})$ for all $x \in B$ and $\mu \in \NN$. Thus, if $\ell \in \N$ is large enough so that $B_\ell \supset \tilde{B}$, then
\begin{equation*}
	\kappa(x, \mu_{B_\ell})
	= \kappa\big( x, (\mu_{B_\ell})_{\tilde{B}} \big)
	= \kappa(x, \mu_{\tilde{B}})
	= \kappa(x, \mu)
\end{equation*}
for each $x \in B$ and $\mu \in \NN$, and
\begin{equation*}
	\int_B \sup_{\mu \in \NN} \big| \kappa(x, \mu_{B_\ell}) - \kappa(x, \mu) \big| \, \dd\lambda(x) 
	= 0.
\end{equation*}

If $\X$ is such that $\{ x \} \in \mathcal{X}$ for each $x \in \X$, it is always possible to include $x$ in $B_x$ while maintaining all other properties. In metric spaces the classical finite range property known from the literature trivially implies the given definition. Note that in metric spaces we always choose as a localizing structure a sequence of balls with growing radius around a fixed point in $\X$, so $\mathcal{X}_b$ corresponds to the bounded sets (with respect to the metric).
\begin{lemma} \label{LEMMA finite range property}
	Assume that $\X$ is a (Borel subset of a) metric space with metric denoted by $d$. A measurable map $\kappa : \X \times \NN \to [0, \infty)$ has finite range if there exists a fixed radius of interaction $R > 0$ such that $\kappa(x, \mu) = \kappa\big( x, \mu_{B(x, R)} \big)$ for all $x \in \X$ and $\mu \in \NN$, where $B(x, R)$ denotes the closed ball of radius $R$ around $x$ (with respect to $d$).
\end{lemma}

If we combine local stability and finite range, we can provide a neat existence result which also allows for boundary conditions without any problems concerning the formal prerequisites in the sense that $\tilde{\kappa}(\,\cdot\, , \,\cdot\,) = \kappa(\,\cdot\, , \psi + \cdot\,)$ inherits the assumptions on $\kappa$. Note that the assumptions in the following result are purely on $\kappa$ without any explicit reference to the processes involved.

\begin{corollary} \label{COR existence result with local stability and finite range}
	Let $(\X, \mathcal{X})$ be a substandard Borel space with localizing structure $B_1 \subset B_2 \subset \dotso$, let $\psi \in \NN$, and let $\lambda$ be a locally finite measure on $\X$. Let $\kappa : \X \times \NN \to [0, \infty)$ be a measurable map such that $\kappa$ satisfies the cocycle assumption \eqref{cocycle assumption} and is locally stable as in \eqref{local stability of kappa}. Furthermore, assume that $\kappa$ has finite range. Then there exists a Gibbs process $\eta$ with PI $\tilde{\kappa}(\,\cdot\, , \,\cdot\,) = \kappa(\,\cdot\, , \psi + \cdot\,)$.
\end{corollary}

\begin{example}[\bf Strauss processes] \label{EXAM Strauss process}
	Let $\X$ be a (Borel subset of a) complete separable metric space. Let $R > 0$, $c \in [0, 1]$, and $\vartheta : \X \to [0, \infty)$ a locally ($\lambda$-)integrable function. Consider the measurable map $\kappa : \X \times \NN \to [0, \infty)$ given by
	\begin{equation*}
		\kappa(x, \mu) = \vartheta(x) \cdot c^{\mu(B(x, R))} .
	\end{equation*}
	Then Corollary \ref{COR existence result with local stability and finite range} provides the existence of a Gibbs process in $\X$ with PI $\kappa$. Such a process is called \textit{Strauss process}. A special case arises for $c = 0$, where
	\begin{equation*}
		\kappa(x, \mu)
		= \vartheta(x) \cdot \mathds{1}\big\{ \mu\big( B(x, R) \big) = 0 \big\}, \quad x \in \X, ~ \mu \in \NN ,
	\end{equation*} 
	is the PI of a \textit{hard spheres} process. Of course we can also consider Strauss processes with boundary conditions.
\end{example}

Apart from Gibbs processes with a finite interaction range, the generality of Theorem \ref{THM abstract existence result} can also be utilized to obtain existence results for pair interaction processes. With a bit of technical work (involving the choice of a suitable $\pi$-system $\mathcal{X}_b^*$), roughly related to what is used in Appendix B of \cite{J:2019} but generalizing preliminary results to arbitrary measurable spaces and pair potentials with negative part, it is possible to derive the result displayed in the following. As these techniques are not new (but only refined) and do not yield significant novelties in the pair potential setting, we skip the proof. For details we refer to Chapter 4 of the authors Ph.D. thesis, \cite{B:2022}.

Suppose that $v : \X \times \X \to (- \infty, \infty]$ is measurable and symmetric, and define $\kappa$ through
\begin{equation*}
	\kappa(x, \mu)
	= \exp\bigg( - \int_\X v(x, y) \, \dd \mu(y) \bigg) \cdot \mathds{1}\bigg\{ \int_\X v^-(x, y) \, \dd\mu(y) < \infty \bigg\}, \quad x \in \X, ~ \mu \in \NN.
\end{equation*}
Notice that $\kappa$ is well-defined, measurable, and satisfies the cocycle assumption \eqref{cocycle assumption}. 
\begin{theorem} \label{THM existence result in the pair potential setting}
	Let $(\X, \mathcal{X})$ be a substandard Borel space with localizing structure $B_1 \subset B_2 \subset \dotso$, and let $\lambda$ be a locally finite measure on $\X$. Let $v : \X \times \X \to (- \infty, \infty]$ be measurable and symmetric such that $\inf_{x, y \in \X} v(x, y) \geq - A$ for a constant $A \geq 0$. Let $\kappa$ be given through $v$ as above and assume that $Z_{B_n}(\mathbf{0}) < \infty$ for each $n \in \N$. Suppose that
	\begin{equation*}
		\sup_{n \in \N} \, \frac{1}{Z_{B_n}(\mathbf{0})} \sum_{k = m}^\infty \frac{1}{(k - m)!} \int_{B_n^{k - m}} \kappa_k(x_1, \dots, x_k, \mathbf{0}) \, \dd\lambda^{k - m}(x_{m + 1}, \dots, x_k)
		\leq \vartheta(x_1) \cdot \dotso \cdot \vartheta(x_m)
	\end{equation*}
	for $\lambda^m$-a.e.\ $(x_1, \dots, x_m) \in \X^m$, all $m \in \N$, and a measurable, locally $\lambda$-integrable map $\vartheta : \X \to [0, \infty)$ which is such that, for $\lambda$-a.e.\ $x \in \X$,
	\begin{equation*}
		\int_\X \big| \mathrm{e}^{- v(x, y)} - 1 \big| \, \vartheta(y) \, \dd\lambda(y)
		< \infty
		\qquad \text{as well as} \qquad
		\int_{\X} v^-(x, y) \, \vartheta(y) \, \dd\lambda(y)
		< \infty .
	\end{equation*}
	Then there exists a Gibbs process in $\X$ with pair potential $v$ and reference measure $\lambda$.
\end{theorem}

If one restricts to non-negative pair potentials, Theorem B.1 of \cite{J:2019} is seen to be slightly generalized to arbitrary substandard Borel spaces. If one restricts to translation invariant pair potentials on $\R^d$, it is straightforward to see that superstability and lower regularity according to \cite{R:1970} also suffice for an existence result \cite[that is, we recover the existence result due to][]{R:1970}. With these last two specific settings, the existence of the classical hard core process and of various non-negative soft-core processes \citep[like the Gaussian core process due to][]{S:1976} as well as of the Lennard-Jones and Morse potential are guaranteed. These are not new results but it is certainly reassuring that our very general approach recovers them.

\section{Cluster-dependent interactions in the subcritical regime}
\label{SEC Gibbs point processes with cluster-dependent interaction}

The duality of the Boolean model and the hard sphere model, see \cite{Ht:2019}, is a point in case that the interaction between points of a Gibbs process is often intrinsically linked to some (binary) relation on the state space. In this section we recall this concept in the generality of \cite{LO:2021} and prove that the corresponding Gibbs processes exist and are unique in distribution as soon as the clusters of a Poisson process with respect to the corresponding relation are finite. We thus provide a general superstructure for results from \cite{Ht:2019}, \cite{HtH:2019}, \cite{D:2019}, and \cite{BHLV:2020}.

In this section let $(\X, d)$ be a complete separable metric space with localizing structure $B_n = B(x_0, n)$, $n \in \N$, for some fixed $x_0 \in \X$. Denote by $\mathcal{X}$ the Borel $\sigma$-field of $\X$ and by $\mathcal{X}_b$ the bounded Borel sets, and let $\lambda$ be a locally finite measure on $\X$. The restriction to metric spaces has technical reasons, namely that we need to deal with weak convergence of measures in Lemma \ref{LEMMA domination of Gibbs processes via Poisson processes}. 

Let $\sim$ be a symmetric relation on $\X$ such that $\{(x, y) \in \X^2 : x \sim y\} \in \mathcal{X}^{\otimes 2}$. We call $x, y \in \X$ connected via $D \subset \X$ if there exist $n \in \N_0$ and $z_1, \dots, z_n \in D$ such that $z_j \sim z_{j + 1}$ for each $j \in \{ 0, \dots, n \}$, where we put $z_0 = x$ and $z_{n + 1} = y$. This last terminology is also used for counting measures $\mu \in \NN$, think of $D = \mathrm{supp}(\mu) = \{ x \in \X : \mu(\{ x \}) > 0 \}$. For formal completeness observe that $\mathrm{supp}(\mu) \in \mathcal{X}$ by Lemma \ref{LEMMA Appen. prelim. statement, implications of meas. diag.}. If $x, y$ are connected via $\mu$, we write $x \stackrel{\mu}{\sim} y$. Also if $x$ is connected via $\mu$ to some point in $\psi \in \NN$, we write $x \overset{\mu}{\sim} \psi$. We define
\begin{equation} \label{mu-cluster of x}
C(x, \mu)
= \int_\X \mathds{1}\{ y \in \cdot \, \} \, \mathds{1}\big\{ x \stackrel{\mu}{\sim} y \big\} \, \dd\mu(y) 
\end{equation}
and call $C(x, \mu)$ the \textit{($\mu$-)cluster} of $x$. In the given setting the points of any counting measure can be chosen measurably with the help of Lemma 1.6 of \cite{K:2017} which allows for a measurable construction of these clusters. In other words, the mapping
\begin{equation*}
	(x, \mu) \mapsto C(x, \mu) \in \NN
\end{equation*}
is measurable and so is the map $(x, \mu, \psi) \mapsto \mathds{1}\big\{ x \stackrel{\mu}{\sim} \psi \big\}$. In the context of disagreement couplings or when considering Gibbsian particle processes it is instrumental to consider interactions which only depend on corresponding clusters, that is, we have
\begin{equation*}
\kappa(x, \mu)
= \kappa\big( x, C(x, \mu) \big) .
\end{equation*}
Particularly interesting is the case where the clusters, and hence the interaction range, are infinite. Theorem \ref{THM abstract existence result} allows us to prove existence in that setting whenever we are in a \textit{subcritical} regime with respect to a suitable Poisson process. In the same regime we also prove uniqueness.

\begin{theorem} \label{THM existence and uniqueness for subcritical cluster-dependent interactions}
	Let $(\X, d)$ be a complete separable metric space as above, and let $\lambda$ be a locally finite measure on $\X$. Let $\kappa : \X \times \NN \to [0, \infty)$ be a measurable map which satisfies the cocycle relation \eqref{cocycle assumption} and is locally stable as in \eqref{local stability of kappa} with bound $\vartheta$. Moreover, suppose that
	\begin{equation} \label{cluster dependence of kappa}
	\kappa(x, \mu)
	= \kappa\big( x, C(x, \mu) \big), \quad x \in \X, \, \mu \in \NN ,
	\end{equation}
	and that, for $\lambda$-a.e.\ $x \in \X$,
	\begin{equation*}
	\Pi_{\vartheta \lambda}\big( \big\{ \mu \in \NN \, : \, C(x, \mu)(\X) < \infty \big\} \big)
	= 1 .
	\end{equation*}
	Then, up to equality in distribution, there exists exactly one Gibbs process in $\X$ with PI $\kappa$.
\end{theorem}

\subsection{Notation and preliminaries for the proof of Theorem \ref{THM existence and uniqueness for subcritical cluster-dependent interactions}}
\label{SUBSEC Notation and preliminaries for cluster-dependent interaction processes}

For the proof of both the existence and uniqueness part of Theorem \ref{THM existence and uniqueness for subcritical cluster-dependent interactions}, it is essential that any Gibbs process is dominated by a suitable Poisson process. For finite Gibbs processes this is known from \cite{GK:1997}, but the extension to infinite processes takes additional effort. Conceptually we follow \cite{LO:2021} in doing so. However, we need further properties related to the local convergence in the context of the existence result. As mentioned at the beginning of this section, the assumption of $\X$ being a complete separable metric space becomes necessary for the use of weak convergence in order to derive the following lemma. A similar result could also be obtained by combining Lemma 5.3 of \cite{LO:2021} and the appendix of \cite{SS:2014} (in fact, we conjecture that such an approach would allow to replace the $\Phi_n$'s by a single Poisson process, but as this is no harm to us, the additional technical details seem unnecessary).
\begin{lemma} \label{LEMMA domination of Gibbs processes via Poisson processes}
	In the setting given at the beginning of Section \ref{SEC Gibbs point processes with cluster-dependent interaction}, let $\kappa : \X \times \NN \to [0, \infty)$ be a measurable function that satisfies the cocycle relation \eqref{cocycle assumption} and is locally stable as in \eqref{local stability of kappa} with bound $\vartheta$. For each $n \in \N$, let $\xi_n$ be a finite Gibbs process with PI $\kappa^{(B_n, \mathbf{0})}$, and let $\eta$ be (one of) the local limit process(es) obtainable from Theorem \ref{THM existence of locally convergent subsequences}. Then (possibly after extending the underlying probability space) there exist Poisson processes $\Phi_n$ with intensity measures $\vartheta \lambda_{B_n}$ ($n \in \N$) and a Poisson process $\Phi$ with intensity measure $\vartheta \lambda$ as well as point processes $\tilde{\xi}_n$ ($n \in \N$) and $\tilde{\eta}$ such that
	\begin{equation*}
		\tilde{\xi}_n \stackrel{d}{=} \xi_n \quad \text{and} \quad \tilde{\xi}_n \leq \Phi_n \quad (\PP\text{-a.s.})
	\end{equation*}
	for all $n \in \N$,
	\begin{equation*}
		\tilde{\eta} \stackrel{d}{=} \eta \quad \text{and} \quad \tilde{\eta} \leq \Phi \quad (\PP\text{-a.s.}),
	\end{equation*}
	and $\tilde{\xi}_{n_k} \stackrel{loc}{\longrightarrow} \tilde{\eta}$ (as $k \to \infty$) along a suitable subsequence. Here $\stackrel{d}{=}$ denotes equality in distribution and two measures $\mu, \nu$ on $\X$ satisfy $\mu \leq \nu$ if $\mu(B) \leq \nu(B)$ for all $B \in \mathcal{X}$.
\end{lemma}
\begin{prf}
	Fix $n \in \N$. By Example 2.1 of \cite{GK:1997} (available only for the finite processes) and the Strassen theorem \citep[cf.][]{L:1992}, there exists a point process $\tilde{\xi}_n$ with $\tilde{\xi}_n \stackrel{d}{=} \xi_n$ and a Poisson process $\Phi_n$ with intensity measure $\vartheta \lambda_{B_n}$ such that $\tilde{\xi}_n \leq \Phi_n$ almost surely. Let $(\xi_{n_k})_{k \in \N}$ be the subsequence for which $\xi_{n_k} \stackrel{loc}{\longrightarrow} \eta$ as $k \to \infty$. By Proposition \ref{PROP basic properties of local convergence} we have $\tilde{\xi}_{n_k} \stackrel{d}{\longrightarrow} \eta$, where $\stackrel{d}{\longrightarrow}$ denotes convergence in distribution, implicitly using the metric structure that $\NN$ inherits from $\X$ as detailed in Appendix 2 of \cite{DV:2005}. Standard reasoning \citep[e.g. via Theorem 16.3 of][]{K:2002} yields that $(\tilde{\xi}_{n_k}, \Phi_{n_k})_{k \in \N}$ is a tight sequence. Thus, there exists a subsequence which converges in distribution to a limit element $(\tilde{\eta}, \Phi)$, where $\Phi$ is a Poisson process with intensity measure $\vartheta \lambda$ and $\tilde{\eta} \stackrel{d}{=} \eta$. Using Proposition A.2.6 of \cite{DV:2005} to reduce to $\NN_f$ and arguing by contradiction it is easy to see that the set $\{(\mu, \nu) \in \NN \times \NN : \mu \leq \nu\}$ is closed. The weak convergence together with $\tilde{\xi}_n \leq \Phi_n$ (a.s.) and the Portmanteau theorem \citep[Theorem 4.25 of][]{K:2002} yield
	\begin{equation*}
	\PP(\tilde{\eta} \leq \Phi) = 1 .
	\end{equation*}
	Moreover, for each measurable, local, and tame map $F : \NN \to [0, \infty)$ we have
	\begin{equation*}
	\E\big[ F\big( \tilde{\xi}_{n_k} \big) \big]
	= \E\big[ F(\xi_{n_k}) \big]
	\longrightarrow \E\big[ F(\eta) \big]
	= \E\big[ F(\tilde{\eta}) \big]
	\end{equation*}
	as $k \to \infty$ which implies the asserted local convergence along the constructed subsequence.
\end{prf}

\vspace{3mm}

The tools we use to prove the uniqueness part of Theorem \ref{THM existence and uniqueness for subcritical cluster-dependent interactions} are only available for diffuse reference measures. A simple \textit{randomization} property, a special case of which was already used by \cite{SS:2014} (in a general manner) and \cite{BL:2021} (in the context of uniqueness results), allows for a workaround. We generalize this result slightly and give a proof.

Recall that any point process in a Borel space, in the sense of \cite{K:2017} and \cite{LP:2017}, is \textit{proper} meaning that it can be written as a sum of random Dirac measures.
\begin{lemma} \label{LEMMA randomization of Gibbs processes}
	Let $(\X, \mathcal{X})$ be a localized Borel space, $\lambda$ a locally finite reference measure on $\X$, and let $(\mathbb{Y}, \mathcal{Y})$ be a measurable space endowed with a probability measure $Q$. Let $\kappa : \X \times \NN \to [0, \infty)$ be a measurable function that satisfies the cocycle relation \eqref{cocycle assumption}. Put $\tilde{\kappa} : \X \times \mathbb{Y} \times \NN(\X \times \mathbb{Y}) \to [0, \infty)$,
	\begin{equation*}
		\tilde{\kappa}(x, r, \mu)
		= \kappa\big( x, \mu(\,\cdot \times \mathbb{Y}) \big) ,
	\end{equation*}
	and let $\eta = \sum_{j = 1}^{\eta(\X)} \delta_{X_j}$ be a Gibbs process in $\X$ with PI $\kappa$ and reference measure $\lambda$. Let $R_1, R_2, \dotso$ be independent random variables in $\mathbb{Y}$ distributed according to $Q$, with the whole sequence independent of $\eta$. Then the randomization $\tilde{\eta} = \sum_{j = 1}^{\eta(\X)} \delta_{(X_j, R_j)}$ of $\eta$ is a Gibbs process in $\X \times \mathbb{Y}$ with PI $\tilde{\kappa}$ and reference measure $\lambda \otimes Q$.
\end{lemma}

Note that the product space $\X \times \mathbb{Y}$ is localized by the sets $B_1 \times \mathbb{Y} \subset B_2 \times \mathbb{Y}_2 \subset \ldots$.

\vspace{2mm}

\begin{prf}
	For $\mu = \sum_{j = 1}^{\mu(\X)} \delta_{x_j} \in \NN$, with $x_1, x_2, \ldots \in \X$, and $s = (s_j)_{j \in \N} \in \mathbb{Y}^\infty = \bigotimes_{j = 1}^\infty \mathbb{Y}$ define
	\begin{equation*}
		T(\mu, s)
		= \sum_{j = 1}^{\mu(\X)} \delta_{(x_j, s_j)}
	\end{equation*}
	which constitutes a measurable map $\NN \times \mathbb{Y}^\infty \to \NN(\X \times \mathbb{Y})$ by Lemma 1.6 of \cite{K:2017}. Write $R = (R_j)_{j \in \N}$ and $\mathbb{Q} = \bigotimes_{j = 1}^\infty Q$. For any measurable function $f : \X \times \mathbb{Y} \times \NN(\X \times \mathbb{Y}) \to [0, \infty]$ we have
	\begin{equation*}
		\E\bigg[ \int_{\X \times \mathbb{Y}} f(x, r, \tilde{\eta}) \, \dd \tilde{\eta}(x, r) \bigg]
		= \E\bigg[ \int_{\X \times \mathbb{Y}} f\big( x, r, T(\eta, R) \big) \, \dd \big( T(\eta, R) \big)(x, r) \bigg] .
	\end{equation*}
	Due to the given independence properties, this term equals
	\begin{equation} \label{eq in the proof of the randomization property}
		\int_{\mathbb{Y}^\infty} \E\bigg[ \int_{\X \times \mathbb{Y}} f\big( x, r, T(\eta, s) \big) \, \dd \big( T(\eta, s) \big)(x, r) \bigg] \dd\mathbb{Q}(s)
		= \E\Bigg[ \sum_{j = 1}^{\eta(\X)} \int_{\mathbb{Y}^\infty} f\big( X_j, s_j, T(\eta, s) \big) \, \dd\mathbb{Q}(s) \Bigg] .
	\end{equation}
	In order to write the term in the sum in \eqref{eq in the proof of the randomization property} as $g(X_j, \eta)$ for a suitable map $g$, we have to ensure that $g$ assigns the marks correctly in the sense that the fact that $s_j$ is the mark of $X_j$ is also relevant in $T(\eta, s)$. Upon defining the measurable map $g : \X \times \NN \to [0, \infty]$ by
	\begin{equation*}
		g(x, \mu)
		= \mathds{1}\big\{ \mu(\{ x \}) > 0 \big\} \int_{\mathbb{Y}^\infty} \int_\mathbb{Y} f\big( x, r, T(\mu \setminus \delta_x, s) + \delta_{(x, r)} \big) \, \dd Q(r) \, \dd \mathbb{Q}(s) ,
	\end{equation*}
	this is ensured as the infinite product structure implies (almost surely)
	\begin{equation*}
		g(X_j, \eta)
		= \int_{\mathbb{Y}^\infty} \int_\mathbb{Y} f\big( X_j, r, T(\eta \setminus \delta_{X_j}, s) + \delta_{(X_j, r)} \big) \, \dd Q(r) \, \dd \mathbb{Q}(s)
		= \int_{\mathbb{Y}^\infty} f\big( X_j, s_j, T(\eta, s) \big) \, \dd\mathbb{Q}(s) .
	\end{equation*}
	Hence, picking up the previous calculation from \eqref{eq in the proof of the randomization property}, we have
	\begin{equation*}
		\E\bigg[ \int_{\X \times \mathbb{Y}} f(x, r, \tilde{\eta}) \, \dd \tilde{\eta}(x, r) \bigg]
		= \E\Bigg[ \sum_{j = 1}^{\eta(\X)} g(X_j, \eta) \Bigg]
		= \E\bigg[ \int_\X g(x, \eta) \, \dd \eta(x) \bigg] .
	\end{equation*}
	By the GNZ equation for $\eta$ and the given independence properties, this term further equals
	\begin{align*}
		\E\bigg[ \int_\X g(x, \eta + \delta_x) \, \kappa(x, \eta) \, \dd \lambda(x) \bigg]
		&= \E\bigg[ \int_\X \bigg( \int_{\mathbb{Y}^\infty} \int_\mathbb{Y} f\big( x, r, T(\eta, s) + \delta_{(x, r)} \big) \, \dd Q(r) \, \dd \mathbb{Q}(s) \bigg) \kappa(x, \eta) \, \dd\lambda(x) \bigg] \\
		&= \E\bigg[  \int_\X \int_\mathbb{Y} f\big( x, r, T(\eta, R) + \delta_{(x, r)} \big) \, \kappa(x, \eta) \, \dd Q(r) \, \dd\lambda(x) \bigg] \\
		&= \E\bigg[ \int_{\X \times \mathbb{Y}} f(x, r, \tilde{\eta} + \delta_{(x, r)}) \, \tilde{\kappa}(x, r, \tilde{\eta}) \, \dd (\lambda \otimes Q)(x, r) \bigg] ,
	\end{align*}
	which concludes the proof.
\end{prf}

\vspace{3mm}

We recall the following \textit{disagreement coupling} from \cite{LO:2021} (their Theorem 6.3), adapted to our notation and in a strictly less general version focusing on locally stable Papangelou intensities. For two finite counting measures $\mu, \mu' \in \NN_f$, we denote by $|\mu - \mu'| \in \NN_f$ the total variation measure of the signed measure $\mu - \mu'$ defined via the Jordan-decomposition, see Corollary 3.1.2 of \cite{B:2007}. It is straight forward to show that the map $\NN_f \times \NN_f \ni (\mu, \mu') \mapsto |\mu - \mu'| \in \NN_f$ is measurable.

\begin{prop}[\bf Disagreement coupling due to \citeauthor{LO:2021}] \label{PROP disagreement coupling}
	In the setting given at the beginning of Section \ref{SEC Gibbs point processes with cluster-dependent interaction} assume that $\lambda$ is diffuse and let $\kappa : \X \times \NN \to [0, \infty)$ be a measurable function that satisfies the cocycle relation \eqref{cocycle assumption} and is locally stable as in \eqref{local stability of kappa} with bound $\vartheta$. Moreover, assume that
	\begin{equation*}
	\kappa(x, \mu)
	= \kappa\big( x, C(x, \mu) \big), \quad x \in \X, \, \mu \in \NN .
	\end{equation*}
	Let $W \in \mathcal{X}_b$ and $\psi, \psi' \in \NN_{W^\mathsf{c}}$. There exists a Gibbs process $\xi$ with PI $\kappa^{(W, \psi)}$ and a Gibbs process $\xi'$ with PI $\kappa^{(W, \psi')}$ such that $\xi \leq \Psi$ and $\xi' \leq \Psi$ (almost surely), where $\Psi$ is a Poisson process in $\X$ with intensity measure $\vartheta \lambda$, and such that every point in $|\xi - \xi'|$ is connected via $\xi + \xi'$ to some point in $\psi + \psi'$.
\end{prop}

Note that if $\psi = \psi'$, the two so constructed Gibbs processes are identical (path-wise). The concept of disagreement percolation was first introduced by \cite{vdBM:1994} in the discrete setting and transferred to the continuum by \cite{Ht:2019} and \cite{HtH:2019} though the latter works contain technical imprecisions. The above result by \cite{LO:2021} provides a rigorous, up-to-date, and very general version of the underlying technique.

\subsection{A proof of Theorem \ref{THM existence and uniqueness for subcritical cluster-dependent interactions}}
\label{SUBSEC Proof of existence and uniqueness of cluster-dependent interaction processes}

We first prove existence. By Remark \ref{RMK notes on the local stability assumption}, the local stability assumption on $\kappa$ covers all prerequisites of Theorem \ref{THM abstract existence result} except the limit relations. To show that these are also satisfied in the given setting, let $\xi_n$ be a Gibbs process with PI $\kappa^{(B_n, \mathbf{0})}$ ($n \in \N$) and denote by $\eta$ a corresponding limit process from Theorem \ref{THM existence of locally convergent subsequences}. We assume, without loss of generality (by Lemma \ref{LEMMA domination of Gibbs processes via Poisson processes}), that $\xi_n \leq \Phi_n$ for each $n\in \N$ and $\eta \leq \Phi$ almost surely, where $\Phi_n$ is a Poisson process in $\X$ with intensity measure $\vartheta \lambda_{B_n}$ ($n \in \N$) and $\Phi$ is a Poisson process with intensity measure $\vartheta \lambda$. We denote the corresponding combined $\PP$-null set by $\Omega_0$. If, for $\ell \in \N$, $x \in \X$, and $\omega \in \Omega_0^\mathsf{c}$, we have $C\big( x, \Phi_n(\omega) \big)(B_\ell^\mathsf{c}) = 0$, then $C\big( x, \xi_n(\omega) \big)(B_\ell^\mathsf{c}) = 0$ and therefore
\begin{equation} \label{from the existence proof in the cluster-dependence setting}
	C\big( x, (\xi_n)_{B_\ell}(\omega) \big)
	= C\big( x, \xi_n(\omega) \big)_{B_\ell}
	= C\big( x, \xi_n(\omega) \big) .
\end{equation}
For $n, \ell \in \N$ and $\lambda$-a.e.\ $x \in \X$ we have, by \eqref{from the existence proof in the cluster-dependence setting} and the assumptions on $\kappa$,
\begin{align*}
	\E \big| \kappa\big( x, (\xi_n)_{B_\ell} \big) - \kappa(x, \xi_n) \big|
	&= \E \Big| \Big( \kappa\big( x, (\xi_n)_{B_\ell} \big) - \kappa(x, \xi_n) \Big) \mathds{1}\big\{ C(x, \Phi_n)(B_\ell^\mathsf{c}) > 0 \big\} \Big| \\
	&\leq \vartheta(x) \cdot \PP\big( C(x, \Phi_{B_n})(B_\ell^\mathsf{c}) > 0 \big) \\
	&\leq \vartheta(x) \cdot \PP\big( C(x, \Phi)(B_\ell^\mathsf{c}) > 0 \big) .
\end{align*}
Combining this with the observation that
\begin{equation} \label{convergence of connection probability}
	\limsup_{\ell \to \infty} \PP\big( C(x, \Phi)(B_\ell^\mathsf{c}) > 0 \big)
	= 1 - \PP\bigg( \bigcup_{\ell = 1}^\infty \Big\{ \omega \in \Omega : C\big(x, \Phi(\omega) \big)(B_\ell^\mathsf{c}) = 0 \Big\} \bigg)
	= 1 - \PP\big( C(x, \Phi)(\X) < \infty \big)
	= 0 
\end{equation}
for $\lambda$-a.e.\ $x \in \X$, dominated convergence (using the local integrability of $\vartheta$) gives
\begin{equation*}
	\limsup_{\ell \to \infty} \, \sup_{n \in \N} \int_B \E\big| \kappa\big( x, (\xi_{n})_{B_\ell} \big) - \kappa(x, \xi_{n}) \big| \, \dd\lambda(x)
	\leq \limsup_{\ell \to \infty} \int_B \vartheta(x) \cdot \PP\big( C(x, \Phi)(B_\ell^\mathsf{c}) > 0 \big) \, \dd\lambda(x)
	= 0 
\end{equation*}
for every $B \in \mathcal{X}_b$. Similarly, $\lim_{\ell \to \infty} \int_B \E| \kappa(x, \eta_{B_\ell}) - \kappa(x, \eta) | \, \dd\lambda(x) = 0$, and Theorem \ref{THM abstract existence result} yields the existence result. 

Next, we argue that for the uniqueness part we can assume, without loss of generality, that $\lambda$ is diffuse. The reasoning is similar to that of \cite{BL:2021}. Consider $\tilde{\X} = \X \times [0, 1]$ equipped with some complete metric that induces the product topology and let $\tilde{\lambda} = \lambda \otimes \mathcal{L}^1_{[0, 1]}$, where $\mathcal{L}^d$ denotes the $d$-dimensional Borel-Lebesgue measure. For $x \in \X$, $r \in [0, 1]$, and $\mu \in \NN(\tilde{\X})$ let $\tilde{C}(x, r, \mu)$ be the $\mu$-cluster of $(x, r)$, where $(x, r), (y, s) \in \tilde{\X}$ are connected simply when $x \sim y$. We have $\tilde{C}(x, r, \mu)(\tilde{\X}) = C\big( x, \mu(\, \cdot \times [0, 1]) \big)(\X)$, so if $\Phi$ is a Poisson process with intensity measure $\vartheta \lambda$ and $\tilde{\Phi}$ is a uniform randomization of $\Phi$ (as in Lemma \ref{LEMMA randomization of Gibbs processes}), which is a Poisson process in $\tilde{\X}$ with intensity measure $\vartheta \lambda \otimes \mathcal{L}^1_{[0, 1]}$ by the marking theorem \cite[Theorem 5.6 of][]{LP:2017}, then the condition $\PP\big( C(x, \Phi)(\X) < \infty \big) = 1$ immediately implies
\begin{equation*}
	\PP\big( \tilde{C}(x, r, \tilde\Phi)(\tilde\X) < \infty \big) = 1
\end{equation*}
for $\lambda$-a.e.\ $x \in \X$ and any $r \in [0,1]$. For $x \in \X$, $r \in [0, 1]$, and $\mu \in \NN(\tilde{\X})$ we also define
\begin{equation*}
	\tilde{\kappa}(x, r, \mu)
	= \kappa\big( x, \mu(\, \cdot \times [0, 1]) \big) ,
\end{equation*}
which inherits the cocycle and local stability property of $\kappa$. Moreover, $\tilde{\kappa}\big( x, r, \tilde{C}(x, r, \mu) \big) = \tilde{\kappa}(x, r, \mu)$. If $\eta, \eta'$ are Gibbs processes with PI $\kappa$ and reference measure $\lambda$, then by Lemma \ref{LEMMA randomization of Gibbs processes} the uniform randomizations $\tilde\eta, \tilde\eta'$ of $\eta$ and $\eta'$ are Gibbs processes with PI $\tilde\kappa$ and reference measure $\tilde\lambda$. Thus, if the uniqueness result holds for diffuse reference measures, then $\tilde\eta \overset{d}{=} \tilde\eta'$ as $\tilde\lambda$ is diffuse, and we obtain
\begin{equation*}
	\eta
	\, = \, \tilde\eta(\, \cdot \times [0, 1])
	\, \overset{d}{=} \, \tilde\eta'(\, \cdot \times [0, 1])
	\, = \, \eta' .
\end{equation*}
We conclude that the general result holds if it is proven for diffuse reference measures. Hence, let us assume that $\lambda$ is diffuse.

For $B, W \in \mathcal{X}_b$ with $B \subset W$ and $\psi, \psi' \in \NN_{W^\mathsf{c}}$, we have
\begin{equation*}
	\big| \mathrm{P}_{W, \psi}(E) - \mathrm{P}_{W, \psi'}(E) \big|
	= \big| \PP(\xi \in E) - \PP(\xi' \in E) \big|
	\leq \E\big| \mathds{1}\{ \xi_B \in E \} - \mathds{1}\{ \xi'_B \in E \} \big|
	\leq \PP\big( \xi_B \neq \xi'_B \big)
	= \PP\big( |\xi - \xi'|(B) > 0 \big) .
\end{equation*}
for all $E \in \mathcal{N}_B$, where $\xi, \xi'$ are the Gibbs processes from the disagreement coupling in Proposition \ref{PROP disagreement coupling}. Since each point in $|\xi - \xi'| \, (\leq \Psi_W)$ is connected via $\xi + \xi'$ (hence via $\Psi_W$) to some point in $\psi + \psi'$, the probability on the right hand side is bounded by
\begin{equation*}
	\E\bigg[ \int_B \mathds{1}\Big\{ x \overset{\Psi_W}{\sim} (\psi + \psi') \Big\} \, \dd\Psi(x) \bigg]
\end{equation*}
which, due to Mecke's equation, equals
\begin{equation} \label{bound on the disagreement}
	\int_B \PP\Big( x \overset{\Psi_W + \delta_x}{\sim} (\psi + \psi') \Big) \, \vartheta(x) \, \dd\lambda(x)
	= \int_B \PP\Big( x \overset{\Psi_W}{\sim} (\psi + \psi') \Big) \, \vartheta(x) \, \dd\lambda(x) .
\end{equation}
Now, let $\eta, \eta'$ be two Gibbs processes in $\X$ with PI $\kappa$. According to Lemma 5.3 of \cite{LO:2021} we can assume, without loss of generality, that there exist two Poisson processes $\Phi, \Phi'$ with intensity measure $\vartheta \lambda$ such that $\eta \leq \Phi$ and $\eta' \leq \Phi'$ (almost surely) and such that $(\eta, \eta', \Phi, \Phi')$ is independent of $\Psi$. Let $B \in \mathcal{X}_b$ be arbitrary and choose $\ell$ large enough so that $B \subset B_\ell$. Take $E \in \mathcal{N}_B$. By Corollary \ref{COR DLR reformulation in terms of finite Gibbs proc. distr.} and the bound \eqref{bound on the disagreement} applied to $W = B_\ell$, we obtain
\begin{equation*}
	\big| \PP(\eta \in E) - \PP(\eta' \in E) \big|
	\leq \E \Big| \mathrm{P}_{B_\ell, \eta_{B_\ell^\mathsf{c}}}(E) - \mathrm{P}_{B_\ell, \eta'_{B_\ell^\mathsf{c}}}(E) \Big|
	\leq \int_B \PP\Big( x \overset{\Psi_{B_\ell}}{\sim} \big( \eta_{B_\ell^\mathsf{c}} + \eta'_{B_\ell^\mathsf{c}} \big) \Big) \, \vartheta(x) \, \dd\lambda(x) .
\end{equation*}
The right hand side is clearly bounded by
\begin{equation*}
	\int_B \PP\Big( x \overset{\Psi_{B_\ell}}{\sim} \eta_{B_\ell^\mathsf{c}} \Big) \, \vartheta(x) \, \dd\lambda(x)
	+ \int_B \PP\Big( x \overset{\Psi_{B_\ell}}{\sim} \eta'_{B_\ell^\mathsf{c}} \Big) \, \vartheta(x) \, \dd\lambda(x)
\end{equation*}
and as $\eta_{B_\ell^\mathsf{c}} \leq \Phi_{B_\ell^\mathsf{c}}$ and $\eta'_{B_\ell^\mathsf{c}} \leq \Phi'_{B_\ell^\mathsf{c}}$ ($\PP$-a.s.), the independence properties of the Poisson process yield
\begin{equation*}
	\big| \PP(\eta \in E) - \PP(\eta' \in E) \big|
	\leq 2 \int_B \PP\Big( x \overset{\Phi_{B_\ell}}{\sim} \Phi_{B_\ell^\mathsf{c}} \Big) \, \vartheta(x) \, \dd\lambda(x)
	\leq 2 \int_B \PP\big( C(x, \Phi)(B_\ell^\mathsf{c}) > 0 \big) \, \vartheta(x) \, \dd\lambda(x) .
\end{equation*}
It follows from \eqref{convergence of connection probability} and dominated convergence that $\PP^{\eta} = \PP^{\eta'}$ on the algebra $\mathcal{Z}$ that generates $\mathcal{N}$, and therefore $\PP^{\eta} = \PP^{\eta'}$. This finishes the proof. \qed

\begin{remark}
	Notice that in equation \eqref{from the existence proof in the cluster-dependence setting}, in general, we only have
	\begin{equation*}
		C\big( x, (\xi_n)_{B_\ell}(\omega) \big)
		\leq C\big( x, \xi_n(\omega) \big)_{B_\ell}
	\end{equation*}
	as there could be points in $(\xi_n)_{B_\ell}(\omega)$ which are connected to $x$ via points in $(\xi_n)_{B_\ell^\mathsf{c}}$. However, if $C\big( x, \xi_n(\omega) \big)(B_\ell^\mathsf{c}) = 0$ this is not possible.
\end{remark}

\section{Gibbs particle processes}
\label{SEC existence and uniqueness of Gibbs particle processes}

Let $(\mathbb{B}, d)$ be a complete separable metric space. Denote by $\mathcal{B}(\mathbb{B})$ the Borel subsets of $\mathbb{B}$ and by $\mathcal{B}_b(\mathbb{B})$ those Borel sets which are bounded with respect to the metric $d$. Let $\mathcal{C}(\mathbb{B})$ be the space of compact subsets (particles) of $\mathbb{B}$ and equip $(\X =) \, \mathcal{C}^*(\mathbb{B}) = \mathcal{C}(\mathbb{B}) \setminus \{ \varnothing \}$ with the Hausdorff metric $d_H$. For a definition of $d_H$ and the fact that $\mathcal{C}^*(\mathbb{B})$ (equipped with $d_H$) is a complete separable metric space, we refer to Appendix D of \cite{M:2017}. Denote the corresponding Borel $\sigma$-field by $(\mathcal{X} =) \, \mathcal{B}( \mathcal{C}^*(\mathbb{B}) )$ and write $(\mathcal{X}_b =) \, \mathcal{B}_b( \mathcal{C}^*(\mathbb{B}) )$ for those Borel sets which are bounded with respect to $d_H$. Let $\lambda$ be a locally finite measure on $\mathcal{C}^*(\mathbb{B})$. In line with our setting in general measurable spaces, we consider the locally finite counting measures
\begin{equation*}
	\NN( \mathcal{C}^*(\mathbb{B}) )
	= \big\{ \mu \text{ measure on } \mathcal{C}^*(\mathbb{B}) : \mu\big( B_H(K, r) \big) \in \N_0 \text{ for all } K \in \mathcal{C}^*(\mathbb{B}), \, r \geq 0 \big\} ,
\end{equation*}
where $B_H(K, r)$ denotes a ball of radius $r$ in $\mathcal{C}^*(\mathbb{B})$ around $K$ with respect to $d_H$. Clearly our general setup from previous sections covers this setting and yields a definition and properties of Gibbs processes in $\mathcal{C}^*(\mathbb{B})$, which are called \textit{Gibbsian particle processes}. Corollary \ref{COR existence result with local stability and finite range} immediately applies in this setting and guarantees the existence of a large class of Gibbs particle processes. As an explicit example we use the result to prove existence of the processes considered by \cite{BHLV:2020}, a question left unanswered in their work. Any issues of measurability can be cleared with straight forward standard arguments and we omit the details.

\begin{example}[\bf Admissible Gibbs particle processes due to \citeauthor{BHLV:2020}] \label{EXA Gibbs particle processes}
	In their work, \cite{BHLV:2020} conjectured that the Gibbs particle processes they consider exist. In the following, we show the existence in a slightly more general setting which fully includes theirs, thus proving their conjecture. We restrict our attention to $\mathbb{B} = \R^d$, the $d$-dimensional space with the Euclidean metric. We consider the space $\mathcal{C}^{(d)} = \mathcal{C}(\R^d) \setminus \{ \varnothing \}$ of non-empty compact subsets of $\R^d$ equipped with the Borel $\sigma$-field induced by the Hausdorff metric. We let $c : \mathcal{C}^{(d)} \to \R^d$ be a center function, that is, $c$ is measurable and satisfies $c(K + x) = c(K) + x$ for all $K \in \mathcal{C}^{(d)}$ and $x \in \R^d$. Let $\mathbb{Q}$ be a probability measure on $\mathcal{C}^{(d)}$ with $\mathbb{Q}(\mathcal{C}^{(d)}_0) = 1$, where
	\begin{equation*}
		\mathcal{C}^{(d)}_0 = \{ K \in \mathcal{C}^{(d)} : c(K) = 0 \} .
	\end{equation*}
	Moreover, assume that 
	\begin{equation*}
		\int_{\mathcal{C}^{(d)}} \mathcal{L}^d\big( L + (-C) \big) \, \dd\mathbb{Q}(L) < \infty, \quad C \in \mathcal{C}^{(d)} ,
	\end{equation*}
	where adding two subsets of $\R^d$ means taking their Minkowski sum. This last condition is certainly satisfied if the particle size is bounded, that is, if $\mathbb{Q}\big( \big\{ K \in \mathcal{C}^{(d)} : K \subset B(0, R) \big\} \big) = 1$ for a fixed $R > 0$. For some fixed intensity parameter $z > 0$ consider the measure
	\begin{equation*}
		\lambda(\cdot)
		= z \int_{\R^d} \int_{\mathcal{C}^{(d)}} \mathds{1}\{ L + x \in \cdot \, \} \, \dd\mathbb{Q}(L) \, \dd x
	\end{equation*}
	on $\mathcal{C}^{(d)}$. The measure $\lambda$ is locally finite since, for any $K \in \mathcal{C}^{(d)}$ and $r \geq 0$, there exists a set $C \in \mathcal{C}^{(d)}$ such that $L \subset C$ for all $L \in B_H(K, r)$, and
	\begin{equation*}
		\lambda\big( B_H(K, r) \big)
		\leq z \int_{\R^d} \int_{\mathcal{C}^{(d)}} \mathds{1}\big\{ (L + x) \cap C \neq \varnothing \big\} \, \dd\mathbb{Q}(L) \, \dd x 
		= z \int_{\mathcal{C}^{(d)}} \mathcal{L}^d\big( L + (-C) \big) \, \dd\mathbb{Q}(L)
		< \infty ,
	\end{equation*}
	by assumption on $\mathbb{Q}$. Let $V = \{ V_n : n \in \N, \, n \geq 2 \}$ be a collection of higher-order interaction potentials, meaning that, for each $n \geq 2$, $V_n : (\mathcal{C}^{(d)})^n \to (- \infty, \infty]$ is a measurable and symmetric function. Assume that there exists some $R_V > 0$ such that $V_n(K_1, \dots, K_n) = 0$ for all $K_1, \dots, K_n \in \mathcal{C}^{(d)}$ with $\max\{ d_H(K_i, K_j) : 1 \leq i < j \leq n\} > R_V$ and all $n \geq 2$. Moreover, assume that
	\begin{equation*}
		\sum_{n = 2}^\infty \frac{1}{(n - 1)!} \, \max\bigg\{ - \int_{(\mathcal{C}^{(d)})^{n - 1}} V_n(K, L_1, \dots, L_{n - 1}) \, \dd\mu^{(n - 1)}(L_1, \dots, L_{n - 1}) , \, 0 \bigg\}
		< \infty, \quad K \in \mathcal{C}^{(d)}, \, \mu \in \NN(\mathcal{C}^{(d)}) .
	\end{equation*}
	Define $\kappa : \mathcal{C}^{(d)} \times \NN(\mathcal{C}^{(d)}) \to [0, \infty)$ as
	\begin{equation*}
		\kappa(K, \mu)
		= \exp\bigg( - \sum_{n = 2}^\infty \frac{1}{(n - 1)!} \int_{(\mathcal{C}^{(d)})^{n - 1}} V_n(K, L_1, \dots, L_{n - 1}) \, \dd\mu^{(n - 1)}(L_1, \dots, L_{n - 1}) \bigg) .
	\end{equation*}
	Notice that $\kappa$ is well-defined as the term in the exponential is finite by the summability assumption on the $V_n$. \cite{BHLV:2020} multiply $\kappa$ with an indicator $\mathds{1}\big\{ \mu( \{ K \} ) = 0 \big\}$, a modification that changes nothing about the following observations. The map $\kappa$ is apparently measurable and it satisfies the cocycle relation \eqref{cocycle assumption} \citep[by Exercise 4.3 of][and the symmetry of the $V_n$]{LP:2017}. Furthermore, $\kappa$ has a finite range of interaction as
	\begin{equation*}
		\kappa(K, \mu)
		= \kappa\big( K, \mu_{B_H(K, R_V)} \big)
	\end{equation*}
	by assumption on $V$. Finally, if we assume that $\kappa$ is locally stable, meaning that
	\begin{equation*}
		\sup_{\mu \in \NN(\mathcal{C}^{(d)})} \kappa(K, \mu)
		\leq \vartheta(K)
	\end{equation*}
	for $\lambda$-a.e.\ $K \in \mathcal{C}^{(d)}$ and some locally $\lambda$-integrable function $\vartheta : \mathcal{C}^{(d)} \to [0, \infty)$, then Corollary \ref{COR existence result with local stability and finite range} guarantees the existence of a Gibbs particle process with PI $\kappa$. In fact, it is assumed by \cite{BHLV:2020} that $\kappa \leq 1$ uniformly.
\end{example}

In the context of particle processes there is a natural notion of clusters. In the language of Section \ref{SEC Gibbs point processes with cluster-dependent interaction} consider on $\mathcal{C}^*(\mathbb{B})$ the relation given through
\begin{equation} \label{relation for compact sets}
	K \sim L \quad \text{if, and only if,} \quad K \cap L \neq \varnothing .
\end{equation}
This relation is clearly symmetric and satisfies $\{ (K, L) \in \mathcal{C}^*(\mathbb{B})^2 : K \sim L \} \in \mathcal{B}(\mathcal{C}^*(\mathbb{B}))^{\otimes 2}$. It also leads to a very intuitive interpretation of the clusters in \eqref{mu-cluster of x} which, in this case, describe the connected components of a germ-grain model. Notice that the assumption in Theorem \ref{THM existence and uniqueness for subcritical cluster-dependent interactions} then simply describes absence of percolation in the Poisson-Boolean model. In the remainder of this section we use results on percolation in the Boolean model due to \cite{G:2008} to establish existence and uniqueness of corresponding Gibbs particle processes in the subcritical phase. As such we recover the uniqueness results of \cite{Ht:2019} and \cite{HtH:2019}, which is a benefit as these works contain technical gaps that are now filled by \cite{LO:2021} and the paper at hand. We also extend those results to particle processes with more general grains, similar to \cite{BHLV:2020} but allowing for unbounded grains. Moreover, we newly provide the corresponding existence results for those same processes in the region of their uniqueness.

We consider the setting from Example \ref{EXA Gibbs particle processes}, thus restricting to $\mathbb{B} = \R^d$. For a center function let $c : \mathcal{C}^{(d)} \to \R^d$ denote the map that assign to each compact set in $\R^d$ the center of its circumball, and let $\mathrm{rad} : \mathcal{C}^{(d)} \to [0, \infty)$ denote the map that assigns to a compact set the radius of its circumball. Both of these maps are continuous with respect to the Hausdorff metric. With the relation in \eqref{relation for compact sets} define clusters as in \eqref{mu-cluster of x}.

\begin{corollary} \label{COR existence and uniqueness of Gibbs particle processes in the subcritical phase}
	Let $\mathbb{Q}$ be a probability measure on $\mathcal{C}^{(d)}$ with $\mathbb{Q}(\mathcal{C}^{(d)}_0) = 1$ and
	\begin{equation*}
		\int_{\mathcal{C}^{(d)}} \mathrm{rad}(L)^d \, \dd \mathbb{Q}(L)
		< \infty ,
	\end{equation*}
	and, for $z > 0$, put
	\begin{equation*}
		\lambda_z(\cdot)
		= z \int_{\R^d} \int_{\mathcal{C}^{(d)}} \mathds{1}\{ L + x \in \cdot \, \} \, \dd\mathbb{Q}(L) \, \dd x .
	\end{equation*}
	Let $\kappa : \mathcal{C}^{(d)} \times \NN(\mathcal{C}^{(d)}) \to [0, \infty)$ be a measurable map that satisfies the cocycle property \eqref{cocycle assumption} and which is such that, for any $K \in \mathcal{C}^{(d)}$ and $\mu \in \NN(\mathcal{C}^{(d)})$,
	\begin{equation*}
		\kappa(K, \mu) \leq 1 
		\quad \text{as well as} \quad
		\kappa(K, \mu)
		= \kappa\big( K, C(K, \mu) \big) .
	\end{equation*}
	Then there is a constant $z_c = z_c(\mathbb{Q}, d) > 0$ such that for any $z < z_c$ there exists, up to equality in distribution, exactly one Gibbs particle process with PI $\kappa$, activity $z$, and grain distribution $\mathbb{Q}$, that is, exactly one Gibbs process in $\mathcal{C}^{(d)}$ with PI $\kappa$ and reference measure $\lambda_z$.
\end{corollary}
\begin{prf}
	Let $Z$ be a random element of $\mathcal{C}^{(d)}$ with distribution $\mathbb{Q}$ (the typical grain). We first note that the measure $\lambda_z$ is locally finite. Indeed, it follows as in Example \ref{EXA Gibbs particle processes} that, for any $K \in \mathcal{C}^{(d)}$ and $r > 0$,
	\begin{equation*}
		\lambda_z\big( B_H(K, r) \big)
		\leq z \int_{\mathcal{C}^{(d)}} \mathcal{L}^d\big( L + (-C) \big) \, \dd\mathbb{Q}(L)
		= z \cdot \E\big[ \mathcal{L}^d\big( Z + (-C) \big) \big]
	\end{equation*}
	for some suitable set $C \in \mathcal{C}^{(d)}$. As $C$ is compact, we have $C \subset B(0, s)$ for $s > 0$ large enough, and we also have $Z \subset B\big( 0, \mathrm{rad}(Z) \big)$ almost surely. Thus, by assumption on $\mathbb{Q}$,
	\begin{align*}
		\lambda_z\big( B_H(K, r) \big)
		\leq z \cdot \E\big[ \mathcal{L}^d\big( B(0, \mathrm{rad}(Z) + s) \big) \big]
		&= z \cdot \mathcal{L}^d\big( B(0, 1) \big) \cdot \E\big[ \big( \mathrm{rad}(Z) + s \big)^d \big] \\
		&\leq z \cdot \mathcal{L}^d\big( B(0, 1) \big) \cdot (1 + s)^d \cdot \max\big\{ 1, \, \E\big[ \mathrm{rad}(Z)^d \big] \big\} \\
		&< \infty .
	\end{align*}
	Let $\Phi = \sum_{j = 1}^\infty \delta_{K_j}$ be a Poisson particle process with intensity measure $\lambda_z$. Theorem \ref{THM existence and uniqueness for subcritical cluster-dependent interactions} yields the claim if we can show that the Boolean model based on $\Phi$ (in other words, the Boolean model with intensity $z$ and grain distribution $\mathbb{Q}$) is subcritical. To this end, let us define on $[0, \infty)$ the probability measure
	\begin{equation*}
		m(\cdot) = \int_{\mathcal{C}^{(d)}} \mathds{1}\big\{ \mathrm{rad}(L) \in \cdot \, \big\} \, \dd \mathbb{Q}(L)
		= \PP\big( \mathrm{rad}(Z) \in \cdot \, \big) ,
	\end{equation*}
	which satisfies
	\begin{equation*}
		\int_0^\infty r^d \, \dd m(r)
		= \int_{\mathcal{C}^{(d)}} \mathrm{rad}(L)^d \, \dd \mathbb{Q}(L)
		< \infty .
	\end{equation*}
	The Boolean model $\bigcup_{j = 1}^\infty K_j$ is almost surely contained in $\bigcup_{j = 1}^\infty B\big( c(K_j), \mathrm{rad}(K_j) \big)$, which is itself a Boolean model, since
	\begin{equation*}
		\Psi
		= \int_{\mathcal{C}^{(d)}} \mathds{1}\big\{ B\big( c(L), \mathrm{rad}(L) \big) \in \cdot \, \big\} \, \dd\Phi(L)
	\end{equation*}
	is a (stationary) Poisson particle process with intensity measure 
	\begin{align*}
		\int_{\mathcal{C}^{(d)}} \mathds{1}\big\{ B\big( c(L), \mathrm{rad}(L) \big) \in \cdot \, \big\} \, \dd\lambda_z(L)
		&= z \int_{\R^d} \int_{\mathcal{C}^{(d)}} \mathds{1}\big\{ B\big( x, \mathrm{rad}(L) \big) \in \cdot \, \big\} \, \dd\mathbb{Q}(L) \, \dd x \\
		&= z \int_{\R^d} \int_0^\infty \mathds{1}\big\{ B(x, r) \in \cdot \, \big\} \, \dd m(r) \, \dd x
	\end{align*}
	by the mapping theorem for Poisson processes \cite[Theorem 5.1 of][]{LP:2017}. Theorem 2.1 of \cite{G:2008} provides a constant $z_c(\mathbb{Q}, d) > 0$ such that
	\begin{equation*}
		\PP\big( C(K, \Phi)(\X) < \infty \big) \geq \PP\big( C(K, \Psi)(\X) < \infty \big) = 1
	\end{equation*}
	for all $K \in \mathcal{C}^{(d)}$ and every $z < z_c(\mathbb{Q}, d)$.
\end{prf}

\vspace{3mm}

To extend an example from the literature, we consider segment processes in $\R^d$, cf. Example 2.2 of \cite{FB:2018}. Proceeding as in Example \ref{EXA Gibbs particle processes} gives an existence result which covers the processes discussed by \cite{FB:2018} (even in arbitrary dimension), where a global and deterministic bound on the length of the segments is assumed. We use Corollary \ref{COR existence and uniqueness of Gibbs particle processes in the subcritical phase} to provide an existence and uniqueness result for segment processes which also allow for unbounded length distributions.
\begin{example}[\bf Segment processes] \label{EXA segment processes}
	Let $m$ be a probability measure on $[0, \infty)$ with existing $d$-th moment, which yields (half) the length of the segments, and let $Q$ be a probability measure on $\R^d$ concentrated on the unit sphere $S^{d - 1}$, which yields the orientation of the segments. Consider on $\mathcal{C}_0^{(d)}$ the probability measure
	\begin{equation*}
		\mathbb{Q}(\cdot)
		= \int_0^\infty \int_{S^{d - 1}} \mathds{1}\Big\{ \big\{ s \cdot v : s \in [-r, r] \big\} \in \cdot \, \Big\} \, \dd Q(v) \, \dd m(r) .
	\end{equation*}
	Let $V : \mathcal{C}^{(d)} \cup \{ \varnothing \} \to [0, \infty]$ be measurable with $V(\varnothing) = 0$ (where the $\sigma$-field on $\mathcal{C}^{(d)} \cup \{ \varnothing \}$ is constructed from the one on $\mathcal{C}^{(d)}$ by adding the singleton $\{ \varnothing \}$ as a measurable set). With the PI
	\begin{equation*}
		\kappa(K, \mu)
		= \exp\bigg( - \beta \int_{\mathcal{C}^{(d)}} V(K \cap L) \, \dd \mu(L) \bigg) ,
	\end{equation*}
	where $\beta > 0$, this fits into the setting of Corollary \ref{COR existence and uniqueness of Gibbs particle processes in the subcritical phase}. A specific example is $V(K) = c \cdot \mathds{1}\{ K \neq \varnothing \}$ for $c \in [0, \infty]$, but $V$ could also be the restriction of a locally finite measure on $\R^d$ onto $\mathcal{C}^{(d)} \cup \{ \varnothing \}$ \cite[which gives a measurable map by Proposition E.13 of][]{M:2017}. The existence result for this particular $\kappa$ can be improved by noting that it corresponds to a pair interaction and using Theorem \ref{THM existence result in the pair potential setting} or Theorem B.1 of \cite{J:2019}. But of course more complicated $\kappa$ can be considered as well.
\end{example}

To conclude this section on particle processes we leave as a take-away message that Gibbs particle processes exist as soon as the PI is bounded and has either finite range or depends only on suitable clusters with the underlying intensity parameter being small enough. Particle processes thus provide one large class of examples that show the usefulness of the generality of Corollary \ref{COR existence result with local stability and finite range} and Theorem \ref{THM existence and uniqueness for subcritical cluster-dependent interactions}, and thus, by extension, Theorem \ref{THM abstract existence result}. In particular, the possibility of considering abstract product spaces allows for interesting constructions.

\appendix

\section{Factorial measures}
\label{Appendix factorial measures}

We recall the definition of factorial measures on an arbitrary measurable space. The generality of the construction is due to \cite{LP:2017}. We mostly collect their results, but add some new (and mostly technical) insights. Before we start, we mention that a very slight technical omission is made by \cite{LP:2017}: to guarantee the measurability of the term in their equation (A.15), the localizing structure on the measurable space is explicitly needed, as laid out in Lemma \ref{LEMMA measurability of factorial-measure-integrals} below. Fix a localized measurable space $(\X, \mathcal{X})$ and let the space $\NN$ be defined as in Section \ref{SEC Introduction}. We write $[k] = \{ 1, \dots, k \}$ for $k \in \N$, $[0] = \varnothing$, and $[k] = \N$ if $k = \infty$. For a measure $\mu = \sum_{j = 1}^{k} \delta_{x_j} \in \NN$ with $k \in \N_0 \cup \{ \infty \}$ and $x_j \in \X$, define the $m$-th \textit{factorial measure} of $\mu$ on $(\X^m, \mathcal{X}^{\otimes m})$ as
\begin{equation} \label{fact. measure for sum of Diracs}
	\mu^{(m)} = \sideset{}{^{\neq}}\sum_{j_1, \dots, j_m \in [k]} \delta_{(x_{j_1}, \dots, x_{j_m})} ,
\end{equation}
where the superscript $\neq$ is used to indicate that the indices in the summation are pairwise distinct, and where the term is defined as $\mathbf{0} \in \NN(\X^m)$ if the sum is empty. Clearly, $\mu^{(m)} \in \NN(\X^m)$ and $\mu^{(1)} = \mu$. It is easy to verify that
\begin{equation} \label{def. equation factorial measures}
	\mu^{(m + 1)}(\cdot)
	= \int_{\X^m} \bigg( \int_\X \mathds{1}\big\{ (x_1, \dots, x_{m+1}) \in \cdot \, \big\} \, \dd \mu(x_{m+1}) - \sum_{j = 1}^{m} \mathds{1}\big\{ (x_1, \dots, x_m, x_j) \in \cdot \, \big\} \bigg) \dd \mu^{(m)}(x_1, \dots, x_m) 
\end{equation}
for each $m \in \N$. It is well known that in general measurable spaces not every measure in $\NN$ can be written as a sum of Dirac measures. Still, for each measure $\mu \in \NN$ there exists a unique sequence of symmetric measures $\mu^{(m)} \in \NN(\X^m)$, $m \in \N$, with $\mu^{(1)} = \mu$ and such that \eqref{def. equation factorial measures} is valid for all $m \in \N$. This is guaranteed by the following result, stated as Proposition 4.3 by \cite{LP:2017}.
\begin{prop} \label{PROP existence of factorial measures}
	For $\mu \in \NN$ there exists a unique sequence of symmetric measures $\mu^{(m)} \in \NN(\X^m)$ ($m \in \N$) such that $\mu^{(1)} = \mu$ and the recursion \eqref{def. equation factorial measures} is valid for all $m \in \N$. Moreover, the maps $\mu \mapsto \mu^{(m)} \in \NN(\X^m)$ are measurable.
\end{prop}
A proof of this proposition is given in Appendix A of \cite{LP:2017}. As we have mentioned before, though the authors do not state this explicitly, it is essential to have the localizing structure on $\X$ to obtain measurability of $\mu \mapsto \mu^{(m)}$. Proposition \ref{PROP existence of factorial measures}, and thus recursion \eqref{def. equation factorial measures}, is the definition of the factorial measures that we adopt. The following proposition collects those properties of factorial measures that can be found in Chapter 4 and Appendix A of \cite{LP:2017}. As before we write $\mu_B$ for the restriction of $\mu$ to a set $B$.
\begin{prop} \label{PROP properties of factorial measures due to Last/Penrose}
	Let $\mu, \nu \in \NN$ and fix any $m \in \N$. The following properties are satisfied.
	\begin{itemize}
		\item[(i)] If $D_1, \dots, D_m \in \mathcal{X}$ are pairwise disjoint, then $\mu^{(m)}(D_1 \times \dotso \times D_m) = \prod_{j = 1}^m \mu(D_j)$.
		
		\item[(ii)] For $B \in \mathcal{X}$ it holds that $\mu^{(m)}(B^m) = \mu(B) \cdot \big( \mu(B) - 1 \big) \cdot \dotso \cdot \big( \mu(B) - m + 1 \big)$.
		
		\item[(iii)] For $B \in \mathcal{X}$ it holds that $(\mu^{(m)})_{B^m} = \mu_B^{(m)}$.
		
		\item[(iv)] For $B \in \mathcal{X}$ the relation $\mu_B^{(m)}(B^m) = 0$ holds whenever $\mu(B) < m$.
		
		\item[(v)] If $\mu \leq \nu$, then $\mu^{(m)} \leq \nu^{(m)}$.
	\end{itemize}
\end{prop}

We complement the previous proposition by the following additional property.
\begin{lemma} \label{LEMMA properties of factorial measures newly added}
	Fix $m, k \in \N$ with $k \geq m$, and let $\mu \in \NN$. Let $B \in \mathcal{X}$ and $D \in \mathcal{X}^{\otimes m}$. If $\mu(B) = k$, then
		\begin{equation*}
			\mu_B^{(k)}(D \times B^{k - m}) = (k - m)! \cdot \mu_B^{(m)}(D) .
		\end{equation*}
\end{lemma}
\begin{prf}
	For $\mu \in \NN$ with $\mu(B) = k$ and $k \geq m$, the recursion \eqref{def. equation factorial measures} gives
	\begin{equation*}
		\mu_B^{(m + 1)}(D \times B)
		= \big( \mu(B) - m \big) \int_{\X^m} \mathds{1}\big\{ (x_1, \ldots, x_m) \in D \big\} \, \dd \mu_B^{(m)}(x_1, \ldots, x_m)
		= (k - m) \cdot \mu_B^{(m)}(D)
	\end{equation*}
	and, by applying this relation iteratively $k - m$ times, we obtain $\mu_B^{(k)}(D \times B^{k - m}) = (k - m)! \cdot \mu_B^{(m)}(D)$. 
\end{prf}

\vspace{3mm}

Already in stating the (multivariate) GNZ equations, the following measurability property is essential. We mostly use this result for $\mathbb{Y} = \NN$.
\begin{lemma} \label{LEMMA measurability of factorial-measure-integrals}
	Let $(\mathbb{Y}, \mathcal{Y})$ be an arbitrary measurable space and fix $m \in \N$. The mapping
	\begin{equation*}
		\NN \times \mathbb{Y} \ni (\mu, y) \mapsto
		\int_{\X^m} f(x_1, \dots, x_m, y) \, \dd\mu^{(m)}(x_1, \dots, x_m)
		\in [0, \infty]
	\end{equation*}
	is $(\mathcal{N} \otimes \mathcal{Y})$-measurable for every measurable function $f : \X^m \times \mathbb{Y} \to [0, \infty]$.
\end{lemma}
\begin{prf}
	By Proposition \ref{PROP existence of factorial measures}, the mapping $\mu \mapsto \mu^{(m)}$ is measurable, so for any $D \in \mathcal{X}^{\otimes m}$ and $A \in \mathcal{Y}$ the map
	\begin{equation*}
		(\mu, y) \mapsto
		\mathds{1}_A(y) \int_{\X^m} \mathds{1}_D \, \dd\mu^{(m)}
		= \mathds{1}_A(y) \cdot \mu^{(m)}(D)
	\end{equation*}
	is measurable. Denote by $\mathcal{D}$ the collection of all sets $E \in \mathcal{X}^{\otimes m} \otimes \mathcal{Y}$ for which
	\begin{equation*}
		(\mu, y) \mapsto
		\int_{\X^m} \mathds{1}_E(x_1, \dots, x_m, y) \, \dd\mu^{(m)}(x_1, \dots, x_m)
	\end{equation*}
	is measurable. The $\pi$-system $\{ D \times A : D \in \mathcal{X}^{\otimes m}, \, A \in \mathcal{Y} \}$ is contained in $\mathcal{D}$. Moreover, we have $\X^m \times \mathbb{Y} \in \mathcal{D}$ and $\mathcal{D}$ is clearly closed with respect to countable disjoint unions. Let $E, F \in \mathcal{D}$ with $E \subset F$, and observe that
	\begin{equation*}
		\int_{\X^m} \mathds{1}_{F \setminus E}(x_1, \dots, x_m, y) \, \dd\mu^{(m)}(x_1, \dots, x_m)
		= \lim_{\ell \to \infty} \int_{\X^m} \mathds{1}_{F \setminus E}(x_1, \dots, x_m, y) \, \dd\mu_{B_\ell}^{(m)}(x_1, \dots, x_m)
	\end{equation*}
	by monotone convergence, using that $\bigcup_{\ell = 1}^\infty B_\ell^m = \X^m$ and $(\mu^{(m)})_{B_\ell^m} = \mu_{B_\ell}^{(m)}$. Since $\mu_{B_\ell}^{(m)}(\X^m) \leq \mu(B_\ell)! < \infty$, the following difference of integrals is well defined (for each $\ell \in \N$)
	\begin{align*}
		\int_{\X^m} \mathds{1}_{F \setminus E}(x_1, \dots, x_m, y) \, \dd\mu_{B_\ell}^{(m)}(x_1, \dots, x_m)
		= &\int_{\X^m} \mathds{1}_{F}(x_1, \dots, x_m, y) \, \dd\mu_{B_\ell}^{(m)}(x_1, \dots, x_m) \\
		&- \int_{\X^m} \mathds{1}_{E}(x_1, \dots, x_m, y) \, \dd\mu_{B_\ell}^{(m)}(x_1, \dots, x_m) .
	\end{align*}
	The right hand side is a measurable function of $(\mu, y)$ since $E, F$ are from $\mathcal{D}$ and $\mu \mapsto \mu_{B_\ell}$ is measurable. As limits of measurable functions are measurable, we conclude that
	\begin{equation*}
		(\mu, y) \mapsto
		\int_{\X^m} \mathds{1}_{F \setminus E}(x_1, \dots, x_m, y) \, \dd\mu^{(m)}(x_1, \dots, x_m)
	\end{equation*}
	is measurable, so $F \setminus E \in \mathcal{D}$. Thus, $\mathcal{D}$ is a Dynkin system and Dynkin's $\pi$-$\lambda$-theorem implies $\mathcal{D} = \mathcal{X}^{\otimes m} \otimes \mathcal{Y}$. A standard monotone approximation completes the proof.
\end{prf}

\vspace{3mm}

Equation (4.19) of \cite{LP:2017} indicates how factorial measures can be used to represent any functional on $\NN$ when evaluated on $\NN_f$. Indeed, refining that particular equation and extending it by monotone approximation yields the following result.
\begin{prop} \label{PROP representation of NN to R functions via factorial measures}
	Let $F : \NN \to [0, \infty]$ be a measurable map. Then, for any $\mu \in \NN_f$,
	\begin{equation*}
		F(\mu)
		= \mathds{1}\big\{ \mu(\X) = 0 \big\} \cdot F(\mathbf{0}) + \sum_{m = 1}^\infty \frac{1}{m!} \, \mathds{1}\big\{ \mu(\X) = m \big\} \int_{\X^m} F\Big( \sum_{j = 1}^m \delta_{x_j} \Big) \, \dd\mu^{(m)}(x_1, \dots, x_m) .
	\end{equation*}
\end{prop}

We frequently use Proposition \ref{PROP representation of NN to R functions via factorial measures} to argue that a function which is finite for sums of finitely many Dirac measures is finite on the whole of $\NN_f$.
\begin{lemma} \label{LEMMA finite functions on Dirac measures are finite on NN_f}
	Let $G : \NN \to [0, \infty]$ be a measurable map such that $G(\mathbf{0}) < \infty$ and $G(\delta_{x_1} + \ldots + \delta_{x_m}) < \infty$ for all $x_1, \ldots, x_m \in \X$ and any $m \in \N$. Then $G(\mu) < \infty$ for each $\mu \in \NN_f$.
\end{lemma}
\begin{prf}
	Applying Proposition \ref{PROP representation of NN to R functions via factorial measures} to the map $F(\mu) = \mathds{1}\big\{ G(\mu) < \infty \big\}$ and using the assumption, we obtain, for any $\mu \in \NN_f$,
	\begin{align*}
		\mathds{1}\big\{ G(\mu) < \infty \big\}
		&= \mathds{1}\big\{ \mu(\X) = 0 \big\} \cdot \mathds{1}\big\{ G(\mathbf{0}) < \infty \big\} \\
		&\quad + \sum_{m = 1}^\infty \frac{1}{m!} \, \mathds{1}\big\{ \mu(\X) = m \big\} \int_{\X^m} \mathds{1}\bigg\{ G\Big( \sum_{j = 1}^m \delta_{x_j} \Big) < \infty \bigg\} \, \dd\mu^{(m)}(x_1, \ldots, x_m) \\
		&= \mathds{1}\big\{ \mu(\X) = 0 \big\} + \sum_{m = 1}^\infty \frac{1}{m!} \, \mathds{1}\big\{ \mu(\X) = m \big\} \, \mu^{(m)}(\X^m) \\
		&= \mathds{1}\big\{ \mu(\X) < \infty \big\} \\
		&= 1 .
	\end{align*}
\end{prf}

Exercise 4.3 of \cite{LP:2017} asks the reader to prove that, for any $\mu \in \NN$, $x \in \X$, and $m \in \N$,
\begin{equation*}
	(\mu + \delta_x)^{(m + 1)}
	= \mu^{(m + 1)} + \int_{\X^m} \Big[ \mathds{1}\big\{ (x, x_1, \dots, x_m) \in \cdot \, \big\} + \dotso + \mathds{1}\big\{ (x_1, \dots, x_m, x) \in \cdot \, \big\} \Big] \, \dd\mu^{(m)}(x_1, \dots, x_m) .
\end{equation*}
We formulate and prove the following lemma as an extension of this exercise. Therein, we denote by $S(n)$ the symmetric group of degree $n$ containing all permutations of $[n] = \{ 1, \dots, n \}$. Notice that, for $z_1, \dots, z_n \in \X$ and $\tau \in S(n)$, we write $\tau(z_1, \dots, z_n) = (z_{\tau(1)}, \dots, z_{\tau(n)})$.
\begin{lemma} \label{LEMMA general lemma on factorial measures}
	Let $\mu, \nu \in \NN$ with $\mu(\X) = k \in \N$ and $\nu(\X) = m \in \N$. Then
	\begin{equation*}
		(\mu + \nu)^{(k + m)}
		= \frac{1}{k! \cdot m!} \int_{\X^k} \int_{\X^m} \sum_{\tau \in S(k + m)} \mathds{1}\big\{ \tau(z_1, \dots, z_{k + m}) \in \cdot \, \big\} \, \dd\nu^{(m)}(z_{k + 1}, \dots, z_{k + m}) \, \dd\mu^{(k)}(z_{1}, \dots, z_{k})
	\end{equation*}
	and, for any measurable and symmetric map $f : \X^{k + m} \to [0, \infty]$,
	\begin{align*}
		&\int_{\X^{k + m}} f(z_1, \dots, z_{k + m}) \, \dd(\mu + \nu)^{(k + m)}(z_1, \dots, z_{k + m}) \\
		&\quad = \frac{(k + m)!}{k! \cdot m!} \int_{\X^k} \int_{\X^m} f(x_1, \dots, x_k, y_1, \dots, y_m) \, \dd\nu^{(m)}(y_1, \dots, y_m) \, \dd\mu^{(k)}(x_1, \dots, x_k) .
	\end{align*}
\end{lemma}
\begin{prf}
	We only show the first claim, with the second claim following by monotone approximation. Clearly, both $(\mu + \nu)^{(k + m)}$ and the term in the lemma constitute finite measures on $\X^{k + m}$. To prove that these measures are equal it suffices to show that they agree on sets of the form $C_1 \times \ldots \times C_{k + m} \in \mathcal{X}^{\otimes (k + m)}$. Thus, let $C_1, \ldots, C_{k + m} \in \mathcal{X}$. Denote by $\mathcal{A}$ the field generated by these sets. Lemma A.15 of \cite{LP:2017} yields the existence of $x_1, \ldots, x_k, y_1, \ldots, y_m \in \X$ such that $\mu^\prime = \sum_{i = 1}^k \delta_{x_i}$ and $\nu^\prime = \sum_{j = 1}^m \delta_{y_j}$ satisfy
	\begin{equation*}
		\mu^{(n)}(D) = (\mu^\prime)^{(n)}(D) \quad \text{and} \quad \nu^{(n)}(D) = (\nu^\prime)^{(n)}(D)
	\end{equation*}
	for all $n \in \N$ and $D \in \mathcal{A}^{\otimes n}$, where $\mathcal{A}^{\otimes n}$ is the field generated by the system $\{ D_1 \times \ldots \times D_n : D_1, \ldots, D_n \in \mathcal{A} \}$. Similarly, there exist $z_1, \ldots, z_{k + m} \in \X$ such that $(\mu + \nu)^\prime = \sum_{i = 1}^{k + m} \delta_{z_i}$ satisfies
	\begin{equation*}
		(\mu + \nu)^{(n)}(D) = \big( (\mu + \nu)^\prime \big)^{(n)}(D)
	\end{equation*}
	for all $n \in \N$ and $D \in \mathcal{A}^{\otimes n}$. Hence, \eqref{fact. measure for sum of Diracs} implies
	\begin{align*}
		(\mu + \nu)^{(k + m)}(C_1 \times \ldots \times C_{k + m})
		&= \big( (\mu + \nu)^\prime \big)^{(k + m)}(C_1 \times \ldots \times C_{k + m}) \\
		&= \sideset{}{^{\neq}}\sum_{i_1, \ldots, i_{k + m} \in [k + m]} \mathds{1}\Big\{ (z_{i_1}, \ldots, z_{i_{k + m}}) \in C_1 \times \ldots \times C_{k + m} \Big\} .
	\end{align*}	
	From the construction in Lemma A.15 of \cite{LP:2017} it is obvious that $(\mu + \nu)^\prime = \mu^\prime + \nu^\prime$, so the previous term equals
	\begin{equation*}
		\sideset{}{^{\neq}}\sum_{i_1, \ldots, i_{k} \in [k]} ~~\; \sideset{}{^{\neq}}\sum_{j_1, \ldots, j_{m} \in [m]} ~ \frac{1}{k! \cdot m!} \sum_{\tau \in S(k + m)} \mathds{1}\Big\{ \tau(x_{i_1}, \ldots, x_{i_k}, y_{j_1}, \ldots, y_{j_m}) \in C_1 \times \ldots \times C_{k + m} \Big\} ,
	\end{equation*}
	which, by \eqref{fact. measure for sum of Diracs} and the construction of $\mu'$ and $\nu'$, equals the right hand side of the claim.
\end{prf}

\section{Measurable diagonals} 
\label{Appendix measurable diagonals}

Let $(\X, \mathcal{X})$ be a localized measurable space and let $\NN$ be defined as in Section \ref{SEC Introduction}. We say that the measurable space $\X$ has a \textit{measurable diagonal} if
\begin{equation*}
	D_\X
	= \big\{ (x, y) \in \X \times \X : x = y \big\}
	\in \mathcal{X}^{\otimes 2} .
\end{equation*}

\begin{lemma} \label{LEMMA Appen. prelim. statement, implications of meas. diag.}
	Assume that $\X$ has a measurable diagonal. Then $\X$ is separable, that is, $\{ x \} \in \mathcal{X}$ for each $x \in \X$. Moreover, the map $\X \times \NN \ni (x, \mu) \mapsto \mu(\{x\}) \in \N_0$ is well-defined and measurable.
\end{lemma}
\begin{prf}
	For each $x \in \X$ the map $h_x : \X \to \X^2$, $h_x(y) = (x, y)$ is measurable, so $\{ x \} = h_x^{-1}(D_\X) \in \mathcal{X}$. As there exists some $\ell \in \N$ such that $x \in B_\ell$, we actually have $\{x\} \in \mathcal{X}_b$, so the definition of $\NN$ ensures that $\mu(\{x\}) < \infty$. In particular, the map in consideration is well-defined. Since $\X$ has a measurable diagonal, the mapping $(x, y) \mapsto \mathds{1}_{D_\X}(x, y)$ is measurable and Lemma \ref{LEMMA measurability of factorial-measure-integrals} implies the measurability of
	\begin{equation*}
		(x, \mu)
		\mapsto \int_\X \mathds{1}_{D_\X}(x, y) \, \dd \mu(y)
		= \int_\X \mathds{1}\{ x = y \} \, \dd \mu(y)
		= \int_\X \mathds{1}_{\{ x \}}(y) \, \dd \mu(y)
		= \mu(\{ x \}) .
	\end{equation*}
\end{prf}

The measurable space $(\X, \mathcal{X})$ is called \textit{countably generated} if the $\sigma$-field $\mathcal{X}$ is generated by countably many sets from $\mathcal{X}$.
\begin{prop} \label{PROP Appen. prelim. statement, measurable diagonal}
	If $(\X, \mathcal{X})$ is countably generated and separable, then $\X$ has a measurable diagonal.
\end{prop}

A short proof of this result can be given from scratch, but we refer to \cite{M:1980} for a stronger version, where a slightly weaker assumption than the measurable space being countably generated actually gives equivalence. For our purposes the above statement is well-suited since the substandard Borel spaces we consider (see Appendix \ref{Appendix Kolmog. extension result on N} and Sections \ref{SEC prelinimaries}, \ref{SEC existence result Gibbs processes}) are countably generated. Moreover, notice that Proposition \ref{PROP Appen. prelim. statement, measurable diagonal} includes every second countable Hausdorff space and thus in particular every separable metric space. It also includes every Borel space (as \citealp{K:2017}, and \citealp{LP:2017}, understand them).

For a measure $\mu$ on $(\X, \mathcal{X})$ and $x \in \X$ define
\begin{equation*}
	\mu \setminus \delta_x
	= \mu - \delta_x \, \mathds{1}\big\{ \mu(\{ x \}) > 0 \big\} ,
\end{equation*}
which is itself a measure on $(\X, \mathcal{X})$. The measure $\mu \setminus \delta_{x_1} \setminus \dotso \setminus \delta_{x_m}$, for $x_1, \dots, x_m \in \X$ and $m \in \N$, is defined iteratively.

\begin{lemma} \label{LEMMA Appen. meas. of reducing a counting measure}
	Assume that $\X$ has a measurable diagonal. Then the map $d_m : \X^m \times \NN \to \NN$,
	\begin{equation*}
		d_m(x_1, \dots, x_m, \mu)
		= \mu \setminus \delta_{x_1} \setminus \dotso \setminus \delta_{x_m}
	\end{equation*}
	is measurable for each $m \in \N$.
\end{lemma}
\begin{prf}
	We prove the claim by induction, starting with the initial case $m = 1$. Note that, for any $B \in \mathcal{X}$,
	\begin{equation*}
		(x, \mu) \mapsto d_1(x, \mu)(B) = (\mu \setminus \delta_x) (B)
		= \mu(B) - \mathds{1}\{ x \in B \} \cdot \mathds{1}\big\{ \mu(\{ x \}) > 0 \big\}
		\in \N_0 \cup \{\infty\}
	\end{equation*}
	is measurable by Lemma \ref{LEMMA Appen. prelim. statement, implications of meas. diag.}. If the map $d_m$ is measurable for any fixed $m \in \N$, then $d_{m + 1}$ is also measurable since $d_{m + 1}(x_1, \ldots, x_{m + 1}, \mu) = d_1\big( x_{m + 1}, d_m(x_1, \ldots, x_m, \mu) \big)$.
\end{prf}

\section{Local events and local functions}
\label{Appendix local events and functions}

In this part of the appendix we deal with local events and functions. For a set $B \in \mathcal{X}$ define the map $p_B : \NN \to \NN$, $p_B(\mu) = \mu_B$, which assigns counting measures their restriction onto $B$. Put
\begin{equation*}
	\mathcal{N}_B = \sigma(p_B) = p_B^{-1}(\mathcal{N}),
\end{equation*}
the $\sigma$-field on $\NN$ generated by $p_B$. Moreover, define $\mathcal{Z} = \bigcup_{B \in \mathcal{X}_b} \mathcal{N}_B$.
\begin{definition}[\bf Local events and local functions] \label{DEF local events and functions}
	~\begin{itemize}
		\item A measurable set $A \in \mathcal{N}$ is called a $B$\textit{-local event} (for some $B \in \mathcal{X}$) if $A \in \mathcal{N}_B$. We say that $A \in \mathcal{N}$ is a \textit{local event} if there exists a set $B \in \mathcal{X}_b$ such that $A$ is $B$-local, that is, if $A \in \mathcal{Z}$.
		
		\item An $\mathcal{N}$-measurable map $F : \NN \to [- \infty, \infty]$ is called $B$\textit{-local} (for some $B \in \mathcal{X}$) if $F(\mu) = F(\mu_B)$ for all $\mu \in \NN$. The map $F$ is called a \textit{local function} if there exists a set $B \in \mathcal{X}_b$ such that $F$ is $B$-local.
	\end{itemize}
\end{definition}

Notice that we define $B$-locality for arbitrary sets $B \in \mathcal{X}$, but when we call an event or a function local, we specify to bounded sets $B \in \mathcal{X}_b$. The following properties are either obvious or easy exercises.
\begin{lemma} \label{LEMMA properties of local events and functions}
	In parts \textit{(i)}--\textit{(iv)}, fix a set $B \in \mathcal{X}$.
	\begin{itemize}
		\item[(i)] The $\sigma$-field $\mathcal{N}_B$ is generated by the evaluation maps $\pi_D : \mu \mapsto \mu(D)$, $D \in \mathcal{X}$ with $D \subset B$. In particular, $\mathcal{N}_B \subset \mathcal{N}$.
		
		\item[(ii)] For $A \in \mathcal{N}_B$ and $\mu \in \NN$, it is true that $\mu \in A$ if, and only if, $\mu_B \in A$.
		
		\item[(iii)] An $\mathcal{N}$-measurable map $F : \NN \to [- \infty, \infty]$ is $B$-local if, and only if, $F$ is $\mathcal{N}_B$-measurable.
		
		\item[(iv)] For each $B_1, B_2 \in \mathcal{X}$, it holds that $\mathcal{N}_{B_1 \cup B_2} = \sigma(\mathcal{N}_{B_1} \cup \mathcal{N}_{B_2})$.
		
		\item[(v)] The collection $\mathcal{Z} = \bigcup_{B \in \mathcal{X}_b} \mathcal{N}_B$ is an algebra of subsets of $\NN$ with $\sigma(\mathcal{Z}) = \mathcal{N}$.
	\end{itemize}
\end{lemma}

Similar to part \textit{(v)}, $\bigcup_{j = 1}^\infty \mathcal{N}_{B_j}$ is an algebra with $\sigma\big( \bigcup_{j = 1}^\infty \mathcal{N}_{B_j} \big) = \mathcal{N}$.

\section{General facts about Janossy and factorial moment measures}
\label{Appendix facts about Jan. and Cor. meas.}

We first recall the definition of factorial moment and Janossy measures, following Chapter 4 of \cite{LP:2017}, and then discuss basic properties and their mutual relations, generalizing results from Chapter 5.4 of \cite{DV:2005}.

Let $\eta$ be a point process in a localized measurable space $(\X, \mathcal{X})$ and let $m\in \N$. The $m$\textit{-th factorial moment measure} of $\eta$ is $\alpha_{\eta, m}(\cdot) = \E \big[ \eta^{(m)}(\cdot) \big]$, which is a measure on $(\X^m, \mathcal{X}^{\otimes m})$. The $m$-th factorial moment measure of a Poisson process with intensity measure $\lambda$ is simply $\lambda^m$ and the Poisson process is the only point process which satisfies this property for every $m \in \N$. If the factorial moment measure $\alpha_{\eta, m}$ of a point process $\eta$ is absolutely continuous with respect to $\lambda^m$ with Radon-Nikodym density $\rho_{\eta, m}$, then $\rho_{\eta, m}$ is called \textit{correlation function} of order $m$ (of $\eta$ with respect to $\lambda$).

The \textit{Janossy measure} of order $m \in \N$ of $\eta$ restricted to $B$ is the measure $J_{\eta, B, m}$ on the product space $(\X^m, \mathcal{X}^{\otimes m})$ defined as
\begin{equation*}
	J_{\eta, B, m}(\cdot) = \frac{1}{m!} \, \E \Big[ \mathds{1}\big\{ \eta(B) = m \big\} \, \eta_B^{(m)}(\cdot) \Big] .
\end{equation*}
The Janossy measures are symmetric and satisfy
\begin{equation*}
	J_{\eta, B, m}(\X^m) = \PP\big( \eta(B) = m \big) .
\end{equation*}
In line with this last observation we put $J_{\eta, B, 0} = \PP\big( \eta(B) = 0 \big)$. With our choice of the space $\NN$, Theorem 4.7 of \cite{LP:2017} states that if $\eta$ and $\eta^\prime$ are two point processes with $J_{\eta, B, m} = J_{\eta^\prime, B, m}$ for all $m \in \N_0$ and some set $B \in \mathcal{X}_b$, then $\eta_B$ and $\eta^\prime_B$ have the same distribution. By Example 4.8 of \cite{LP:2017}, the Janossy measures of a Poisson process with intensity measure $\lambda$ are
\begin{equation*}
	J_{B, m} = \frac{\mathrm{e}^{- \lambda(B)}}{m!} \, \lambda_B^m .
\end{equation*}
If, for fixed $B \in \mathcal{X}$ and $m \in \N$, the Janossy measure $J_{\eta, B, m}$ of some point process $\eta$ is absolutely continuous with respect to $\lambda_B^m$, then the density function $j_{\eta, B, m}$ is called \textit{Janossy density}.

With Proposition \ref{PROP representation of NN to R functions via factorial measures} it is possible to express the expectation of functionals of $\eta_B$ via the Janossy measures.
\begin{lemma} \label{LEMMA expec. repres. via Janossy measures}
	Let $\eta$ be a point process in $\X$, and fix $B \in \mathcal{X}$. Then
	\begin{equation*}
		\E \Big[ F(\eta_B) \cdot \mathds{1}\big\{ \eta(B) < \infty \big\} \Big]
		= F(\mathbf{0}) \cdot \PP\big( \eta(B) = 0 \big) + \sum_{m = 1}^\infty \int_{\X^m} F\Big( \sum_{j = 1}^m \delta_{x_j} \Big) \, \dd J_{\eta, B, m}(x_1, \dots, x_m)
	\end{equation*}
	for every measurable map $F : \NN \to [0, \infty]$.
\end{lemma}

The statement of Exercise 3.7 from \cite{LP:2017} is a simple corollary of the previous lemma.
\begin{corollary} \label{COR Janossy representation of Poisson processes}
	Let $\Phi$ be a Poisson process in $\X$ with intensity measure $\lambda$. Then, for any set $B \in \mathcal{X}$ with $\lambda(B) < \infty$ and every measurable map $F : \NN \to [0, \infty]$,
	\begin{equation*}
		\E\big[ F(\Phi_B) \big]
		= \mathrm{e}^{- \lambda(B)} F(\mathbf{0}) + \mathrm{e}^{- \lambda(B)} \sum_{m = 1}^\infty \frac{1}{m!} \int_{B^m} F\Big( \sum_{j = 1}^m \delta_{x_j} \Big) \, \dd\lambda^m(x_1, \dots, x_m) .
	\end{equation*}
\end{corollary}

With the formalism of local events and functions from Appendix \ref{Appendix local events and functions}, a converse of Lemma \ref{LEMMA expec. repres. via Janossy measures} reads as follows.
\begin{lemma} \label{LEMMA characterization of Janossy measures}
	Let $\eta$ be a point process in $\X$ and let $B \in \mathcal{X}$. Assume there exists a collection of symmetric measures $J_{B, m}$ on $(\X^m, \mathcal{X}^{\otimes m})$ which vanish outside $B^m$, for each $m \in \N$, and $J_{B, 0} \in [0, \infty)$ such that, for all $A \in \mathcal{N}_B$,
	\begin{equation*}
		\PP(\eta \in A)
		= \mathds{1}_A(\mathbf{0}) \, J_{B, 0} + \sum_{m = 1}^\infty \int_{B^m} \mathds{1}_A\Big( \sum_{j = 1}^m \delta_{x_j} \Big) \, \dd J_{B, m}(x_1, \dots, x_m) .
	\end{equation*}
	Then the $(J_{B, m})_{m \in \N_0}$ are the Janossy measures of $\eta$ restricted to $B$, and $\PP\big( \eta(B) < \infty \big) = 1$.
\end{lemma}
\begin{prf}
	First, choose $A = \big\{ \mu \in \NN : \mu(B) = 0 \big\} = \pi_B^{-1}(\{ 0 \}) \in \mathcal{N}_B$ and notice that
	\begin{equation*}
		\PP\big( \eta(B) = 0 \big)
		= \PP( \eta \in A )
		= J_{B, 0} ,
	\end{equation*}
	so $J_{B, 0}$ is the Janossy measure of order $0$ of $\eta$. Now, let $F : \NN \to [0, \infty]$ be a measurable function. By part \textit{(iii)} of Lemma \ref{LEMMA properties of local events and functions}, the map $F \circ p_B$ is $\mathcal{N}_B$-measurable. Hence, by the assumption and monotone approximation, we have
	\begin{equation*}
		\E \big[ F(\eta_B) \big]
		= F(\mathbf{0}) \, J_{B, 0} + \sum_{m = 1}^\infty \int_{B^m} F\Big( \sum_{j = 1}^m \delta_{x_j} \Big) \, \dd J_{B, m}(x_1, \dots, x_m) .
	\end{equation*}
	Let $k \in \N$ and $D \in \mathcal{X}^{\otimes k}$. With the specific choice $F(\mu) = \frac{1}{k!} \, \mathds{1}\big\{ \mu(B) = k \big\} \, \mu_B^{(k)}(D)$ we obtain
	\begin{align*}
		\frac{1}{k!} \, \E \Big[ \mathds{1}\big\{ \eta(B) = k \big\} \, \eta_B^{(k)}(D) \Big]
		= \E \big[ F(\eta_B) \big]
		&= \frac{1}{k!} \int_{B^k} \mathds{1}\bigg\{ \sum_{j = 1}^k \delta_{x_j}(B) = k \bigg\} \, \Big( \sum_{j = 1}^k \delta_{x_j} \Big)_B^{(k)}(D) \, \dd J_{B, k}(x_1, \ldots, x_k) \\
		&= \frac{1}{k!} \quad \sideset{}{^{\neq}}\sum_{j_1, \ldots, j_k \in [k]} ~ \int_{B^k} \mathds{1}_D(x_{j_1}, \ldots, x_{j_k}) \, \dd J_{B, k}(x_1, \ldots, x_k) .
	\end{align*}
	By the symmetry of $J_{B, k}$, the right hand side equals
	\begin{equation*}
		\frac{1}{k!} \quad \sideset{}{^{\neq}}\sum_{j_1, \ldots, j_k \in [k]} \, J_{B, k}(D \cap B^k)
		= J_{B, k}(D \cap B^k)
		= J_{B, k}(D) .
	\end{equation*}
	We have thus verified that $J_{B, k}$ is the Janossy measure of order $k$ of $\eta_B$. Finally, choose $A = \big\{ \mu \in \NN : \mu(B) = \infty \big\} \in \mathcal{N}_B$ to conclude that $\PP\big( \eta(B) = \infty \big) = \PP(\eta \in A) = 0$.
\end{prf}

\vspace{3mm}

In the remainder of this appendix section we discuss the connections between factorial moment and Janossy measures. The following two theorems are also stated in principle in Chapter 5.4 of \cite{DV:2005}, but we generalize the results to potentially infinite point processes on arbitrary measurable spaces. We also translate them into up-to-date notation and give elegant proofs using the properties of abstract factorial measures. We start by showing that the factorial moment measures can be expressed locally via the Janossy measures without any additional assumptions.
\begin{theorem} \label{THM repres. of factorial mom. meas. via Janossy meas.}
	Let $\eta$ be a point process in $\X$. Fix $B \in \mathcal{X}$ with $\PP\big( \eta(B) < \infty \big) = 1$. For each $m \in \N$  and every measurable function $f : \X^m \to [0, \infty]$ with $f = 0$ on $\X^m \setminus B^m$, it holds that
	\begin{equation*}
		\int_{\X^m} f \, \dd \alpha_{\eta, m}
		= \sum_{k = m}^\infty \frac{k!}{(k - m)!} \int_{\X^m} f(x_1, \dots, x_m) \, J_{\eta, B, k}\big( \dd(x_1, \dots, x_m) \times B^{k - m} \big) .
	\end{equation*}
\end{theorem}
\begin{prf}
	Fix $m \in \N$. For $D \in \mathcal{X}^{\otimes m}$ with $D \subset B^m$, we have
	\begin{align} \label{formula factorial mom. meas. via Janossy meas.}
		\alpha_{\eta, m}(D)
		= \E \big[ \eta_B^{(m)}(D) \big]
		= \sum_{k = m}^\infty \E \Big[ \mathds{1}\big\{ \eta(B) = k \big\} \, \eta_B^{(m)}(D) \Big]
		&= \sum_{k = m}^\infty \frac{1}{(k - m)!} \, \E \Big[ \mathds{1}\big\{ \eta(B) = k \big\} \, \eta_B^{(k)}(D \times B^{k - m}) \Big] \notag \\
		&= \sum_{k = m}^\infty \frac{k!}{(k - m)!} \, J_{\eta, B, k}(D \times B^{k - m}) ,
	\end{align}
	by Proposition \ref{PROP properties of factorial measures due to Last/Penrose} and Lemma \ref{LEMMA properties of factorial measures newly added}. The theorem follows from equation \eqref{formula factorial mom. meas. via Janossy meas.} and standard monotone approximation by step functions.
\end{prf}

\begin{corollary} \label{COR repres. of correlation func. via Janossy densities}
	Let $\eta$ be a point process in $\X$. Assume that the Janossy measures $J_{\eta, B, k}$ of $\eta$ admit density functions $j_{\eta, B, k}$ with respect to $\lambda^k$, for each $k \in \N$ and $B \in \mathcal{X}_b$. Then the correlation functions of $\eta$ exist and satisfy
	\begin{equation*}
		\rho_{\eta, m}(x_1, \dots, x_m)
		= \sum_{k = m}^\infty \frac{k!}{(k - m)!} \int_{B^{k - m}} j_{\eta, B, k}(x_1, \dots, x_k) \, \dd \lambda^{k - m}(x_{m+1}, \dots, x_k)
	\end{equation*}
	for $\lambda^m$-a.e.\ $(x_1, \dots, x_m) \in B^m$, all $m \in \N$, and each $B \in \mathcal{X}_b$.
\end{corollary}
\begin{prf}
	Let $B \in \mathcal{X}_b$ and let $f : \X^m \to [0, \infty]$ be a measurable function such that $f = 0$ on $\X^m \setminus B^m$. By Theorem \ref{THM repres. of factorial mom. meas. via Janossy meas.} we have
	\begin{align*}
		\int_{\X^m} f \, \dd \alpha_{\eta, m}
		&= \sum_{k = m}^\infty \frac{k!}{(k - m)!} \int_{\X^k} f(x_1, \dots, x_m) \, \mathds{1}_{B^{k - m}}(x_{m+1}, \dots, x_k) \, j_{\eta, B, k}(x_1, \dots, x_k) \, \dd \lambda^k(x_1, \dots, x_k) \\
		&= \int_{\X^m} f(x_1, \dots, x_m) \sum_{k = m}^\infty \frac{k!}{(k - m)!} \int_{B^{k - m}} j_{\eta, B, k}(x_1, \dots, x_k) \, \dd \lambda^{k - m}(x_{m+1}, \dots, x_k) \, \dd \lambda^m(x_1, \dots, x_m) .
	\end{align*}
	We conclude that $(\alpha_{\eta, m})_{B^m}$ is absolutely continuous with respect to $\lambda_B^m$ for all $B \in \mathcal{X}_b$. Now, let $D \in \mathcal{X}^{\otimes m}$ with $\lambda^m(D) = 0$. Then $\lambda_{B_\ell}^m(D) = \lambda^m(D \cap B_\ell^m) = 0$ for each $\ell \in \N$, and therefore
	\begin{equation*}
		\alpha_{\eta, m}(D)
		= \lim_{\ell \to \infty} \alpha_{\eta, m}(D \cap B_\ell^m)
		= \lim_{\ell \to \infty} (\alpha_{\eta, m})_{B_\ell^m} (D)
		= 0 .
	\end{equation*}
	Thus, $\alpha_{\eta, m}$ is absolutely continuous with respect to $\lambda^m$, so the correlation functions of $\eta$ exist and satisfy the claim.
\end{prf}

\vspace{3mm}

We proceed by providing a converse to Theorem \ref{THM repres. of factorial mom. meas. via Janossy meas.}. Notice that this time we need to impose an additional assumption.
\begin{theorem} \label{THM repres. of Janossy meas. via factorial mom. meas.}
	Let $\eta$ be a point process in $\X$. Let $B \in \mathcal{X}$ and $m \in \N$ be such that
	\begin{equation*}
		\sum_{k = 0}^\infty \frac{\alpha_{\eta, m + k}(B^{m + k})}{k!}
		< \infty.
	\end{equation*}
	Then, for all measurable and bounded functions $f : \X^m \to [0, \infty)$ with $f = 0$ on $\X^m \setminus B^m$,
	\begin{equation*}
		\int_{\X^m} f \, \dd J_{\eta, B, m}
		= \frac{1}{m!} \sum_{k = 0}^\infty \frac{(-1)^k}{k!} \int_{\X^m} f(x_1, \dots, x_m) \, \alpha_{\eta, m + k}\big( \dd(x_1, \dots, x_m) \times B^k \big) .
	\end{equation*}
	Furthermore, for each $B \in \mathcal{X}$ with $\E \big[ 2^{\eta(B)} \big] < \infty$, the Janossy measure of order $0$ is given as
	\begin{equation*}
		J_{\eta, B, 0} = 1 + \sum_{k = 1}^\infty \frac{(-1)^k}{k!} \, \alpha_{\eta, k}(B^k) .
	\end{equation*}
\end{theorem}
\begin{prf}
	First note that the assumption gives $\E\big[ \eta^{(m)}(B^m) \big] = \alpha_{\eta, m}(B^m) < \infty$ which implies $\PP\big( \eta(B) < \infty \big) = 1$, by Proposition \ref{PROP properties of factorial measures due to Last/Penrose}. For $D \in \mathcal{X}^{\otimes m}$ with $D \subset B^m$, we have
	\begin{align} \label{formula Janossy meas. via factorial mom. meas.}
		\sum_{k = 0}^\infty \frac{(-1)^{k}}{k!} \, \alpha_{\eta, m + k}(D \times B^k)
		&= \sum_{k = 0}^\infty \frac{(-1)^{k}}{k!} \sum_{\ell = m + k}^\infty \frac{\ell!}{(\ell - m - k)!} \, J_{\eta, B, \ell}(D \times B^{\ell - m}) \notag \\
		&= \sum_{k = 0}^\infty \frac{(-1)^{k}}{k!} \sum_{\ell = k}^\infty \frac{(\ell + m)!}{(\ell - k)!} \, J_{\eta, B, m + \ell}(D \times B^{\ell}) \notag \\
		&= \sum_{\ell = 0}^\infty (\ell + m)! \, J_{\eta, B, m + \ell}(D \times B^\ell) \sum_{k = 0}^\ell \frac{(-1)^k}{k! \, (\ell - k)!} \notag \\
		&= m! \, J_{\eta, B, m} (D) ,
	\end{align}
	where for the first equality we used Theorem \ref{THM repres. of factorial mom. meas. via Janossy meas.}, for the third equality we applied Fubini's theorem, which is possible since
	\begin{equation*}
		\sum_{k = 0}^\infty \sum_{\ell = k}^\infty \frac{(\ell + m)!}{k! \, (\ell - k)!} \, J_{\eta, B, m + \ell}(D \times B^\ell)
		\leq \sum_{k = 0}^\infty \frac{\alpha_{\eta, m + k}(B^{m + k})}{k!}
		< \infty,
	\end{equation*}
	and the forth equality in \eqref{formula Janossy meas. via factorial mom. meas.} follows from
	\begin{equation*}
		\sum_{k = 0}^\ell \frac{(-1)^k}{k! \, (\ell - k)!}
		= \frac{1}{\ell!} \sum_{k = 0}^\ell {\ell \choose k} (-1)^k
		= \mathds{1}\{ \ell = 0 \} .
	\end{equation*}
	The general result follows from equation \eqref{formula Janossy meas. via factorial mom. meas.} and monotone approximation. Finally, consider the case $m = 0$. As in the assumption, let $B \in \mathcal{X}$ satisfy $\E \big[ 2^{\eta(B)} \big] < \infty$. It follows that $\PP\big( \eta(B) < \infty \big)$ and
	\begin{equation*}
		1 + \sum_{k = 1}^\infty \frac{\alpha_{\eta, k}(B^k)}{k!}
		= \E \Bigg[ 1 + \sum_{k = 1}^\infty \frac{1}{k!} \, \eta^{(k)}(B^k) \Bigg]
		= \E \Bigg[ \sum_{k = 0}^{\eta(B)} {\eta(B) \choose k} \Bigg]
		= \E \big[ 2^{\eta(B)} \big] 
		< \infty.
	\end{equation*}
	In particular, we may use Fubini's theorem to conclude that
	\begin{equation*}
		1 + \sum_{k = 1}^\infty \frac{(-1)^k}{k!} \, \alpha_{\eta, k}(B^k)
		= \E \Bigg[ 1 + \sum_{k = 1}^\infty \frac{(-1)^k}{k!} \, \eta^{(k)}(B^k) \Bigg]
		= \E \Bigg[ \sum_{k = 0}^{\eta(B)} {\eta(B) \choose k} (-1)^k \Bigg]
		= \E \Big[ \mathds{1}\big\{ \eta(B) = 0 \big\} \Big]
		= J_{\eta, B, 0}.
	\end{equation*}
\end{prf}

\begin{corollary} \label{COR repres. of Janossy densities via correlation func.}
	Let $\eta$ be a point process in $\X$. Let $B \in \mathcal{X}$ and $m \in \N$ be such that
	\begin{equation*}
		\sum_{k = 0}^\infty \frac{\alpha_{\eta, m + k}(B^{m + k})}{k!}
		< \infty.
	\end{equation*}
	Assume that the correlation functions $\rho_{\eta, m + k}$ of $\eta$ with respect to $\lambda^{m + k}$ exist for all $k \in \N_0$. Then $J_{\eta, B, m}$ is absolutely continuous with respect to $\lambda^m$ and the corresponding Janossy density is given through
	\begin{equation*}
		j_{\eta, B, m}(x_1, \dots, x_m)
		= \frac{1}{m!} \sum_{k = 0}^\infty \frac{(-1)^k}{k!} \int_{B^k} \rho_{\eta, m + k}(x_1, \dots, x_{m + k}) \, \dd \lambda^k(x_{m+1}, \dots, x_{m + k})
	\end{equation*}
	for $\lambda^m$-a.e.\ $(x_1, \dots, x_m) \in B^m$.
\end{corollary}

The corollary immediately follows from Theorem \ref{THM repres. of Janossy meas. via factorial mom. meas.}. We conclude the discussion with a remark on sufficient conditions that ensure the convergence of the series featuring in Theorem \ref{THM repres. of Janossy meas. via factorial mom. meas.} and Corollary \ref{COR repres. of Janossy densities via correlation func.}.
\begin{remark} \label{RMK convergence of fact. moment meas. series}
	Fix a set $B \in \mathcal{X}$ with $\PP\big( \eta(B) < \infty \big) = 1$. Notice that when considering constructions on the full space $\X$, the bounding constants and maps appearing in this remark may all depend on $B$.
	\begin{enumerate}
		\item The obvious assumption is to require the existence of a constant $c \geq 0$ such that $\alpha_{\eta, \ell}(B^\ell) \leq c^\ell$ for each $\ell \in \N$. Then we have
		\begin{equation*}
			\sum_{k = 0}^\infty \frac{\alpha_{\eta, m + k}(B^{m + k})}{k!}
			\leq c^m \sum_{k = 0}^\infty \frac{c^k}{k!}
			= c^m \cdot e^{c}
			< \infty
		\end{equation*}
		for every $m \in \N$. Similarly, we get $\E \big[ 2^{\eta(B)} \big] \leq e^{c} < \infty$, corresponding to the case $m = 0$.
		
		\item The assumption in part 1 certainly holds if, for each $\ell \in \N$, we have $\alpha_{\eta, \ell}(\cdot) \leq (\vartheta \lambda)^\ell(\cdot)$ for some measurable, $\lambda_B$-integrable map $\vartheta : \X \to [0, \infty)$. If correlation functions exist, the condition is satisfied if
		\begin{equation*}
			\rho_{\eta, \ell}(x_1, \dots, x_\ell)
			\leq \vartheta(x_1) \cdot \dotso \cdot \vartheta(x_\ell)
		\end{equation*}
		for $\lambda^\ell$-a.e.\ $(x_1, \dots, x_\ell) \in B^\ell$, and all $\ell \in \N$. This corresponds to Ruelle's condition in Definition \ref{DEF Ruelle's condition}.
		
		\item For any $\varepsilon > 0$ and $m \in \N$ we find $k_{\varepsilon, m} \in \N$, $k_{\varepsilon, m} > m$, such that for every $k \in \N$, $k \geq k_{\varepsilon, m}$, we have $\big( \tfrac{k}{2} \big)^m \leq \big( 1 + \tfrac{\varepsilon}{2} \big)^k$. This we use to calculate
		\begin{align*}
			\sum_{k = 0}^\infty \frac{\alpha_{\eta, m + k}(B^{m + k})}{k!}
			&= \E \Bigg[ \sum_{k = 0}^\infty \frac{1}{k!} \, \eta(B) \big( \eta(B) - 1 \big) \cdot \dotso \cdot \big( \eta(B) - k - m + 1 \big) \Bigg] \\
			&= \E \Bigg[ \frac{\eta(B)!}{\big( \eta(B) - m \big)!} \, \mathds{1}\big\{ \eta(B) \geq m \big\} \sum_{k = 0}^{\eta(B) - m} {\eta(B) - m \choose k} \Bigg] \\
			&\leq \E \bigg[ \bigg( \frac{\eta(B)}{2} \bigg)^m \cdot 2^{\eta(B)} \, \mathds{1}\big\{ \eta(B) \geq m \big\} \bigg] \\
			&\leq \sum_{k = m}^{k_{\varepsilon, m} - 1} \Big( \frac{k}{2} \Big)^m \cdot 2^k \cdot \PP\big( \eta(B) = k \big) + \sum_{k = k_{\varepsilon, m}}^\infty (2 + \varepsilon)^k \cdot \PP\big( \eta(B) = k \big) \\
			&\leq 2^{k_{\varepsilon, m}} \cdot \bigg( \frac{k_{\varepsilon, m}}{2} \bigg)^m + \E \Big[ (2 + \varepsilon)^{\eta(B)} \Big] .
		\end{align*}
		Thus, if there exists some $b \in (2, \infty)$ such that $\E \big[ b^{\eta(B)} \big] < \infty$, the series converges for every $m \in \N_0$.
		
		\item Assume there exists a map $c : \N \to [0, \infty)$ and a constant $b \in (2, \infty)$ with $\sum_{k = 1}^\infty b^k \cdot c(k) < \infty$, and $J_{\eta, B, k}(B^k) \leq c(k)$ for all $k \in \N$. Then
		\begin{equation*}
			\E \big[ b^{\eta(B)} \big]
			= \sum_{k = 0}^\infty b^k \, \PP\big( \eta(B) = k \big)
			= \sum_{k = 0}^\infty b^k J_{\eta, B, k}(B^k)
			\leq \sum_{k = 1}^\infty b^k \cdot c(k) 
			< \infty.
		\end{equation*}
		By part 3, the assumption of Theorem \ref{THM repres. of Janossy meas. via factorial mom. meas.} is satisfied for every $m \in \N_0$.
		
		\item In case Janossy densities exist, the previous condition is satisfied if
		\begin{equation*}
			j_{\eta, B, k}(x_1, \dots, x_k)
			\leq \frac{\vartheta(x_1) \cdot \dotso \cdot \vartheta(x_k)}{k!}
		\end{equation*}
		for $\lambda^k$-a.e.\ $(x_1, \dots, x_k) \in B^k$, all $k \in \N$, and some $\lambda_B$-integrable function $\vartheta : \X \to [0, \infty)$.
	\end{enumerate}
\end{remark}

\section{A Kolmogorov extension result for measures on $\NN$}
\label{Appendix Kolmog. extension result on N}

The results in this part of the appendix are presented in full detail in the lecture notes by \cite{P:2005:meas, P:2005:specif}. Though these are not peer reviewed publications in scientific journals, the math is sound and we have checked all results we present here as well as the preliminaries leading up to them, and, aside from a few typos, they are perfectly correct. We refrain from giving proofs here and only cite and sketch the essential results.

First, consider the space $\mathbb{S} = \{ 0, 1 \}^\infty$ of all sequences of $0$'s and $1$'s. Endowed with the metric $d : \mathbb{S} \times \mathbb{S} \to [0, 1]$,
\begin{equation*}
	d\big( (s_n)_{n \in \N}, \, (s^\prime_n)_{n \in \N} \big)
	= \sum_{n = 1}^\infty \frac{|s_n - s^\prime_n|}{2^n} ,
\end{equation*}
the space $\mathbb{S}$ is a compact metric space and we denote by $\mathcal{B}(\mathbb{S})$ its Borel $\sigma$-field. It is easy to verify that the $\sigma$-field of a measurable space $(\X, \mathcal{X})$ is countably generated if, and only if, there exists a map $f : \X \to \mathbb{S}$ such that $f^{-1}\big( \mathcal{B}(\mathbb{S}) \big) = \mathcal{X}$. 

A measurable space $(\X, \mathcal{X})$ is called a \textit{substandard Borel space} if there exists a map $f : \X \to \mathbb{S}$ such that $f^{-1}\big( \mathcal{B}(\mathbb{S}) \big) = \mathcal{X}$ and $f(\X) \in \mathcal{B}(\mathbb{S})$. Thus, substandard Borel spaces are countably generated measurable spaces that satisfy some additional property. In particular, Borel spaces, and thus also complete separable metric spaces, are substandard Borel, which can be shown with standard constructions in the lines of Theorem 6.1 of \cite{P:2007:borel}. Substandard Borel spaces provide enough structure to supply existence results for probability kernels and extension results as the ones we state below. Note that \cite{P:2007:borel} discusses these spaces in much detail and calls them type-$\mathcal{B}$ spaces.

The results culminating in Proposition 18.2 (4) of \cite{P:2005:meas} and Lemma 2.2 of \cite{P:2005:specif} show that if $(\X, \mathcal{X})$ is a substandard Borel space, then $(\NN, \mathcal{N}_B)$ is also substandard Borel for each $B \in \mathcal{X}_b$. Proposition 2.3 of \cite{P:2005:specif}, which is proven via the general Kolmogorov extension result presented in Theorem 19.1 of \cite{P:2005:meas}, shows that $(\NN, \mathcal{N})$ is substandard Borel and an extension theorem holds for probability laws on the space of counting measures. The precise result reads as follows.

\begin{prop} \label{PROP Kolmogorov extension theorem, countable version}
	Let $(\X, \mathcal{X})$ be a localized substandard Borel space. For each $\ell \in \N$ let $\mathrm{P}_\ell$ be a probability measure on $(\NN, \mathcal{N}_{B_\ell})$ such that $\mathrm{P}_\ell(A) = \mathrm{P}_i(A)$ for all $A \in \mathcal{N}_{B_i}$ whenever $i < \ell$ (Kolmogorov consistency property). Then there exists a unique probability measure $\mathrm{P}$ on $(\NN, \mathcal{N})$ such that
	\begin{equation*}
		\mathrm{P}(A)
		= \mathrm{P}_\ell(A), \quad A \in \mathcal{N}_{B_\ell} , ~ \ell \in \N .
	\end{equation*}
\end{prop}

The reader who is not fully convinced by the reference to the lecture notes can proceed along the following lines to obtain Proposition \ref{PROP Kolmogorov extension theorem, countable version} in the case where $(\X, \mathcal{X})$ is a Borel space. For each $D \in \mathcal{X}$ the space $(D, \mathcal{X} \cap D)$ is a Borel space and so is the space $\big( \NN(D), \mathcal{N}(D) \big)$ by Theorem 1.5 of \cite{K:2017}. The map $q_D : \NN \to \NN(D)$, defined by $q_D(\mu)(B) = \mu(B)$ (for $B \in \mathcal{X} \cap D$), is surjective and $q_D^{-1}\big( \mathcal{N}(D) \big) = \mathcal{N}_{D}$, so $(\NN, \mathcal{N}_D)$ is a substandard Borel space. Next, one has to observe that the proof of Theorem 4.1 in Chapter V of \cite{P:1967} (omitting the last part of that theorem as we already know that $(\NN, \mathcal{N})$ is a Borel space) actually only requires the spaces $(\NN, \mathcal{N}_D)$ to be substandard Borel. Thus, it remains for the reader to convince himself that the atom-condition of Theorem 4.1 in Chapter V of \cite{P:1967} is satisfied. To give a sketch, assume that $A_1 \supset A_2 \supset \dotso$ are subsets of $\NN$ such that $A_n$ is an atom of $(\NN, \mathcal{N}_{B_n})$ \citep[in the sense of Definition 2.1 in Chapter V of][]{P:1967} for each $n \in \N$. Then, using the easy-to-prove fact that the intersection of two atoms from (potentially) different countably generated $\sigma$-fields over the same base set is either empty or an atom of the $\sigma$-field generated by the union of the two $\sigma$-fields, it is straight forward (if a bit technical) to explicitly construct a sequence of measures $(\mu_n)_{n \in \N}$ in $\NN$ such that $\mu_n \in A_n$ with $\mu_n(\X) = \mu_n(B_n)$ and $(\mu_{n + 1})_{B_n} = \mu_n$ for every $n \in \N$. Hence there exists a measure $\mu \in \NN$ so that $\mu_{B_n} = \mu_n$ for each $n \in \N$ and so $\mu \in \bigcap_{n = 1}^\infty A_n \neq \varnothing$. For the full construction of the sequence $(\mu_n)_{n \in \N}$ we refer to Lemma 2.4 of \cite{P:2005:specif}, where all considerations are laid out. With this sketch, referring to the lecture notes only as a convenient check up and not as a necessary reference, the reader can obtain Proposition \ref{PROP Kolmogorov extension theorem, countable version} for Borel spaces.

\section*{Acknowledgements}

The author is grateful to G\"unter Last for valuable discussions on every version of this manuscript. The author is also thankful to Jan Philipp Neumann who read an earlier version of this manuscript and provided countless detailed suggestions for improvement. The results of this paper stem from the author's Ph.D.\ research, \cite{B:2022}.




\bibliography{gibbsreferences}   
\bibliographystyle{apalike}

\end{document}